\documentclass[10pt]{article}
\usepackage{amsmath,amssymb,amsfonts,amsthm,amscd,graphicx,psfrag,epsfig, arcs}
\usepackage{color}
\definecolor{blue}{rgb}{0,0,0.7}
\definecolor{red}{rgb}{0.75, 0, 0}
\usepackage{stmaryrd}
\usepackage[title]{appendix}
\usepackage{hyperref}
\hypersetup{colorlinks, linkcolor=blue, citecolor=cyan}
\usepackage[all]{xy}           
\usepackage{marginnote}  
\usepackage{float}
\usepackage{soul}

\newtheorem{theorem}{Theorem}[section]
\newtheorem{lemma}[theorem]{Lemma}
\newtheorem{proposition}[theorem]{Proposition}
\newtheorem{corollary}[theorem]{Corollary}
\newtheorem{conjecture}[theorem]{Conjecture}
\newtheorem{definition}[theorem]{Definition}

\newtheorem{remark}[theorem]{Remark}

\newcommand{\MeijerG}[7]{G \begin{smallmatrix} #1 & #2 \\ #3 & #4 \end{smallmatrix} \left( \begin{smallmatrix} #5 \\ #6 \end{smallmatrix} \middle\vert #7 \right) }

\newcommand{\bpf}{\begin{proof}}
\newcommand{\epf}{\end{proof}}
\newcommand{\bs}{\begin{split}}
\newcommand{\es}{\end{split}}
\newcommand{\be}{\begin{equation}}
\newcommand{\ee}{\end{equation}}
\newcommand{\bt}{\begin{theorem}}
\newcommand{\et}{\end{theorem}}
\newcommand{\bd}{\begin{definition}}
\newcommand{\ed}{\end{definition}}
\newcommand{\bp}{\begin{proposition}}
\newcommand{\ep}{\end{proposition}}
\newcommand{\bl}{\begin{lemma}}
\newcommand{\el}{\end{lemma}}
\newcommand{\bc}{\begin{corollary}}
\newcommand{\ec}{\end{corollary}}
\newcommand{\bcon}{\begin{conjecture}}
\newcommand{\econ}{\end{conjecture}}
\newcommand{\la}{\label}

\newcommand{\Z}{{\mathbb Z}}
\newcommand{\R}{{\mathbb R}}

\newcommand{\C}{{\mathbb C}}
\newcommand{\G}{{\rm G}}

\newcommand{\lra}{\longrightarrow}
\newcommand{\lms}{\longmapsto}
\newcommand{\bS}{{\Bbb S}}

\newcommand{\bg}{\begin{equation}\begin{gathered}}
\newcommand{\eg}{\end{gathered}\end{equation}}

\newcommand{\F}{{\rm F}}

\newcommand{\Y}{{\rm Y}}

\newcommand{\K}{{\rm K}}

\newcommand{\N}{{\rm N}}

\newcommand{\X}{{\rm X}}

\newcommand{\old}[1]{}

\newcommand{\E}{{\mathcal E}}
\newcommand{\D}{{\rm  D}}
\newcommand{\rM}{{\rm  M}}

\newcommand{\T}{{\rm T}}
\newcommand{\V}{{\rm V}}
\newcommand{\A}{{\rm A}}
\newcommand{\B}{{\rm B}}
\newcommand{\U}{{\rm {U}}}
\renewcommand{\P}{\mathbb{P}}

\renewcommand{\H}{{\rm H}}
\renewcommand{\P}{{\rm P}}

\newcommand{\bspl}{\begin{split}}
\newcommand{\espl}{\end{split}}
\newcommand{\epr}{\end{proof}}
\newcommand{\bpr}{\begin{proof}}

\renewcommand{\L}{{\rm L}}

\input xy
\xyoption{all}

\begin{document}

\date{}

\title {Mellin transformed exponential volumes and probabilistic moduli spaces}
\author{Alexander B. Goncharov, Zhe Sun}

\maketitle

 \tableofcontents

\abstract

A decorated surface $\bS$ is a surface with boundary and a finite set of marked boundary points up to  isotopy,   called {cusps}. 
The intervals on the boundary between the cusps are called boundary intervals. 
\vskip 1mm

We continue to study the {\it exponential volumes} ${\cal E}_\bS(\K_\bullet,\L_\bullet)$ of moduli spaces of ideal hyperbolic structures on a decorated surface $\bS$, with the given values of the exponents of   lengths of 
 boundary geodesic circles $\L_\bullet$ and boundary geodesic intervals $\K_\bullet$, introduced in our previous work. 

  \vskip 2mm
  1. We prove that the Mellin--Laplace transform of the  function ${\cal E}_\bS(\K_\bullet,\L_\bullet)$ admits {\em Barnes integral presentations}, 
  and satisfies  recursions,  generalising Mirzakhani's recursions.%
  \vskip 2mm
  
2. We specify an arbitrary subset ${\cal I}$ of the boundary intervals, called blue intervals, and assign to them any collection of weights $\omega_e\in \C$, $e\in {\cal I}$. 
We show that the Mellin transform of the function ${\cal E}_\bS(\K_\bullet,\L_\bullet)$ along the variables $\K_e$ at the blue intervals, evaluated at the  $\omega_e$, 
calculates  {\it weighted exponential  volumes} of various moduli spaces. 
If all boundary intervals are blue and $\omega_e=\frac{1}{2}$, we recover Chekhov's volumes of moduli spaces of hyperbolic surfaces with crowns, discretizing the Schwarzian action. 
We derive recursions for the weighted exponential volumes of such moduli spaces. 
\vskip 2mm

 3. We introduce {\it probabilistic Teichm\"uller spaces} related to decorated surfaces $\bS$. They are infinite dimensional and carry canonical Wiener type measures.
 The definition of the measures builds on the works of Bauerschmidt - Losev - Wildemann and Losev. 
 Their quotients by the action of the  modular group of $\bS$ are  {\it probabilistic moduli spaces}. We show that their volumes  are finite. 
 We prove recursion formulas for these measures, and hence for the volumes. 

 The emerging theory 
  generalizes the Schwarzian field theory, studied by Stanford - Witten and others, from the disc 
   to arbitrary decorated surfaces, i.e. to  open string format. In particular,  we get recursions for the Feynman integrals in   Schwarzian field theory on decorated surfaces.    
   
    \section{Introduction}

\subsection{Mellin--Laplace transform of  exponential volumes via Barnes integrals}

\subsubsection{Exponetial volumes of moduli spaces of ideal hyperbolic structures \cite{GS24}}

A {\it decorated surface $\bS$} is an oriented topological surface with boundary circles, and a finite number of marked points on the boundary, considered modulo isotopy. 
 
The {\it punctured boundary} $\widehat \partial \bS$ is defined as   the boundary $\partial \bS$ of $\bS$ minus the marked boundary points:
\be \la{PB}
\widehat \partial \bS:=  \partial \bS- \{\mbox{marked boundary points}\}.
\ee
 It is a union of {\it boundary circles} and {\it boundary intervals}, called  {\it punctured boundary components of $\bS$}. 
We refer to   marked points as {\it cusps}, and to boundary circles without cusps as {\it holes} or {\it punctures}.    
\vskip 1mm
An {\em ideal hyperbolic structure} on a decorated surface $\bS$ is a Riemannian metric of the curvature $-1$ on $\bS$ with the following boundary structure:  

\begin{enumerate}

\item Each boundary circle and boundary interval is a geodesic.

\item The hyperbolic metric near each  
marked  point $p$ is a cusp..   

\item We choose a horocycle $h_p$ for each cusp $p$. 
 \end{enumerate}
 
The {\em Teichm\"uller space} ${\cal T}_{{\bS}}$ parametrises ideal hyperbolic structures on $\bS$ up to orientation and boundary preserving homeomorphisms of $\bS$ isotopic to the identity. 

The {\em moduli space} $\mathcal{M}_{\bS}$ is the  quotient of the {space} ${\cal T}_{{\bS}}$ by the pure mapping class group ${\rm Mod}(\bS)$:\footnote{The {pure mapping class group} ${\rm Mod}(\bS)$  is the quotient of the group of orientation $\&$ boundary preserving homeomorphisms of $\bS$ modulo the ones   isotopic to the  identity. }
$$
\mathcal{M}_{\bS}:= {\cal T}_{{\bS}}/{\rm Mod}(\bS).
$$

There are the following functions on the moduli space ${\cal M}_\bS$:

\begin{itemize}

\item 
For each cusp  $p$ on $\bS$, the {\it local potential} $W_p$, given by the length of the horocycle $h_p$ at  $p$.

 The sum of the  {local potentials} $W_p$ over all cups $p$ of $\bS$ is called the {\it potential}:
$$
W_\bS:= \sum_p W_p.
$$
The potentials $W_p$ and  $W_\bS$ are ${\rm Mod}(\bS)-$invariant. So they are functions on the  space $\mathcal{M}_\bS$.

 \item For each boundary interval
 $e$ on $\bS$ with the vertices at the cusps $p, q$, let $\kappa_e$ be the  distance between the horoarcs $h_p$ and $h_q$,  taken with the minus sign if they overlap.   Set  
 \be \la{Ke}
 \K_e:= e^{-\kappa_e}.
 \ee

\item For each hole $h$ on $\bS$, we denote by $l_{h}$ its geodesic boundary, and  set 
$$
\L_{h}:= e^{l_{h}}.
$$

\end{itemize} 

Denote by ${\cal M}_\bS(\K, \L)$ the subspace of the moduli space ${\cal M}_\bS$  with the fixed collections of parameters $\K = \{\K_e\}$ and $\L=\{\L_{h}\}$. Any volume form on
 the space ${\cal M}_\bS$ induces a volume form on the  ${\cal M}_\bS(\K, \L)$. 
 \vskip 1mm
 
 Given an orientation of the Teichm\"uller space $\mathcal{T}_\bS$, the space $\mathcal{T}_\bS$ carries  canonical  
 volume form $\Omega_\bS$. 
 It has a simple expression in the cluster Poisson coordinates. 
It is ${\rm Mod}(\bS)-$invariant,  and thus induces a volume form $\Omega_\bS$ on the moduli space ${\cal M}_\bS$, which
  induces a volume form $\Omega_\bS({{\K}, \L})$ on   $\mathcal{M}_\bS({\K}, \L)$.
When $\bS$ has no marked points, it is propotional to the  Weil--Petersson volume form. 
\vskip 1mm

The space  ${\cal M}_\bS$ carries the {\it exponential volume form} \cite{GS24}: 
$$
{\Bbb E}_\bS:= e^{-W_\bS}\Omega_\bS.
$$
It induces the  {\em exponential volume form} on the moduli space ${\cal M}_\bS(\K, \L)$:
\be \la{6}
\mathbb{E}_{\bS}(\K, \L):=e^{-W_\bS}\Omega_\bS({{\K}, \L}).
\ee
 The  {\em exponential  volume}  
of the moduli space $\mathcal{M}_\bS({\K}, \L)$ is its integral 
over  ${\cal M}_\bS(\K, \L)$:
\be \la{EXV}
\mathcal{E}_\bS({\K}, \L):= \int_{\mathcal{M}_\bS({\K}, \L)} e^{-W_\bS}\Omega_\bS({{\K}, \L}). 
\ee
By  \cite[Theorem 1.5]{GS24}  the integral is finite. It provides  a function on the  variables $\K, \L$. 

\subsubsection{The Mellin--Laplace transform of the exponential volumes}

Recall the {\em Mellin and Laplace transforms} of a function $f$: 
\be
\begin{split}
&\mathcal{M}(f,s)=\int_{0}^{+\infty} f(\K)\K^{s} d \log \K, \\
&\mathcal{L}(f,z)=\int_{0}^{+\infty} f(l)e^{-lz} dl.\\
\end{split}
\ee
 
 \bd 
\label{defn:LM}
The  {Mellin--Laplace transform  of the exponential volume function}  
 $ {\cal E}_\bS(\K, \L)$
\be \la{LMTR}
 \mathcal{ML}({\cal E}_\bS; s_1, ..., s_p; z_1, ..., z_q):= \int {\cal E}_\bS(\K, \L)  \bigwedge_{i=1}^p \K^{s_i}  d\log \K_i\wedge \bigwedge^q_{j=1}e^{-l_j z_j} dl_j
 \ee
 is  the Mellin transform along  each 
boundary interval variable $\K_e$, and the Laplace transform along each boundary geodesic  length $l_{h}\geq 0$. It is a function of the dual parameters $s_1, ..., s_p$ and $z_1, ..., z_q$, where $p$ is the  number of the boundary intervals, and $q$  is the number of holes.
\ed

  We  refer to $s_1, ..., s_p$ as the {\it Mellin variables}, and to $z_1, ..., z_q$ as the {\it Laplace variables}. 
  \vskip 1mm

  Our first main result  is the following theorem,   which is  stated precisely in Theorem \ref{thmint}.  
  \bt  \la{MTHML}
  For any decorated surface $\bS$, the Mellin--Laplace transform function (\ref{LMTR}) has a {\em Barnes integral presentation}. 
  \et

 Here is a sample result, see Section \ref{Sec3.2} for details.
  
\begin{figure}[ht]
	\centering
	\includegraphics[scale=0.17]{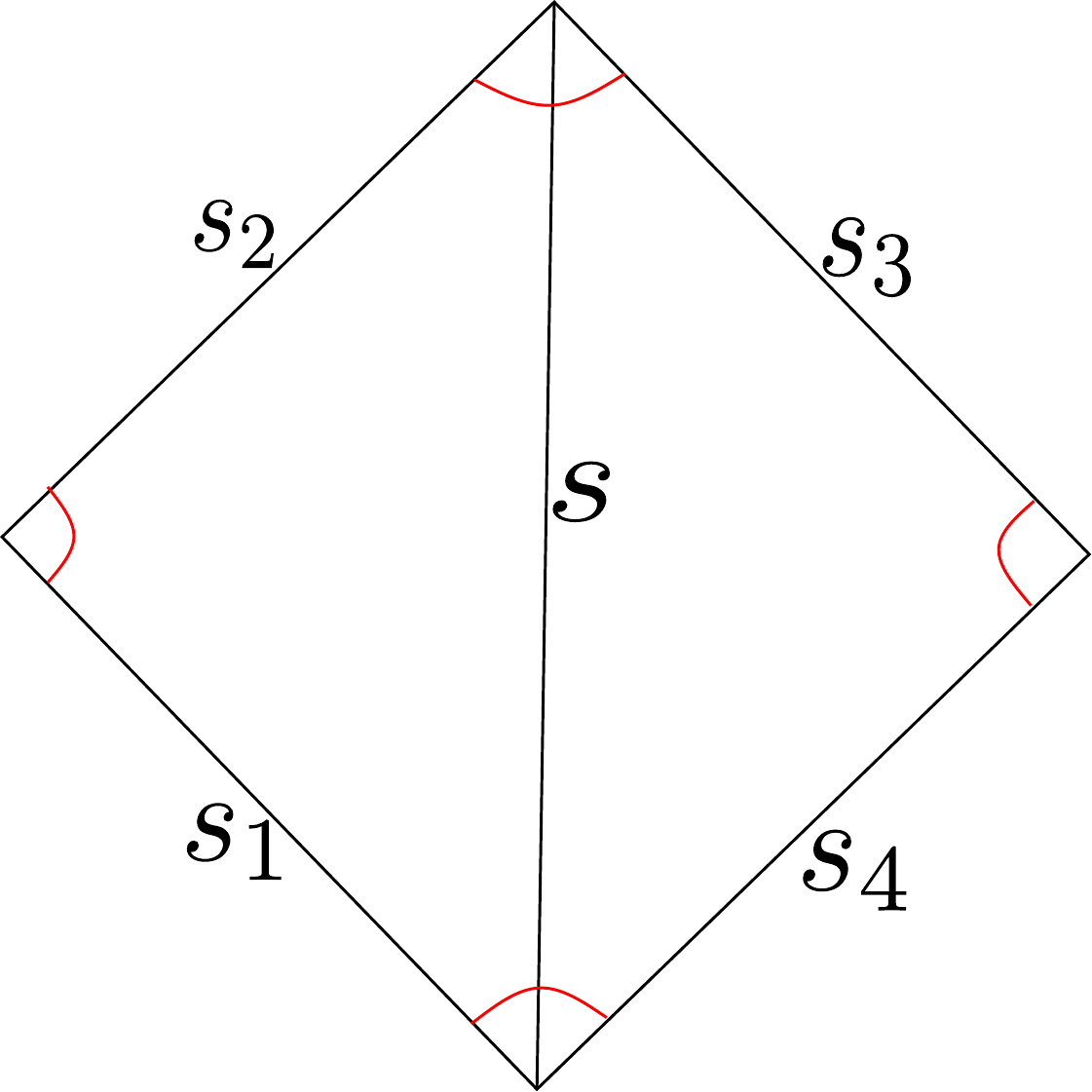}
	\small
	\caption{The  quadrilateral with the sides $\&$  diagonal labeled by the Mellin variables $s_1, ..., s_4, s$.}
	\label{figure:D4s}
\end{figure}

  \bp Let $\bS=\square$ be a quadrilateral as on Figure \ref{figure:D4s}. Then the Mellin transform 
  of the exponential volume function $\mathcal{M}(\mathcal{E}_{\square}; s_1,s_2,s_3,s_4)$ is  calculated by a Barnes $\Gamma-$integral:
\be
\label{eq:334x}
\begin{aligned}
\mathcal{M}(\mathcal{E}_{\square}; s_1,s_2,s_3,s_4) 
&=2\Gamma(s_1+s_2)\Gamma(s_3+s_4) \int_{-i \infty}^{+i \infty}  \Gamma(s_1+s)  \Gamma(s_2+s)  \Gamma(s_3-s)  \Gamma(s_4-s) ds 
\\&\overset{\rm Barnes ~ Lemma}{=}\frac{2\prod_{1\leq i<j\leq 4}\Gamma(s_i+s_j)}{\Gamma(s_1+s_2+s_3+s_4)}.\end{aligned}
\ee
The  integration  contour is a vertical line, denoted by $(-i\infty, i\infty)$,  which separates the poles of $\Gamma(a+s)$ and $\Gamma(b+s)$ from those of $\Gamma(c-s)$ and $\Gamma(d-s)$, 
 and oriented so that the points $a, b$ are on the left.
 \ep
 
 The first equality shows how the general Barnes integral for the exponential volumes looks like. The second equality in (\ref{eq:334x}) is the Barnes first Lemma, which  is specific for the quadrilateral. 
  
\vskip 1mm

 By \cite[Theorem 1.11]{GS24}   exponential volumes satisfy recursion formulas, generalising Mirzakhani's recursions. 
  We prove  that the Mellin--Laplace transformed exponential volumes satisfy similar recursions, given by   Theorem \ref{thm:girre}.

\subsubsection{Mellin--Laplace transform for $u-$deformed exponential volumes}

We can consider more general potentials. Let $u=\{u_p\}$ be a collection of positive real numbers assigned to the cusps $p$ of $\bS$. Take a linear combination of the local potentials $W_p$ weighted by $u_p$:
$$
u\cdot W_\bS:= \sum_p u_p W_p.
$$
We define the {\it $u-$deformed exponential volume} by
 \be \la{EVF1}
 {\cal E}_{\bS, u}(\K, \L ):= \int_{{\cal M}_\bS(\K, \L)} e^{-u\cdot W_\bS}\Omega_\bS(\K, \L). 
 \ee 
 The following result shows that the  Mellin--Laplace transform of the $u-$deformed exponential volume is easily reduced to the initial one. Precisely,  denote by $s'_p$ and $s''_p$ the Mellin 
 $s-$variables assigned to the boundary edges sharing the cusp $p$. Following Remark \ref{rem:deform}, we have
\bl \la{tLMEVL}
The Mellin--Laplace transform of the $u-$deformed exponential volume is calculated by 
\be \la{tLMEV}
\mathcal{ML}({\cal E}_{\bS, u}; \K, \L ) = \prod_{p} u_p^{s'_p + s''_p}  \cdot    \mathcal{ML}({\cal E}_\bS; \K, \L ).
\ee
\el
 Lemma \ref{tLMEVL} is proved in the end of Section \ref{SEC4.1}.
 
 Combining Lemma \ref{tLMEVL} with Theorem \ref{MTHML}, and performing the inverse Mallin-Laplace transform, we can get asymptotics of the $u-$deformed exponential volumes as $u_p \to \infty$. 
  

 \subsection{Weighted exponential volumes via Mellin transform}
 
 We consider the Mellin transform  
 of the exponential volume ${\cal E}_\bS(\K_\bullet, \L_\bullet)$ along some of the $\K-$variables, and specialize it at the given values of Mellin parameters. 
 To record this data, we introduce {\it weighted surfaces}. 
 
 \bd A {weighted surface} is a triple $\bS_\omega= (\bS, {\cal I}_{\omega}, \omega)$,  given by a decorated surface $\bS$, a  subset 
  ${\cal I}_{\omega}$ of the set of  boundary intervals ${\cal I}_\bS$, and a weight function 
  \be \la{wt}
\omega: {\cal I}_{\omega} \lra \C.
\ee
 \ed
 
 We assign to a weighted surface $(\bS, {\cal I}_{\omega}, \omega)$ the Mellin transform of the exponential volume   function ${\cal E}_\bS(\K, \L)$ 
 for the $\K-$variables $\K_e$, where $e \in {\cal I}_\omega$, specialized at the Mellin variables $s_e:=\omega_e$, provided by the weights. 
 We keep intact the rest of the parameters: the $\K-$variables $\K_f$ at the boundary intervals $f \in {\cal I}_\bS-{\cal I}_\omega$ and $\L$. We denote the obtained function by 
  \be \la{PMTR}
  \mathcal{M}({\cal E}_\bS; s_{\bf e}, \K_{\bf f}, \L), \ \ \  s_{\bf e}:= \{s_e\}, e \in {\cal I}_\omega, \ \ \  \K_{\bf f}:=\{\K_f\}, f \in {\cal I}_\bS-{\cal I}_\omega.
  \ee
  
  Our key point is that each function (\ref{PMTR}) can be interpreted as the {\it weighted exponential volume} of certain moduli space, see Sections \ref{sec1.2.1}-\ref{sec1.2.2}. In fact usually there are more than one, although finitely many,  
  of  such moduli spaces. To describe these moduli spaces,  and to define their wighted exponential volumes, we introduce {\it colored surfaces} in Definition \ref{DEF2.7}. 
  
  \begin{figure}[ht]
	\centering
	\includegraphics[scale=0.18]{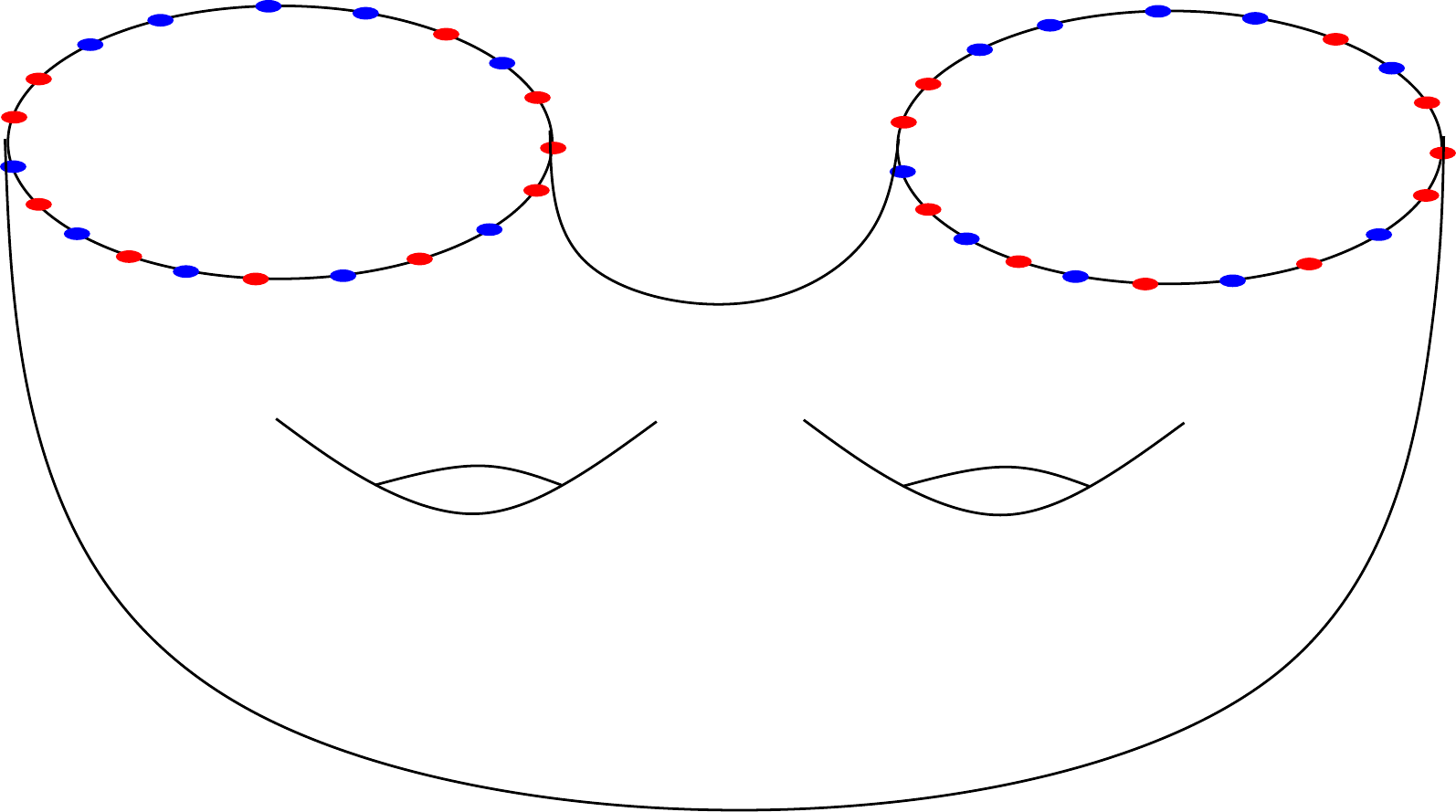}
	\small
	\caption{A colored surface.}
	\label{figure:mixa}
\end{figure}

 
Given a decorated surface $\bS$, let us denote by ${\cal V}_\bS$ the set of its marked boundary point (the cusps).

\subsubsection{Colored surfaces and related moduli spaces} \la{sec1.2.1}
\begin{definition} \la{DEF2.7}
A {\em colored surface} $\bS_{\rm{c}}=(\bS; {\cal V}_{\bf r})$ is a decorated surface $\bS$ with a subset   ${\cal V}_{\bf r}\subset {\cal V}_\bS$ of its red cusps. 
\end{definition}

We color by blue  the other cusps 
on 
 $\bS$, see Figure \ref{figure:mixa}, and denote by ${\cal V}_{\bf b}$ the set of  blue cusps. 

\begin{figure}[ht]
	\centering
	\includegraphics[scale=0.14]{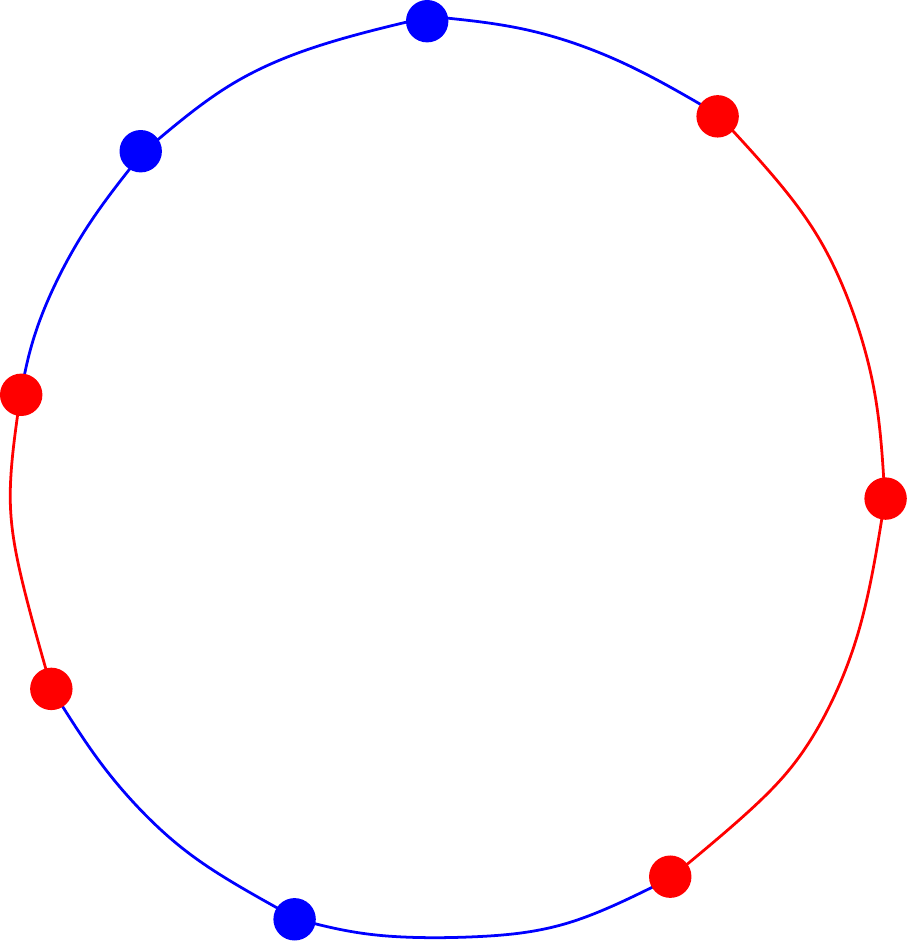}
	\small
	\caption{Corloring boundary edges on a colored surface.}
	\label{figure:coloredgea}
\end{figure}

 We color a boundary interval by red if and only if both of its vertices are  red, and denote by ${\cal I}_{\bf r}\subset {\cal I}_\bS$ the subset of all red boundary intervals, see Figure \ref{figure:coloredgea}. We color the rest of the boundary intervals by blue, and denote by ${\cal I}_{\bf b}\subset {\cal I}_\bS$ the set of all blue boundary intervals. 
  So  we get  partitions
\be \la{VD}
\begin{split}
&{\cal V}_\bS= {\cal V}_{\bf r} \cup {\cal V}_{\bf b}.\\
&{\cal I}_\bS= {\cal I}_{\bf r} \cup {\cal I}_{\bf b}.\\
\end{split}
\ee 
We note that not all partitions of the set of boundary intervals into red and blue ones come this way.

\bd
An {\em ideal hyperbolic structure} on a colored surface $\bS_{\rm c}$  is obtained  from an ideal hyperbolic structure on the underlying decorated surface $\bS$ by 
 forgetting  horocycles   at the blue cusps. \ed
 
The {\em Teichm\"uller space} ${\cal T}_{{\bS}_c}$ parametrises ideal hyperbolic structures on $\bS_{\rm{c}}$ up to orientation and boundary preserving homeomorphisms of $\bS$ isotopic to the identity. 

The {\em moduli space} $\mathcal{M}_{\bS_{\rm{c}}}$ is the  quotient of the {space} ${\cal T}_{{\bS}_c}$ by the pure mapping class group ${\rm Mod}(\bS)$. 
\vskip 1mm

 Given a colored surface $\bS_{\rm{c}}$, there are the following functions  on the Teichm\"uller space ${\cal T}_{\bS_{\rm{c}}}$:
 
 \begin{itemize}
 
 \item The local potentials $W_j$ at the {\it red cusps}.  The potential  given by 
 the sum of  the local potentials \be
 W_{\bS_{\rm{c}}}:= \sum_{j \in {\cal V}_{\bf r}} W_j.
 \ee

\item The $\K-$invariant $\K_e$ at each red boundary interval ${\rm I}_e$, $e \in {\cal I}_{\bf r}$, see (\ref{Ke}).

\end{itemize}

The local potentials $W_j$,  the potential $W_{\bS_{\rm{c}}}$, and the functions $\K_e$  are ${\rm Mod}(\bS)$-invariant. So they descends to functions on the moduli space $\mathcal{M}_{\bS_{\rm{c}}}$.

\subsubsection{Weighted exponential volumes} \la{sec1.2.2}

\begin{definition} \la{D1.10} A {\it weighted colored surface}  $\bS_{\rm{c}; \omega} = (\bS_{\rm c}; \omega) = (\bS_{\rm c}; {\cal V}_{\bf r}, \omega) $  is a colored surface $\bS_{\rm{c}}$ 
equipped with   a   function on the set of blue  intervals, called the   {\em weight $\omega$}:
$$
\omega: {\cal I}_{\bf b} \lra \C.
$$ 
\ed
 
A {\it half-interval} on $\bS$  is a pair $(v,e)$ where $v\in {\cal V}_\bS$   and $e\in {\cal I}_\bS$. 
There is a  map $\kappa: (v,e) \lms e$ from the set of half-intervals on $\bS$ to the set  ${\cal I}_\bS$. So a function  $\omega$ on the set  ${\cal I}_{\bf b}$ of blue intervals  induces a function $\omega':= \kappa^*\omega$  on the  set of blue half-intervals,  and hence 
a function $t$ on the set of blue cusps:
\be \la{v}
 t_b= \omega'_e+\omega'_{e'}= \omega_e+\omega_{e'}, \ \ \ \forall b\in {\cal V}_{\bf b},
\ee
where  $e$ and $e'$ are the half-intervals sharing  the cusp $b$.

\bd  \la{LD.11} The {\em weight factor} $\mathcal{Q}_{{\bS}_{\rm c}; \omega}$ is a function on the Teichm\"uller space ${\cal T}_{\bS_{\rm{c}}}$  is given by
\be \la{Gx}
\mathcal{Q}_{\bS_{\rm{c}}}(\omega):=\prod_{b \in {\cal V}_{\bf b}} \Gamma(t_b)\frac{\prod_{e\in \mathcal{I}_{\bf b}}\K_{e}^{\omega_e}}{ \prod_{b\in \mathcal{V}_{\bf b}} W_b^{t_b}}.
\ee

The {weighted exponential volume form} on the Teichm\"uller space ${\cal T}_{\bS_{\rm{c}}}$  is
$$
{\Bbb E}_{\bS_{\rm{c}}}(\omega):= {\cal Q}_{\bS_{\rm{c}}}(\omega) e^{-W_{\bS_{\rm{c}}}}\Omega_{\bS_{\rm{c}}}.
$$ 
\end{definition}

 \bl \la{L5.1a} i) The  weight factor ${\cal Q}_{\bS_{{\rm{c}}}}(\omega)$ does not depend on the  horocycles $h_b$ at the blue cusps $b$ on $\bS_{\rm{c}}$. 
 So it is a well-defined function on the Teichm\"uller space ${\cal T}_{\bS_{\rm{c}}}$.

ii) The weighted exponential volume form  on the  space ${\cal T}_{\bS_{\rm{c}}}$  descends to the moduli space ${\cal M}_{\bS_{\rm{c}}}$. 
\el

Lemma \ref{L5.1a} is proved  in Proposition \ref{L5.1}.
 
 The weight factor is a function of the variables $t_b$ in (\ref{v}), and  has a pole  if  and only if 
$
 t_b \in \Z_{\leq 0}.
$ 
  \vskip 2mm
  \bd \la{D1.14} Given a weighted colored surface  $(\bS_{\rm{c}};\omega)$, let $\K_{\bf r}$ be the $\K$-parameters at the  red boundary intervals. The weighted exponential volume of the moduli space ${\cal M}_{\bS_{\rm{c}}}$ 
is defined by 
\be
{\cal E}_{\bS_{{\rm{c}}}}(\omega; \K_{\bf r}, \L):=\int_{{\cal M}_{\bS_{\rm{c}}}(\K_{\bf r}, \L)}{\Bbb E}_{\bS_{\rm{c}}}(\omega; \K_{\bf r}, \L). 
\ee
\ed
In Section \ref{SEC10} we proved recursions for the weighted exponential volumes for colored surfaces with at least one red cusp.

\subsubsection{The weighted exponential volumes  via the Mellin transform}

Given a weighted colored surface  $(\bS_{{\rm{c}}};\omega)$, denote by $\mathcal{M}_{\rm blue}(\cdot)$ the Mellin transform of the exponential volume function along the $\K-$variables at the blue edges. The resulting function is 
$$
\mathcal{M}_{\rm blue}(\mathcal{E}_{\bS}(\K_{\bf r}, \L); s_1, ..., s_b).
$$
It depends on the Mellin parameters $s_e$, $e \in {\cal I}_{\bf b}$  at the blue edges of ${\Bbb S}_{\rm c}$, on the $\K-$parameters $\K_f$, $f\in {\cal I}_{\bf r}$ at the red edges, and on  the set $\L$ of the exponents of geodesic lengths  of the neck geodesics. 
Specializing this function to the values of the Mellin parameters given by the weights $s_e=\omega_e$, we get a function 
$$
\mathcal{M}_{\rm blue}(\mathcal{E}_{\bS}(\K_{\bf r}, \L); \omega_1, ..., \omega_b).
$$
On the other hand, the weighted exponential volume of the moduli space ${\cal M}_{\bS_{\rm c}}$ from Definition \ref{D1.14} provides another function depending on the same variables: 
$$
\mathcal{E}_{{\bS}_{\rm c}}(\omega_1, ..., \omega_b, \K_{\bf r}, \L).
$$

\bt
\label{thm:fc}
For any weighted colored surface  $(\bS_{{\rm{c}}};\omega)$, 
 the weighted exponential volume of the moduli space ${\cal M}_{\bS_{\rm c}}$ coincides with the 
 Mellin transform of the exponential volume of the moduli space ${\cal M}_\bS$ for the underlying decorated surface $\bS$, evaluated at the weights $ \omega_1, ..., \omega_b$:
  \be
\label{eq:fc}
\mathcal{E}_{{\bS}_{\rm c}}(\omega_1, ..., \omega_b, \K_{\bf r}, \L)=\mathcal{M}_{\rm blue}(\mathcal{E}_{\bS}(\K_{\bf r}, \L); \omega_1, ..., \omega_b).
\ee
\et

 \paragraph{Examples} 
 
\begin{enumerate}
\item When all the marked points and all the boundary intervals are red, we get $\mathcal{Q}_{\bS_{\rm{c}}}=1$. Then the weighted exponential volume ${\cal E}_{\bS_{\rm{c}}}(\omega;\K_{\bf r},\L)={\cal E}_{\bS}(\K,\L)$ is the 
plain exponential volume. 
\item When all the marked points  are blue, and $\omega_e=\frac{1}{2}$ for any edge $e$, then $-\log \mathcal{Q}_{\bS_{\rm{c}}}$ coincides with the action considered by Chekhov in \cite[Formula 2.2]{Ch24}. Therefore the 
 weighted exponential volume ${\cal E}_{\bS_{\rm{c}}}(\omega,\L)$ is the volume studied  in \cite{Ch24}.
\end{enumerate}


  \subsection{Decompositions of  weighted exponential volumes} 

\begin{figure}[ht]
	\centering
	\includegraphics[scale=0.7]{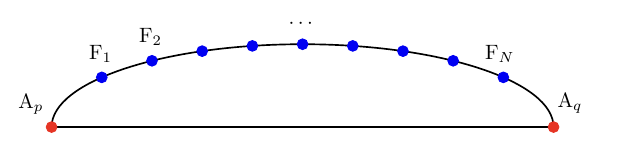}
	\small
	\caption{A red-based half-disc ${\D}^{-}_{\N}$ carries  two adjacent red points and $\N$ blue points .}
	\label{figure:interval11}
\end{figure}

\paragraph{1. A red-based half-disc.}  Consider the colored surface ${\D}^{-}_{\N}$ given by the half-disc  with the red base provided by the 
two adjacent red points, and $\N$ blue points, see Figure \ref{figure:interval11}. 
It carries $\N$ blue intervals and a single red interval. We denote by $\omega = (\omega_1, ..., \omega_\N)$ the weights on the blue intervals, and by $\K$ the $\K-$coordinate on the red edge. 
 Its exponential 
volume is denoted by 
$$
{\cal E}_{{\rm D}^{-}_{\N}}(\omega_1, ..., \omega_\N, \K):=\int_{{\cal M}_{{\D}^{-}_{\N}}}{\Bbb E}_{{\rm R}_{\N}}(\omega, \K). $$
If $\N=0$ we have ${\cal E}_{{\rm D}^{-}_{0}}(\K)=1$. 
We  denote by $\D^\bullet_\N$ the disc with a single red point and $\N$ blue points. 

\begin{figure}[ht]
	\centering
	\includegraphics[scale=0.7]{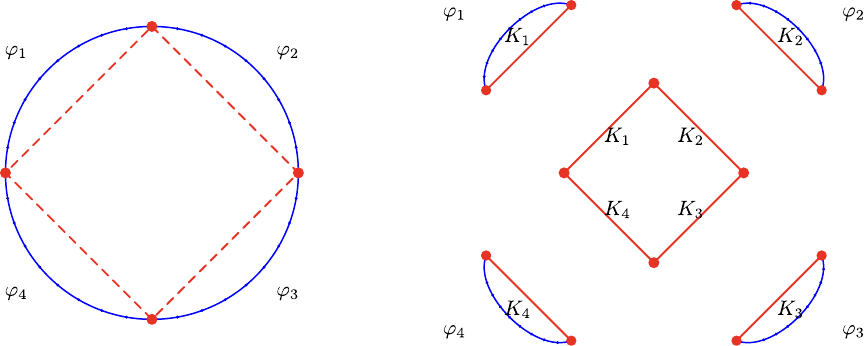}
	\small
	\caption{The weighted exponential volume for a disc with $4$ red points by cutting along the red rectangle.}
	\label{figure:circledec}
\end{figure}

\paragraph{2. Colored discs via  red-based half-discs.} Consider a colored disc $\D_{\N+ m}$ with $m>1$ red points. Cut out from the disc a red-colored  $m-$gon $\D_m$, 
see Figure \ref{figure:circledec}. Denote by $(\K_1, ..., \K_m)$ the $\K-$parameters on its sides. We get $m$ red-based half-discs with the same  $\K-$parameters $\K_1, ..., \K_m$. Denote by $\omega_{(i)}$ the ordered collection of weights on the $i-$th red-based half-disc. So the cyclically ordered collection $(\omega_1, ..., \omega_\N)$ of  weights for the disc $\D_{\N+m}$ is presented as the cyclically 
ordered union of the weights for the red-based half-discs:
\be \la{Nweights}
(\omega_1, ..., \omega_\N) = (\omega_{(1)}, ..., \omega_{(m)}).
\ee
\bp \la{T1.16} The weighted exponential volume for a colored disc $\D_{\N+m}$, see Figure \ref{figure:circledec}, is 
\be \la{CRD}
\begin{split}
&{\cal E}_{{\rm D}_{\N+m}}(\omega_1, ..., \omega_\N)\\
=&\int {\cal E}_{{\rm D}_{m}}(\K_1, ..., \K_m) \cdot {\cal E}_{{\rm D}^-_{n_1}}(\omega_{(1)}, \K_1)\cdot \ldots \cdot  {\cal E}_{{\rm D}^-_{n_m}}(\omega_{(m)}, \K_m) d\log \K_1 \wedge \ldots \wedge d\log \K_m.\\
\end{split}
\ee
\ep

\bpr 

Applying the inverse Mellin transform along the Mellin variables   (\ref{Nweights}), we reduce the claim to the similar claim about the exponential volumes, which follows from their cutting and gluing 
properties  \cite[Section 2.2]{GS24}. 
\epr

\begin{figure}[ht]
	\centering
	\includegraphics[scale=0.7]{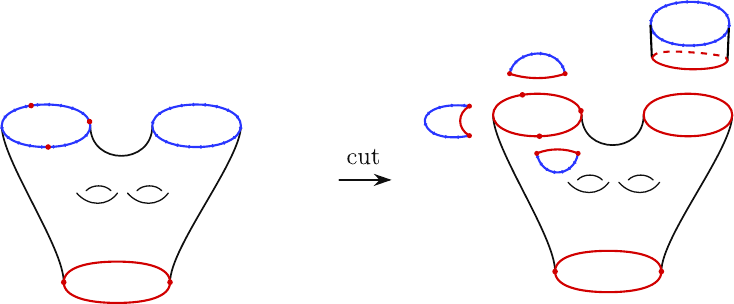}
	\small
	\caption{Cutting red-based half-discs or single-red-cusped discs  from any colored surface.}
	\label{figure:cutlimitXY}
\end{figure}

\paragraph{3. Cutting out red-based half-discs from  any colored surface.} Take any colored surface  surface $\bS_{\rm c}$, and cut out from it all red-based half-discs $\D^-_\N$, and a single-red-point discs $\D^\bullet_\N$, as shown on Figure \ref{figure:cutlimitXY}.
Denote by $\bS_{\bullet}$ the obtained  decorated surface. It has all boundary intervals colored in red. So we have a decomposition
\be \la{CAD}
\bS_{\rm c} = \bS_{\bullet}\ \cup \ {\D}^{-/\bullet}_{(1)}\cup \ldots \cup {\D}^{-/\bullet}_{(\rM)}.
\ee
Here $\rM$ is the total number of red-based half-discs / single-red-point discs which we cut out from the colored surface $\bS_{\rm c}$. It coincides with the total number of red points on $\bS_{\rm c}$. 
The decomposition (\ref{CAD}) is canonical. It gives rise to canonical integral presentation of the weighted exponential volume for arbitrary colored surface via the ones for the red-colored surfaces and 
red-based half-discs and single-red-point discs, which generalizes Proposition \ref{T1.16}, and proved the same way.

\bt \la{T1.17} The weighted exponential volume of the moduli space related to  any colored surface can be written as, see Figure. \ref{figure:cutlimitXY}:
\be
\begin{split}
&{\cal E}_{{\bS}_{\rm c}}(\omega_1, ..., \omega_\N)\\
=&\int {\cal E}_{{\bS}_{\bullet}}(\K_1, ..., \K_\rM) \cdot {\cal E}_{{\D}^{-/\bullet}_{n_1}}(\omega_{(1)}, \K_1)\cdot \ldots \cdot  {\cal E}_{{\D}^{-/\bullet}_{n_\rM}}(\omega_{(\rM)}, \K_\rM) d\log \K_1 \wedge \ldots \wedge d\log \K_\rM.\\
\end{split}
\ee
\et

\paragraph{5. The neck recursion formula.} 
Look at each crown ${\rm C}_j$ on the underlying  surface $\bS$.  Given an ideal hyperbolic structure on the surface,  denote by $l_j$ the length of the neck geodesic $\gamma_j$ on the crown ${\rm C}_j$. 
Cut out the crown ${\rm C}_j$ from the  surface by the neck geodesic $\gamma_j$. 
We get a punctured colored disc $\D^*_{(j)}$, whose puncture match the neck geodesic $\gamma_j$, and whose external part describes the colored crown ${\rm C}_j$.  
Denote by $\bS_\bullet$ the remaining balk surface, which has $m$ holes, so that 
$$
\bS_{\rm c}= \bS_\bullet \cup \cup_{j=1}^m\D^*_{(j)}.
$$
We denote by $\omega_{(j)}, \K_{(j)}$ the weights and  $\K-$coordinates for the crown ${\rm C}_j$, and set $\K =  (\K_{(1)}, ...,  \K_{(m)})$. 
\bt There is the following integral formula for the weighted exponential volume of any colored surface $\bS_{\rm c}$: 
\be \la{OPO}
\begin{split}
&{\cal E}_{\bS_{\rm c}}(\omega_1, ..., \omega_\N, \K)\\
=&\int {\cal E}_{{\bS}_{\bullet}}(l_1, ..., l_m) \cdot {\cal E}_{{\D}^*_{n_1}}(\omega_{(1)}, \K_{(1)}, l_1)\cdot \ldots \cdot  {\cal E}_{{\D}^*_{n_m}}(\omega_{(m)}, \K_{(m)}, l_m) l_1d l_1\wedge \ldots  \wedge l_md l_m.\\
\end{split}
\ee
\et

\bpr Similar to the proof of the neck recursion formula in \cite{GS24}.
\epr

 \subsection{Probabilistic Teichm\"uller and moduli space}

 \subsubsection{Schwarzian structures on $\varphi-$colored surfaces and the space ${\cal M}_{\bS_\varphi}$} 
 
 
 To define probabilistic Teichmuller spaces we  specify which boundary components of $\bS$ become probabilistic, that is carry infinite dimensional spaces of fields $\varphi$.  
 We call them $\varphi-$colored  components.

    \bd 
A ${\varphi}-$colored surface $(\bS, \varphi)$ is   a decorated surface  $\bS$  with a subset  
  of its punctured boundary components, which we  call   {\em $\varphi-$colored}   boundary components. 
 \ed
 
 \vskip 2mm
 
 Given a {${\varphi}-$colored surface} $\bS_\varphi= (\bS, {\varphi})$, denote by 
  ${\cal I}_{\varphi}$ the set of  {$\varphi-$colored boundary intervals  and by $J_{\varphi}$ the set of  {$\varphi-$colored boundary circles. So 
\be \la{SPHI}
(\bS, {\varphi}):=(\bS, \mathcal{I}_\varphi, \mathcal{J}_\varphi), \ \ \ \ {\cal I}_{\varphi}\subset {\cal I}_\bS, \ \ \ \ \mathcal{J}_\varphi\subset {\cal J}_\bS. 
\ee
 We show the $\varphi-$colored punctured boundary components on pictures as colored in green. It is handy to color the rest of the boundary components in red, although this is not an extra data on $\bS_\varphi$. 
 So 
 $$
 {\cal I}_\bS = \mathcal{I}_\varphi \cup {\cal I}_{\bf r},\;\;\;  {\cal J}_\bS = \mathcal{J}_\varphi \cup {\cal J}_{\bf r}.  
 $$
Here ${\cal I}_{\bf r}$ parametrises the red boundary intervals and ${\cal J}_{\bf r}$ parametrises the red boundary circles. 

 \vskip 1mm

Given a $\varphi-$colored surface $\bS_\varphi$, let us perform the following surgeries, see Figure \ref{figure:cutlimitX}:

\begin{figure}[ht]
	\centering
	\includegraphics[scale=0.7]{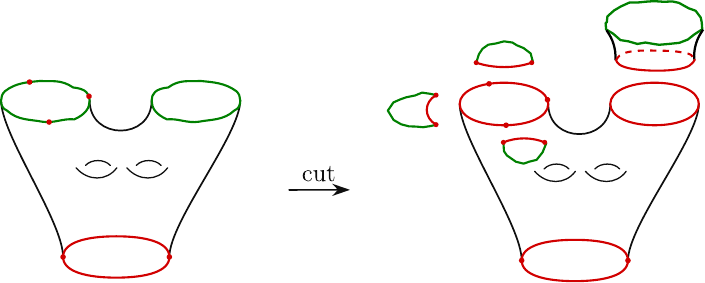}
	\small
	\caption{Cutting red-based half-discs and  discs with red-based holes from a $\varphi-$colored surface $\bS_\varphi$ on the left produces a decorated surface $\bS_\bullet$ with red boundary components on the right.}
	\label{figure:cutlimitX}
\end{figure}

 \begin{enumerate}

\item For each $\varphi-$colored boundary interval ${\rm I}_e$:
 
 Cut out from $\bS$ a red-based half-disc $\D_e$ along  the unique up to isotopy\footnote{Here isotopy refers to an isotopy with fixed endpoints} arc $\alpha_e$ in $\bS$,  
which connects the two endpoints of ${\rm I}_e$ and is   isotopic to ${\rm I}_e$,  as on right hand side of Figure \ref{figure:cutlimitX}. 
 Denote by $\K_e$ the $\K-$parameter at the diameter of the half-disc $\D_e$.

\item 
For each $\varphi-$colored boundary circles $S^1_\gamma$:

Cut the surface $\bS$ along a neck loop  $\gamma$,  isotopic to the boundary circle $S^1_\gamma$. Denote by $\D^*_\gamma$ the obtained  disc with a hole, so that 
 the neck loop becomes the boundary of the  hole. 
 

\end{enumerate}

Denote by $\bS_\bullet$  the resulting decorated surface with red boundary components :
 \be \la{SBUL}
 \bS_\bullet:= \bS_\varphi- \bigcup_{e \in {\cal I}_{\varphi}}\D_{e} \ \cup \ \bigcup_{\gamma \in {\cal J}_{\varphi}} \D^*_{\gamma}. 
\ee
Note that an ideal hyperbolic structure on the decorated  surface $\bS_\bullet$ gives rise to the following data: 

 \begin{itemize}
\item For each  boundary interval ${\rm I}_i$, $i \in {\cal I}_{{\bS}_\bullet}$,  the $\K-$invariant $\K_i$ of ${\rm I}_i$.

\item 
 For each boundary circle $S^1_\gamma$, 
  the neck geodesic  loop length $l_\gamma$.  
Set  $\L_\gamma:= e^{l_\gamma}$. 
\end{itemize} 
 Let us set 
\be  \la{38a}
\begin{split}
&\K_{\bullet}:= \{\K_i\}, \ \  \ \ {\ \ i \in \cal I_{\bS_\bullet}}.\\
&\L_{\bullet} :=\{\L_\gamma\}, \ \ \ \ \ \ {\gamma\in {\cal J}_{\bS_\bullet}}.\\
\end{split}
\ee

 Denote by  $\mathrm{Diff}^1_+([0,\T]\rightarrow [0,1])$ the space of strictly monotone  ${\rm C}^1-$maps $\varphi: [0,\T] \lra [0,1]$.  
Set \be
\label{eq:intdiff}
\mathrm{Diff}^1_+(\lambda, \T):=\{\varphi  \in \mathrm{Diff}^1_+([0,\T]\rightarrow [0,1])\;|\; \varphi(0)=0, \ \varphi(\T)=1, \  \ \varphi '(0)=\varphi '(\T)=\pi/\lambda\}.
\ee

\bd \la{DSCHS}
A {\em Schwarzian structure} on a $\varphi-$colored surface $(\bS, \varphi)$ is given by 
 the following:

\begin{enumerate}

\item  An ideal hyperbolic structure on the decorated surface $\bS_\bullet$ in (\ref{SBUL})  with the  boundary data $(\K_\bullet, \L_\bullet)$. 

\item For each $\varphi-$colored boundary interval ${\rm I}_e$, a 
positive number   $\T_e$, and a ${\rm C}^1$ monotone map 
\be \begin{split}
&\varphi_e\in \mathrm{Diff}^1_+(\lambda_e,\T_e), \ \ \ \lambda_e:=\frac{1}{2}\log \K_e.\\
\end{split}
\ee 

 \item For each  $\varphi-$colored boundary circle $S^1_\gamma$, a positive 
number $\T_\gamma$ and a  $(\T_\gamma, l_\gamma)-$quasiperiodic  monotone ${\rm C}^1$ function $\varphi_\gamma(t)$:
 \be
 \varphi_\gamma \in \mathrm{C}^1(\R);  \ \ \  \varphi_\gamma(t+\T_\gamma)=\varphi_\gamma(t)+l_\gamma,\ \ \ \varphi_\gamma'(t)> 0 \ \ \ \forall t \in \R; \ \ \ \varphi(0)=0.
\ee
Here $l_\gamma$ is  
  the neck geodesic  loop length of the  boundary circle $S^1_\gamma$ on $\bS_\bullet$.  
 \end{enumerate}   
\ed

 \bd  The probabilistic Teichmuller space ${\cal T}_{\varphi}(\bS)$ parametrises the Schwarzian structures on $\bS_\varphi$ on  a $\varphi-$colored surface $\bS_\varphi$. 

 The probabilistic moduli space ${\cal M}_{\varphi}(\bS)$ is its quotient by the action of the mapping class group:
  $$
  {\cal M}_{\varphi}(\bS):={\cal T}_{\varphi}(\bS)/{\rm Mod}_\bS. 
  $$
 \ed

 We introduce a measure ${\bf E}$ on the probabilistic moduli space ${\cal M}_{\bS_\varphi}$. 
 When $\bS$ is a disc, the  ${\cal M}_{\bS_\varphi}$ is identified with  the quotient ${\rm Diff}(S^1)/{\rm PGL}_2(\R)$ of the group of ${\rm C}^1-$diffeomorphisms  of the circle. Then the  measure is the Schwarzian measure, 
 considered 
 by Stanford and Witten \cite{SW17}, and defined rigorously by Bauerschmidt, Losev and Wildemann \cite{BLW24} using the Wiener measure and the Malliavin-Shavgulidze map. 
 The measures for a red-based $\varphi-$colored half-disc and a red-based $\varphi-$colored annulus are cousins of the Schwarzian measure,  defined using the Brownian bridge. The measure on the space 
 ${\cal M}_{\bS_\varphi}$ is defined by gluing these measures with the Weil-Petersson measure on ${\cal M}_{\bS_\bullet}(\K_\bullet, \L_\bullet)$ following the pattern  on Figure \ref{figure:cutlimitX}. 
 So  the  measure  ${\bf E}$ generalise  the Weil-Petersson measure, the exponential measure  \cite{GS24} and the measure  \cite{Ch24},  and the Schwarzian measure \cite{BLW24}.  \vskip 1mm
 
  Both   spaces ${\cal T}_{\varphi}(\bS)$ and  ${\cal M}_{\varphi}(\bS)$ depend on the  parameters
\be  \la{38}
\begin{split}
&\K_{\bf r}:= \{\K_i\}\ \ \ \ \ \ \  \ \ {i \in \cal I_{\bf r}}.\\
&\L_{\bf r} :=\{\L_\gamma\}, \ \ \ \ \ \ \ \ {\gamma\in {\cal J}_{\bf r}}.\\
&\T_\bullet:= \{\T_e, \T_\gamma\}, \ \ \ e \in {\cal I}_\varphi, \gamma \in {\cal J}_\varphi.\\
\end{split}
\ee
Denote by  $\mathcal{M}_{\bS_\varphi}(\T_\bullet,\K_{\bf r}, \L_{\bf r})$ the subspace with the fixed parameters (\ref{38}). 
We prove in Theorem \ref{Th6.14} that the measure of the moduli space $\mathcal{M}_{\bS_\varphi}(\T_\bullet, \K_{\bf r}, \L_{\bf r})$  is finite. 
   It describes a field theory which can be called {\it Schwarzian field theory on decorated surfaces}. The classical Schwarzian field theory  \cite{SW17} corresponds to the case when $\bS$ is a 
   $\varphi-$decorated disc.    
   \vskip 1mm
   
   Section \ref{SEC6} can be viewed as  the continuous analog of Section \ref{SEC4}. Precisely, 
   L. Chekhov argued  \cite{Ch24} that the volume form for the moduli space of hyperbolic polygons   in loc. cit.  is a discretization of the Schwarzian measure for the disc. 
This was proved  rigorously in the forthcoming paper of Budd and Chekhov, cf \cite{BC}. 
The weighted exponential volume forms assigned to a weighted colored surface in Section \ref{SEC4} can be viewed a  discretization of the measure ${\bf E}$ for decorated surfaces, see Section \ref{SEC8}.
  
\vskip 1mm

   In Section \ref{Sect9} we prove    recursions  for the measures on the probabilistic moduli spaces $\mathcal{M}_{\bS_\varphi}(\T_\bullet, \K_{\bf r}, \L_{\bf r})$, and hence their volumes, see Theorem \ref{Th7.8}. It generalizes Mirzakhani's recursions \cite{Mir07b} 
 for the Weil-Petersson volumes, and recursions   for the exponential volumes \cite{GS24}. 
 Their novel feature  is that they give  recursions for a certain class of Feynman integrals. 
 Recall that only some boundary intervals are $\varphi-$colored, and 
 we color the  rest in red. An important feature of the recursions  
  is this:  
  \vskip 1mm
{\it Recursions for probabilistic moduli spaces ${\cal M}_{\bS_\varphi}$ always require the ones with red boundary intervals}.

 \paragraph{\it Conclusion.}     We suggest that {the {probabilistic moduli spaces} ${\cal M}_{\bS_\varphi}$ are the true infinite dimensional analogs of the classical moduli spaces of hyperbolic metrics on surfaces}.

\subsubsection{Acknowledgements}  This work benefited a lot from the support of IHES. 
We enjoyed the hospitality of IHES during the Summer of 2025, when  the  Barnes integral presentation was worked out,   during the August  of 2026. 
We are grateful to IHES for the hospitality and support. 

The work of A.G. was supported by the NSF grant DMS-2153059. 

A.G. is grateful to Leonid Chekhov for several useful discussions in June 2026 at IHES, and to Ilya Losev for discussions at CIRM in August 2026.

\section{Barnes integrals for  decorated surfaces: examples}
\label{sectionex}

We start from the three  fundamental examples 
  of the Mellin--Laplace transforms of exponential volumes, presented in Sections \ref{Tr}, \ref{Sec3.2} $\&$ \ref{sec3.2a}, \ref{Cr}. These examples will be  used in further calculations.

\subsection{Triangle} \la{Tr}
\begin{figure}[ht]
	\centering
	\includegraphics[scale=0.25]{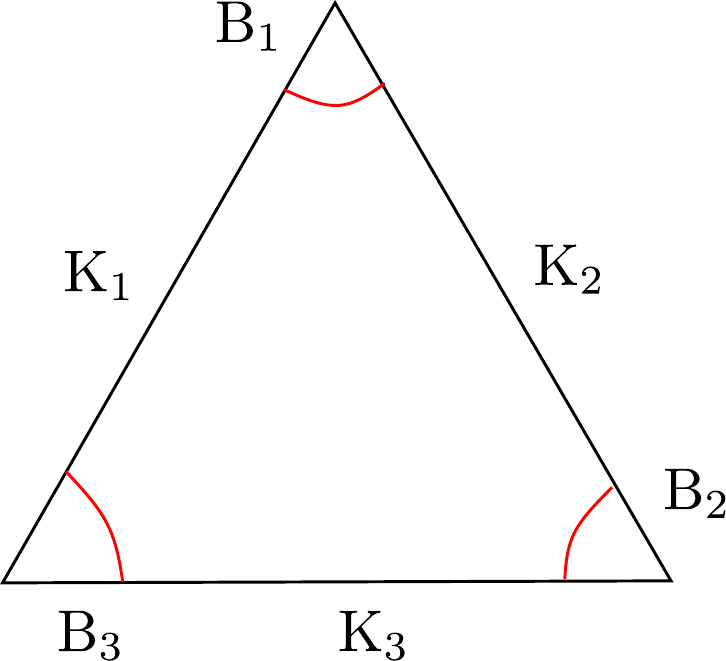}
	\small
	\caption{The ideal triangle with the $\K-$ and $\B-$coordinates.}
	\label{figure:D3}
\end{figure}

Consider an ideal hyperbolic structure on a disc  $\Delta$ with $3$ marked points on the boundary, called a hyperbolic triangle. As in Figure \ref{figure:D3}, the three $\K$-coordinates are $\K_1$, $\K_2$, $\K_3$. 
Following \cite{GS24}, we also use the $\B-$coordinates. In this example the $\B$-coordinates and $\K$-coordinates are related by 
\be \la{BK}
\B_1 \B_2 = \K_2,\;\;\; \B_2 \B_3=\K_3, \;\;\; \B_3 \B_1=\K_1.
\ee
The potential
\[W=\B_1+\B_2+\B_3= \sqrt{\frac{\K_1 \K_2}{\K_3}}+ \sqrt{\frac{\K_2 \K_3}{\K_1}}+ \sqrt{\frac{\K_3 \K_1}{\K_2}}.\]
The exponential volume
\[\mathcal{E}_{\Delta}(\K_1,\K_2,\K_3)=e^{-W}= 
 e^{-(\B_1+\B_2+\B_3)}.\]
Then the Mellin transform of exponential volume is calculated via $\Gamma-$functions as follows:
\be
\label{eq:D3}
\begin{aligned}
&\mathcal{M}(\mathcal{E}_{\Delta}; s_1,s_2,s_3)=\int_0^{+\infty}\int_0^{+\infty}\int_0^{+\infty}e^{-W} \K_1^{s_1} \K_2^{s_2} \K_3^{s_3}  d\log \K_1 d\log \K_2 d\log \K_3
\\&=\int_0^{+\infty}\int_0^{+\infty}\int_0^{+\infty}e^{-(\B_1+\B_2+\B_3)} (\B_3 \B_1)^{s_1} (\B_1 \B_2)^{s_2} (\B_2 \B_3)^{s_3}  2 d \log \B_1 d \log \B_2 d \log \B_3
\\&=2 \int_0^{+\infty} e^{-\B_1} \B_1^{s_1+s_2} d \log \B_1 \int_0^{+\infty} e^{-\B_2} \B_2^{s_2+s_3} d \log \B_2 \int_0^{+\infty} e^{-\B_3} \B_3^{s_3+s_1} d \log \B_3
\\&=2 \Gamma(s_1+s_2) \Gamma(s_2+s_3) \Gamma(s_3+s_1).
\end{aligned}
\ee
\subsection{Quadrilateral} \la{Sec3.2}
\begin{figure}[ht]
	\centering
	\includegraphics[scale=0.23]{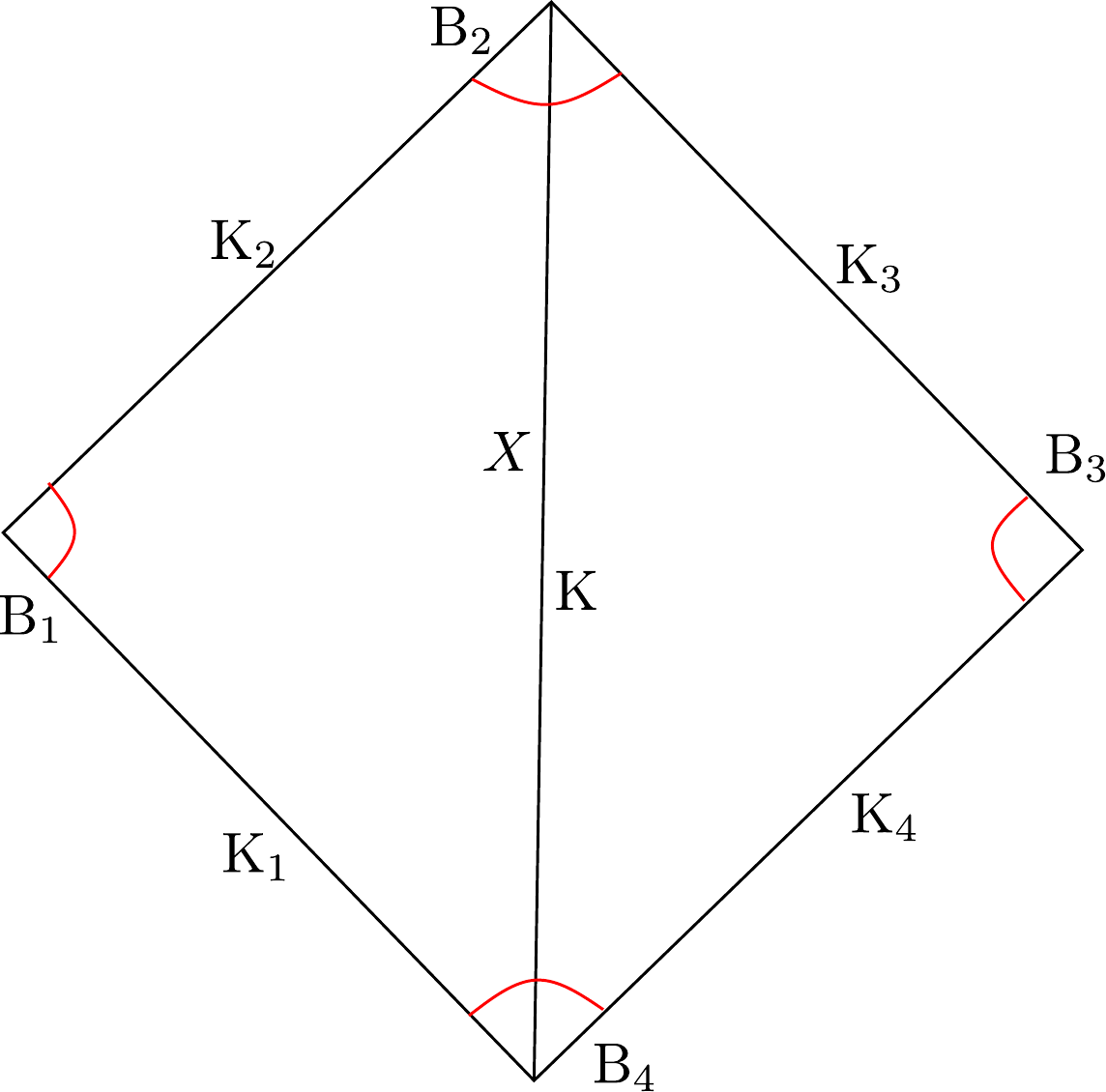}
	\small
	\caption{The ideal quadrilateral.}
	\label{figure:D4}
\end{figure}
Consider an ideal hyperbolic structure on a disc  $\square$ with $4$ marked points on the boundary, called a hyperbolic quadrilateral. As in Figure \ref{figure:D4}, the four $\K$-coordinates are $\K_1$, $\K_2$, $\K_3$, $\K_4$. A diagonal of the quadrilateral cuts the quadrilateral into two triangles $\Delta$ and $\Delta'$, where $\Delta$ carries the $\K$-coordinates $\K_1$, $\K_2$ along the boundary while $\Delta'$ carries the $\K$-coordinates $\K_3$, $\K_4$. Let $\K$ be the $\K$-coordinate along the diagonal. Given the corresponding $\B$-coordinates $\B_1$, $\B_2$, $\B_3$, $\B_4$, we have
\[\B_1 \B_2=\K_2,\;\; \B_2 X \B_3=\K_3,\;\; \B_3 \B_4= \K_4,\;\; \B_4 X \B_1=\K_1.\]
The potential
\[W= W_{\Delta}+W_{\Delta'}= \left(\sqrt{\frac{\K_1 \K_2}{\K}}+ \sqrt{\frac{\K \K_1 }{\K_2}}+ \sqrt{\frac{\K_2 \K}{\K_1}}\right) + \left(\sqrt{\frac{\K_3 \K_4}{\K}}+ \sqrt{\frac{\K \K_3 }{\K_4}}+ \sqrt{\frac{\K_4 \K}{\K_3}}\right).\]
The exponential volume
\[\mathcal{E}_{\square}(\K_1,\K_2,\K_3, \K_4)=\frac{1}{2}\int_0^{+\infty} \mathcal{E}_{\Delta}(\K_1,\K_2,\K) \mathcal{E}_{\Delta}(\K_3,\K_4,\K) d \log \K.\]

\bp The Mellin transform of the exponential volume function is  calculated via $\Gamma-$functions:
\be
\label{eq:334}
\begin{aligned}
\mathcal{M}(\mathcal{E}_{\square}; s_1,s_2,s_3,s_4) 
&=\frac{2\Gamma(s_1+s_2)\Gamma(s_3+s_4)}{2\pi i} \int_{-i \infty}^{+i \infty}  \Gamma(s_1+s)  \Gamma(s_2+s)  \Gamma(s_3-s)  \Gamma(s_4-s) ds 
\\&\overset{\rm Barnes}{=}\frac{2\prod_{1\leq i<j\leq 4}\Gamma(s_i+s_j)}{\Gamma(s_1+s_2+s_3+s_4)}.\end{aligned}
\ee
\ep

\bpr
Let us recall the Parseval's formula and Barnes' lemma. Let $\Re s$ be the real part of $s$.
\bl[Parseval's formula]
Let $f(\K)$ and $g(\K)$ be functions with well-defined Mellin transforms $\mathcal{M}(f,s)$ and $\mathcal{M}(g,s)$ in the fundamental strips 
$$
a_1<\Re s<b_1 \ \mbox{and} \ a_2<\Re s<b_2, \ \mbox{where $\max\{a_1,-b_2\}<\Re c<\min\{b_1,-a_2\}$}.
$$ 
If the function $\K^{c-\frac{1}{2}} f(\K)$ and $\K^{-c-\frac{1}{2}} g(\K)$ are square-integrable over the interval $(0,+\infty)$, we have
\be
\int_0^{+\infty} f(\K) g(\K) d \log \K
=\frac{1}{2\pi i}\int_{c-i \infty}^{c+i \infty} \mathcal{M}(f,s) \mathcal{M}(g,-s) ds
\ee
\el

\bl[Barnes' formula]
We have the formula
\be \la{7}
\frac{1}{2\pi i} \int_{-i\infty}^{+i\infty} \Gamma(a+s)\Gamma(b+s)\Gamma(c-s)\Gamma(d-s) ds=
\frac{\Gamma(a+c)\Gamma(a+d)\Gamma(b+c)\Gamma(b+d)}{\Gamma(a+b+c+d)},
\ee
where the  integration  contour $\gamma$ is a vertical line from infinity to infinity,  denoted symbolically by $(-i\infty, i\infty)$,  which separates the poles of $\Gamma(a+s)$ and $\Gamma(b+s)$ from those of $\Gamma(c-s)$ and $\Gamma(d-s)$, 
and oriented so that the points $a, b$ are on the left, and the points $c,d$ are on the right. 
\el

\paragraph{Remark.} Let $\overline \gamma$ be the countour $\gamma$ with the opposite orientation. Then  the left hand side of formula (\ref{7})  is invariant under the cyclic shift 
$$
(a,b,\gamma, s, c,d) \lms (c,d, \overline \gamma, -s, a,b).
$$
Indeed, the substitution   $s \lms -s$ 
changes both the integrand and the countour orientation, so that 
$$
 \int_\gamma \Gamma(a+s)\Gamma(b+s)\Gamma(c-s)\Gamma(d-s) ds =   \int_{\overline \gamma} \Gamma(a-s)\Gamma(b-s)\Gamma(c+s)\Gamma(d+s) ds.
 $$ 
\vskip 2mm

By Parseval's formula for the Mellin transform, the Mellin transform of exponential volume 
\be
\begin{aligned}
&\mathcal{M}(\mathcal{E}_{\square};s_1,s_2,s_3,s_4)=\frac{1}{2}\int_0^{+\infty} \mathcal{M}(\mathcal{E}_{\Delta}(\K);s_1,s_2) \mathcal{M}(\mathcal{E}_{\Delta'}(\K);s_3,s_4)  d \log \K
\\&=\frac{1}{4\pi i}\int_{-i \infty}^{+i \infty} \mathcal{M}(\mathcal{E}_{\Delta};s_1,s_2,s) \mathcal{M}(\mathcal{E}_{\Delta'};s_3,s_4,-s) ds
\end{aligned}
\ee
where the integration contour  is the line from bottom to top containing in the intersection of the well-defined fundamental strips of $\mathcal{M}(\mathcal{E}_{\Delta};s_1,s_2,s)$ and $\mathcal{M}(\mathcal{E}_{\Delta'};s_3,s_4,-s) $ for $s$. 
Then by  (\ref{eq:D3}) we have
\be
\label{eq:33to4}
\begin{aligned}
&\mathcal{M}(\mathcal{E}_{\square};s_1,s_2,s_3,s_4)=\frac{1}{2}\int_{-i \infty}^{+i \infty} \mathcal{M}(\mathcal{E}_{\Delta};s_1,s_2,s) \mathcal{M}(\mathcal{E}_{\Delta'};s_3,s_4,-s) ds
\\&=\frac{1}{4\pi i}\int_{-i \infty}^{+i \infty} 2 \Gamma(s_1+s) \Gamma(s_1+s_2) \Gamma(s_2+s) 2 \Gamma(s_3-s) \Gamma(s_3+s_4) \Gamma(s_4-s) ds
\\&=\frac{1}{\pi i}\Gamma(s_1+s_2)\Gamma(s_3+s_4) \int_{-i \infty}^{+i \infty}  \Gamma(s_1+s)  \Gamma(s_2+s)  \Gamma(s_3-s)  \Gamma(s_4-s) ds.
\end{aligned}
\ee
Then by Barnes' lemma, we get
\be
\label{eq:square4}
\mathcal{M}(\mathcal{E}_{\square};s_1,s_2,s_3,s_4)
\overset{\rm{Barnes}}{=}\frac{2\prod_{1\leq i<j\leq 4}\Gamma(s_i+s_j)}{\Gamma(s_1+s_2+s_3+s_4)}.
\ee
\epr

\subsection{Polygon} \la{sec3.2a}
\begin{figure}[ht]
	\centering
	\includegraphics[scale=0.3]{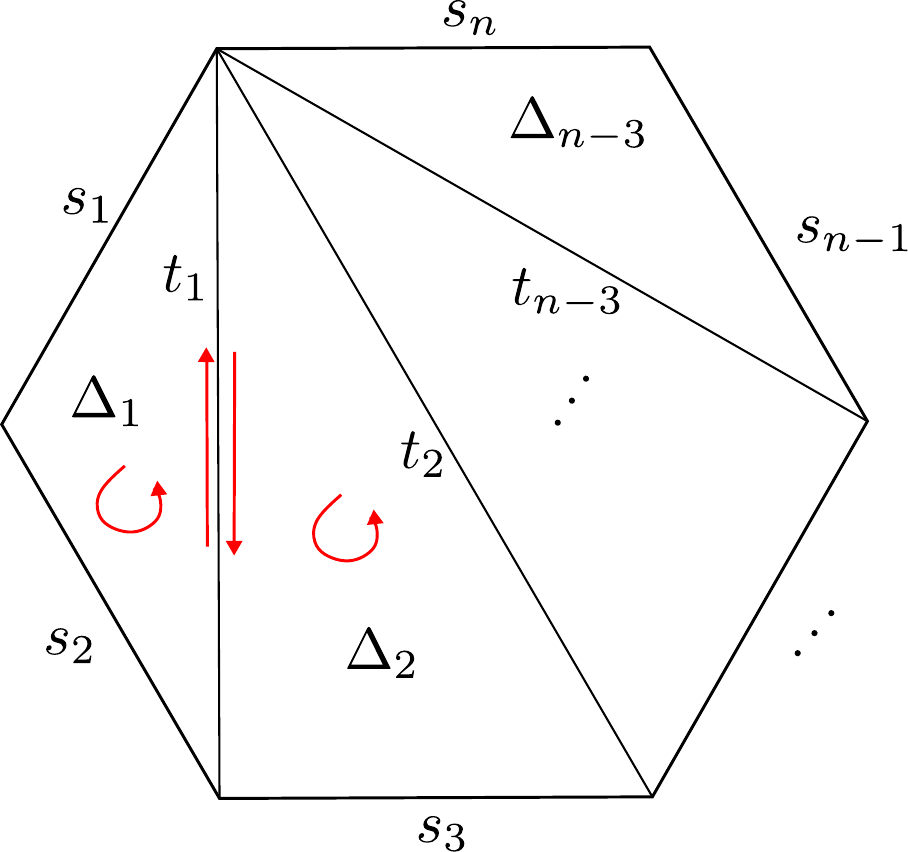}
	\small
	\caption{Triangulated $n$-gon.}
	\label{figure:Dn}
\end{figure}

Consider the Mellin transform of the exponential volume of the counterclockwise oriented polygon $\D_n$.
Cut   the polygon $\D_n$ into $(n-3)$ ideal triangles $\Delta_1$, $\cdots$, $\Delta_{n-3}$ as on Figure \ref{figure:Dn}. Denote by $s_1$, $\cdots$, $s_n$ the Mellin parameters on the boundary  edges,  and by $t_1$, $\cdots$, $t_{n-3}$ the ones on the interior  edges. They are the  Mellin dual parameters  for the $\K$-coordinates on these edges. 

We assign to each angle $\theta$ of an ideal triangle a linear function $\ell(\theta)$ of $s_1$, $\cdots$, $s_n$, $t_1$, $\cdots$, $t_{n-3}$. If the angle $\theta$ is bounded by two ideal edges which carry parameters $a$ and $b$, then $\ell(\theta)=\pm a\pm b$, where the signs can be negative only in the following case: 
\begin{itemize}
\item if $a$ or $b$ equals $t_i$ for some $i$, the anti-clockwise orientation in each ideal triangle induces orientation of each interior ideal edge. Thus when two different angles share a common ideal edge and on two opposite sides, the sign in one is opposite to the other. 
\end{itemize}
Then we can rewrite the second line of Equation (\ref{eq:33to4}) in a uniform way as 
\[\frac{2}{2\pi i} \int_{-i \infty}^{+i \infty}  \prod_\theta \Gamma(\ell(\theta)) ds.
\]
For the $n$-gon $D_n$, by applying Parseval's formula $(n-3)$ times, we get
\be \la{27}
\mathcal{M}(\mathcal{E}_{\D_n};s_1,\cdots,s_n)=\frac{2}{(2\pi i)^{n-3}} \int_{-i \infty}^{+i \infty}\cdots\int_{-i \infty}^{+i \infty}  \prod_\theta \Gamma(\ell(\theta)) d t_1 \cdots d t_{n-3},
\ee
where the contours of integrations are the lines  in the intersections of the well-defined fundamental strips.

 \subsection{Crown with $\N$ cusps $\D_\N^*$} \la{Cr}
\label{subsection:d1}
\subsubsection{Barnes integral presentation of the Mellin--Laplace transform for $\N=1$}
\begin{figure}[ht]
	\centering
	\includegraphics[scale=0.3]{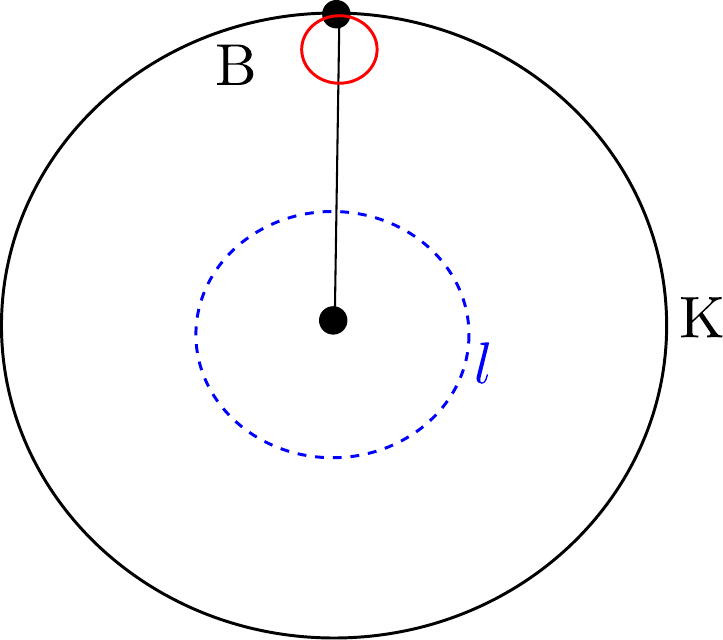}
	\small
	\caption{The trouser leg.}
	\label{figure:D1}
\end{figure}
Consider an ideal hyperbolic structure on the  once-punctured disc with one marked point on the boundary  $\D_1^*$, called a trouser leg. As in Figure \ref{figure:D1}, the unique boundary arc has the $\K$-coordinate $\K$ and the $\B$-coordinate $\B$, the unique simple closed geodesic has length $l$ and $\L=e^l$. Then 
\[\K=\B^2 \L.\]
The potential 
\[W=\B+\B\L=\sqrt{\K}(e^{\frac{l}{2}}+e^{-\frac{l}{2}}).\]
The exponential volume
\[\mathcal{E}_{\D_1^*}(\K,\L)=e^{-W}=e^{-\sqrt{\K}(e^{\frac{l}{2}}+e^{-\frac{l}{2}})}.\]

\bp The Mellin--Laplace transform of exponential volume for $\D_1^*$ can be presented by a Barnes integral:
\be 
\label{eq:d1compx}
\begin{aligned}
&\mathcal{ML}(\mathcal{E}_{\D_1^*};s,z)= 
 \frac{1}{\pi i} \int_{-i \infty}^{+i \infty} \frac{\Gamma(x)\Gamma(2s-x)\Gamma(s+z-x)}{\Gamma(s+z+1-x)} d x.
\end{aligned}
\ee
 The  integration contour  is a vertical line separating the poles of $\Gamma(2s-x)$ and $\Gamma(s+z-x)$ from those of $\Gamma(x)$.
\ep

\bpr
The Mellin--Laplace transform of the exponential volume is calculated as follows:
\be
\label{eq:d1comp}
\begin{aligned}
&\mathcal{ML}(\mathcal{E}_{\D_1^*};s,z)=\int_0^{+\infty}\int_0^{+\infty}e^{-(\sqrt{\K}(e^{\frac{l}{2}}+e^{-\frac{l}{2}}))} \K^{s} e^{-lz}  d\log \K dl
\\&= 2\Gamma(2s) \int_0^{+\infty}(e^{\frac{l}{2}}+e^{-\frac{l}{2}})^{-2s} e^{-lz}   dl
\\&= 2\Gamma(2s) \int_0^{+\infty} (1+e^{-l})^{-2s} e^{-lz-ls}   dl
\\&= 2\Gamma(2s) \int_0^{+\infty}\sum_{k=0}^{+\infty} \left( e^{-l(s+z)} \frac{(2s)_k}{k!}(-e^{-l})^k \right)  dl
\\&= 2\Gamma(2s) \sum_{k=0}^{+\infty} \left(\frac{(2s)_k}{s+z+k}\cdot \frac{(-1)^k}{k!}\right)
\\&= \frac{2\Gamma(2s)}{s+z} {}_2F_1(\{2s,s+z\},\{s+z+1\},-1)
\\&= \frac{2\Gamma(2s)}{s+z} \frac{\Gamma(s+z+1)}{\Gamma(2s)\Gamma(s+z)}\frac{1}{2\pi i} \int_{-i \infty}^{+i \infty} \frac{\Gamma(x)\Gamma(2s-x)\Gamma(s+z-x)}{\Gamma(s+z+1-x)} d x
\\&= \frac{1}{\pi i} \int_{-i \infty}^{+i \infty} \frac{\Gamma(x)\Gamma(2s-x)\Gamma(s+z-x)}{\Gamma(s+z+1-x)} d x.
\end{aligned}
\ee
\epr

\subsubsection{Mellin transform  of the integral of the exponential volume over $l$ for $\N=1$}

\bl The value of the function $\mathcal{ML}(\mathcal{E}_{\D_1^*};s,z)$ at $z=0$ is given by $$
\mathcal{ML}(\mathcal{E}_{\D_1^*};s,0)= \int_{0}^{+\infty} \mathcal{E}_{\D_1^*}(\K,l) dl =\Gamma(s)^2.
$$
\el

\bpr
Recall that 
$
\mathcal{E}_{\D_1^*}(\K,\L) = e^{-\sqrt{\K}(e^{\frac{l}{2}}+e^{-\frac{l}{2}})}.
$ The Mellin transform of its integral over $l$  is:
\be \la{16}
\begin{aligned}
&\mathcal{M}(\int_{0}^{+\infty} \mathcal{E}_{\D_1^*}(\K,\L) dl,s)=\int_0^{+\infty}\int_{0}^{+\infty}e^{-\sqrt{\K}(e^{\frac{l}{2}}+e^{-\frac{l}{2}})} dl \K^{s} d\log \K
\\&= 2\Gamma(2s) \int_0^{+\infty}(e^{\frac{l}{2}}+e^{-\frac{l}{2}})^{-2s}   dl
\\&= 2\Gamma(2s) \int_0^{+\infty}\left(\frac{e^{-\frac{l}{2}}}{1+e^{-l}}\right)^{2s}   dl
\\&= 2\Gamma(2s) \int_0^{+\infty} e^{-ls}\left(\frac{1}{1+e^{-l}}\right)^{2s}   dl.
\end{aligned}
\ee
Now let $x=\frac{1}{1+e^{-l}}$ and integrate over $x$, the right hand side of the above equation is
\be
\begin{aligned}
(\ref{16})&= 2\Gamma(2s) \int_{\frac{1}{2}}^1 \left(\frac{1-x}{x}\right)^s x^{2s} \left(\frac{1-x}{x}\right)^{-1} \frac{1}{x^2}  d x 
\\&= 2\Gamma(2s) \int_{\frac{1}{2}}^1 (1-x)^{s-1} x^{s-1}  d x 
\\&= \Gamma(2s) \int_{0}^1 (1-x)^{s-1} x^{s-1}  d x
\\&= \Gamma(2s) \B(s,s)
=\Gamma(s)^2.
\end{aligned}
\ee
\epr

\subsubsection{ The Mellin--Laplace  transform for $\N=1$ via the function $\gamma(a,b)$}
Let us introduce the following function, discussed in detail in Section \ref{Section7}.
\bd The function $\gamma(a,b)$  is given by 
\be \la{gammaab}
\gamma(a,b)=2^{a+b-1}\cdot 
\int_{ \A \geq \B\geq 0} e^{-(\A+\B)} \A^{a} \B^{b} d \log \A d \log \B.
\ee
\ed

\bl One has 
\be
\label{eq:d1gamma}
\mathcal{ML}(\mathcal{E}_{\D_1^*}; s,z)=2^{2-2s}\gamma(s-z,s+z). 
\ee
\el

\bpr  
Let $\A=\B\L$. Then 
\be \la{gammaba}
\begin{aligned}
&\mathcal{ML}(\mathcal{E}_{\D_1^*}; s,z)=\int_0^{+\infty}\int_0^{+\infty}e^{-W} \K^{s} e^{-lz}  d\log \K dl
\\&=\int_{0\leq \B\leq \A} e^{-(\A+\B)} (\A \B)^{s} \left(\frac{\A}{\B}\right)^{-z}  d\log (\A \B) d \log\left(\frac{\A}{\B}\right)
\\&= 2\int_{0\leq \B\leq \A} e^{-(\A+\B)} \A^{s-z} \B^{s+z} d\log  \A d \log \B
\\&=2^{2-2s}\gamma(s-z,s+z).
\end{aligned}
\ee
\epr



\subsubsection{Crown with two cusps}
\begin{figure}[ht]
	\centering
	\includegraphics[scale=0.27]{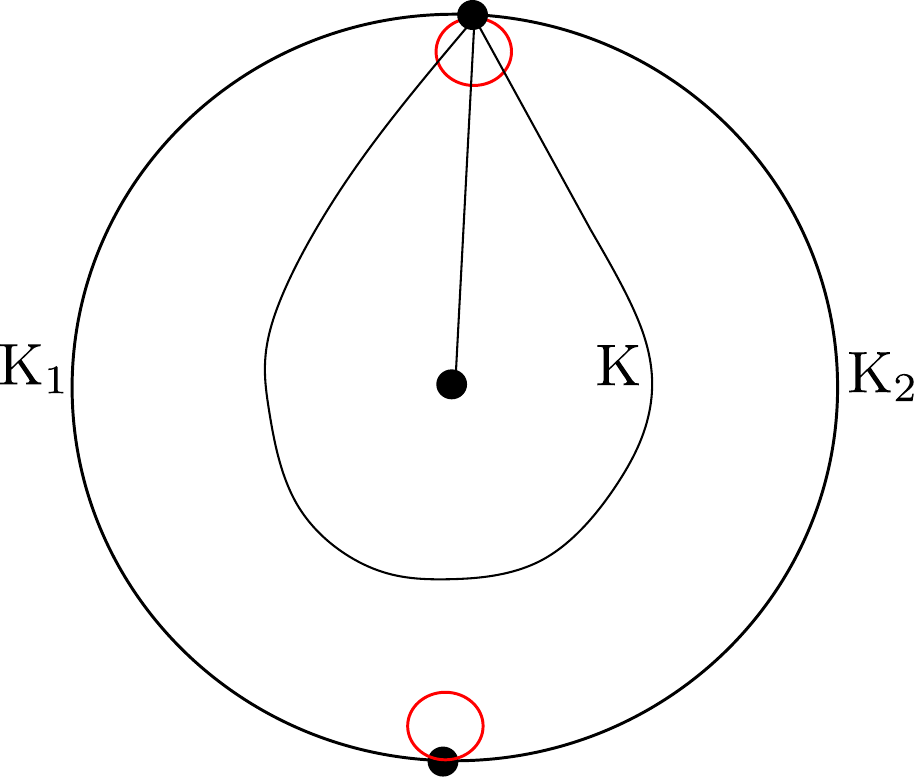}
	\small
	\caption{The crown with two cusps.}
	\label{figure:D2}
\end{figure}
Consider an ideal hyperbolic structure on a crown $\D_2^*$ with two cusps. As in Figure \ref{figure:D2}, the $\K$-coordinates of two boundary arcs are $\K_1$, $\K_2$. The ideal edge of the ideal triangulation cuts $\D_2^*$ into a crown $\D_1^*$ with one cusp and an ideal triangle $\Delta$, where $\D_1^*$ carries the $\K$-coordinate $\K$ and the length $l$ of the unique simple closed geodesic and $\L=e^l$, while $\Delta$ carries the $\K$-coordinates $\K$, $\K_1$, $\K_2$. Then \[\K=\B^2 \L.\]

\bp \la{Prop6.1} One has a Barnes integral presentation of the Mellin-Laplace transform
\be
\begin{aligned}
&\mathcal{ML}(\mathcal{E}_{\D_2^*}; s_1,s_2,z)= 
\\&  \frac{2}{(2\pi i)^2} \Gamma(s_1+s_2) \int_{-i\infty}^{+i\infty} \int_{-i \infty}^{+i \infty} \frac{\Gamma(x)\Gamma(2s-x)\Gamma(s+z-x)}{\Gamma(s+z+1-x)} \Gamma(s_1-s)\Gamma(s_2-s) d x d s.\end{aligned}
\ee
where the integration contour  is a vertical line  containing in the intersection of the well-defined fundamental strips of $\mathcal{ML}(\mathcal{E}_{\D_1^*};s,z)$ and $\mathcal{M}(\mathcal{E}_{\Delta};-s,s_1,s_2)$ for $s$. 
\ep

\bpr
Gluing along the ideal edge with the $\K$-coordinate $\K$, the exponential volume
\[\mathcal{E}_{\D_2^*}(\K_1,\K_2, \L)=\frac{1}{2}\int_0^{+\infty} \mathcal{E}_{\D_1^*}(\K,\L) \E_{\Delta}(\K,\K_1,\K_2)  d \log \K.\]
Then by Parseval's formula, the Mellin--Laplace transform of exponential volume
\be
\begin{aligned}
&\mathcal{ML}(\mathcal{E}_{\D_2^*}; s_1,s_2,z)=\frac{1}{2}\int_0^{+\infty} \mathcal{L}(\mathcal{E}_{\D_1^*}(\K),z) \mathcal{M}(\mathcal{E}_{\Delta}(\K);s_1,s_2) d\log \K
\\&=\frac{1}{4\pi i}\int_{-i\infty}^{+i\infty} \mathcal{ML}(\mathcal{E}_{\D_1^*};s,z) \mathcal{M}(\mathcal{E}_{\Delta};-s,s_1,s_2) d s.\\
\end{aligned}
\ee
Substituting (\ref{eq:D3}) and  (\ref{eq:d1compx}) we get the result.   
\epr

\subsubsection{Crown with $\N$ cusps $\D_\N^*$}

There is  a Barnes integral presentation for the Mellin-Laplace transform 
$\mathcal{ML}(\mathcal{E}_{\D_\N^*}; s_1,\ldots s_\N,z)$ of the exponential volume for the punctured disc $\D_\N^*$ generalising  the proof of Proposition \ref{Prop6.1} as follows. 
We start from a loop $\gamma$ containing one of the cusps $c$ and the puncture $p$. Cutting out the loop $\gamma$ we get a polygon $\D_{\N+1}$ and a punctured disc $\D_1^*$. Equivalently, 
the disc $\D_\N^*$ is obtained  by gluing the punctured disc $\D_1^*$ and the polygon along an ark $\gamma$:
$$
\D_\N^* = \D_1^* \ast_\gamma \D_{\N+1}.
$$
Therefore we get an integral presentation
\be
\begin{aligned}
&\mathcal{ML}(\mathcal{E}_{\D_\N^*}; s_1, \ldots, s_k,z)=\frac{1}{4\pi i}\int_{-i\infty}^{+i\infty} \mathcal{L}(\mathcal{E}_{\D_1^*}; s, z) \mathcal{M}(\mathcal{E}_{\D_{\N+1}}; -s,  s_1, \ldots, s_\N) ds.
\end{aligned}
\ee
 Substituting (\ref{27}) and  (\ref{eq:d1compx}) we get a Barnes integral presentation given by a $(\N-1)-$fold  integral. 
\vskip 2mm

There is another integral presentation for the same function using the double Laplace transform of  formula (\ref{Mei}) from the following Proposition. The particular case for $s_i=\frac{1}{2}$ was computed in \cite{HT25} using the similar method.

 \begin{figure}[ht]
	\centering
	\includegraphics[scale=0.8]{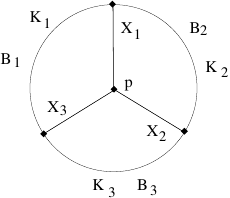}
	\small
	\caption{Triangulated punctured disc $\D^*_{3}$.}
	\label{figure:D13}
\end{figure}
\bp \la{Meij}The Mellin transform of the exponential volume of the crown with $\N$ cusps $\D_\N^*$ is 
\be \la{Mei}
\begin{aligned}
&\mathcal{M}(\mathcal{E}_{\D_\N^*}(e^l),s_1,...,s_\N) =\frac{1}{\pi i}\int_{-i\infty}^{+i\infty} \prod_{i=1}^\N\Gamma(s_{i+1}-t) \Gamma(s_i+t) e^{lt} dt
\\&=2\MeijerG{\N}{\N}{\N}{\N}{1-s_1,\ldots,1-s_\N}{s_1,\ldots,s_\N}{e^l}.
\end{aligned}
\ee 
Here $s_{i+\N}=s_i$ and the last expression is Meijer's function.  \ep

\bpr Take an ideal triangulation of $\D_\N^*$ by the edges sharing  the puncture as in Figure \ref{figure:D13}. Then we have the  cluster Poisson coordinates  $\X_1, ..., \X_\N,\B_1,...,\B_\N$ 
 for the moduli space ${\cal M}_{\D_\N^*}$. In addition to this, there are  $\K$-coordinates $\K_1,...,\K_\N$, and the local potentials $W_1, ..., W_\N$,  related as follows:  
\be
\label{eq:crownkb}
\K_{i+1}=\B_i \X_i \B_{i+1};\;\;\; W_i=\B_i+\B_i \X_i=\B_i+\frac{\K_{i+1}}{\B_{i+1}};\;\;\;\L= \prod_{i=1}^\N \X_i.
\ee
Let $x_i=\log \X_i$ and $l = \log \L$. 
Let us compute $\mathcal{M}(\mathcal{E}_{\D_\N^*}(e^l),s_1,...,s_\N)$. It is handy to consider the hyperplane 
 $\H_l=\{(x_1,...,x_\N)\;|\;x_1+...+x_\N=l\}$ with the induced cluster volume form $2 dx_1...dx_{\N-1}$. Then
\be
\begin{aligned}
&\mathcal{M}(\mathcal{E}_{\D_\N^*}(e^l),s_1,...,s_\N)
\\&:=\int_{(\mathbb{R}_{>0})^\N} \mathcal{E}_{D_\N^*}(\K,\L) \prod_{i=1}^\N \K_i^{s_i} d\log \K_i
\\&=2\int_{\H_l} d x_1 ... dx_{\N-1}\int_{(\mathbb{R}_{>0})^\N} \prod_{i=1}^\N  \K_i^{s_i} \exp(-\B_i(1+e^{x_i}))  d \log \B_i 
\\&=2\int_{\H_l}\prod_{i=1}^\N \left(e^{s_{i+1} x_i}\int_0^\infty \B_i^{s_i+s_{i+1}}e^{-\B_i(1+e^{x_i})}\frac{d \B_i}{\B_i}\right) dx_1...d x_{\N-1}
\\&=2\int_{\H_l}\prod_{i=1}^\N \frac{e^{s_{i+1} x_i} \Gamma(s_i+s_{i+1})}{(1+e^{x_i})^{s_i+s_{i+1}}} dx_1...d x_{\N-1}.
\end{aligned}
\ee 
\bl We have\footnote{Assuming first ${\rm Re}(s_{i+1}-t_i)>0$ and ${\rm Re}(s_{i}+t_i)>0$, and then using the analytic continuation.} 
\be \la{52}
\begin{aligned}
&\int_{-\infty}^{+\infty}\frac{e^{s_{i+1} x_i} \Gamma(s_i+s_{i+1})}{(1+e^{x_i})^{s_i+s_{i+1}}} e^{-t_i x_i} d x_i= \Gamma(s_{i+1}-t_i) \Gamma(s_i+t_i).
\end{aligned}
\ee \el
\bpr The substitution $u=e^{x_i}$ lead us quickly to the Euler $\B-$function integral: 
$$
{\rm LHS}(\ref{52}) = \Gamma(s_i+s_{i+1})\int_0^\infty \frac{u^{s_{i+1}-t_i-1}}{(1+u)^{s_i+s_{i+1}}}du =  \Gamma(s_i+s_{i+1})\frac{\Gamma(s_{i+1}-t_i) \Gamma(s_i+t_i)}{\Gamma(s_{i+1}-t_i + s_i+t_i)}. 
$$
\epr
Therefore
\be
\begin{aligned}
&\int_{-\infty}^{+\infty}\mathcal{M}(\mathcal{E}_{\D_\N^*}(e^l),s_1,...,s_\N) e^{-tl} d l = 2\prod_{i=1}^\N\Gamma(s_{i+1}-t) \Gamma(s_i+t).
\end{aligned}
\ee 
Hence
\be
\begin{aligned}
&\mathcal{M}(\mathcal{E}_{\D_\N^*}(e^l),s_1,...,s_\N) =\frac{1}{\pi i}\int_{-i\infty}^{+i\infty} \prod_{i=1}^\N\Gamma(s_{i+1}-t) \Gamma(s_i+t) e^{lt} dt.
\end{aligned}
\ee 
We arrive at formula (\ref{Mei}).
\epr

\subsection{Cylinder $\A_{1,1}$}
\begin{figure}[ht]
	\centering
	\includegraphics[scale=0.3]{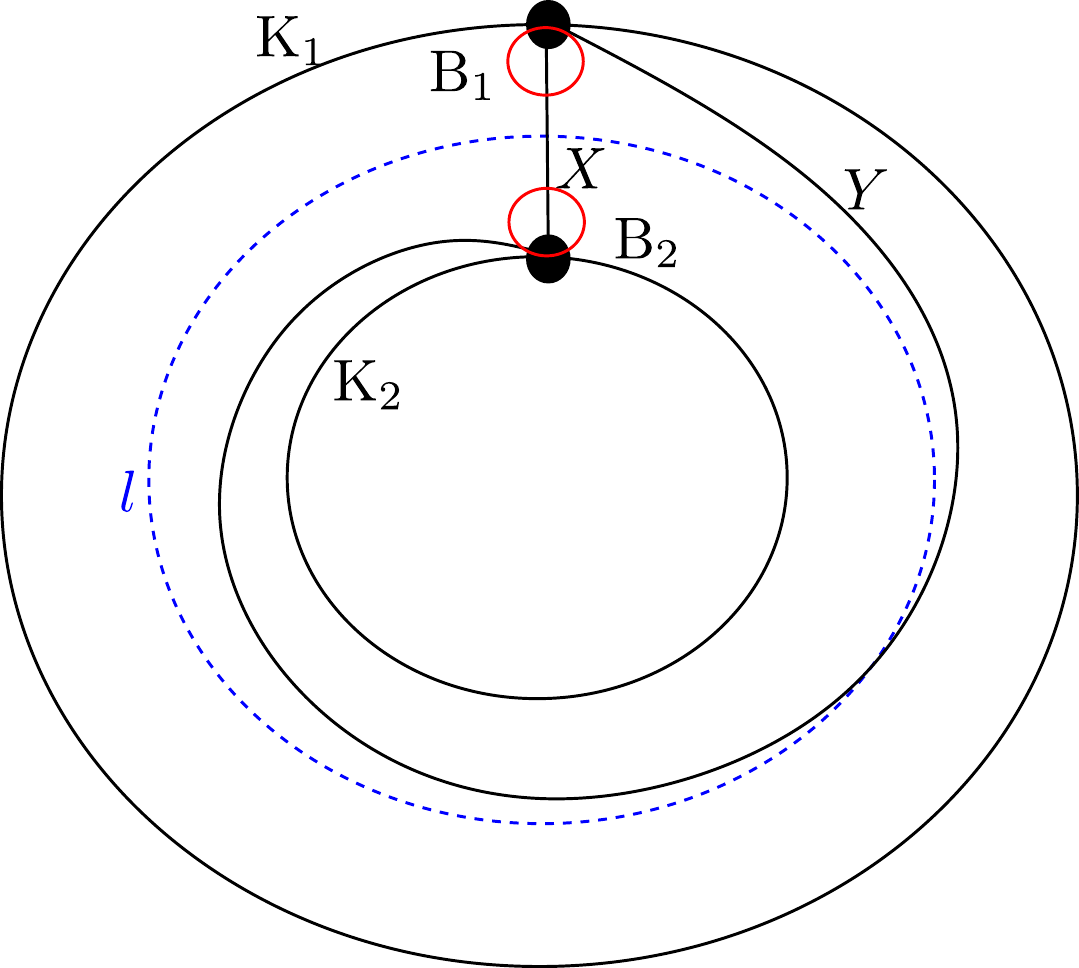}
	\small
	\caption{The cylinder.}
	\label{figure:A11}
\end{figure}
Consider an ideal hyperbolic structure on a cylinder $\A_{1,1}$ with one marked point on each  of the two boundary circles. As in Figure \ref{figure:A11}, the two boundary arcs has the $\K$-coordinates $\K_1$, $\K_2$. Given an ideal triangulation, let $\X$ and $\Y$ be the two $\mathcal{X}$-coordinates on the ideal edges. Then 
$\K_1=\B_1^2 \X \Y, \;\;\; \K_2=\B_2^2 \X \Y.$

\bp One has 
\[\mathcal{M}(\mathcal{E}_{\A_{1,1}}; s_1,s_2)=\frac{\Gamma(2s_1)\Gamma(2s_2)}{\pi i \Gamma(2s_1+2s_2)}  \int_{-i \infty}^{+i \infty} \frac{\Gamma(t)\Gamma(2s_1+2s_2-t)\Gamma(s_1+s_2-t)^{2}}{\Gamma(1+s_1+s_2-t)^{2}} dt. \]
The integral  is along a vertical line  separating the poles of $\Gamma(2s_1+2s_2-t)$ and $\Gamma(s_1+s_2-t)$ from those of $\Gamma(t)$.\ep

\bpr The unique simple closed geodesic has length $l$.
By \cite[Section 7.1.4]{GS24}, 
 the exponential volume
\[\mathcal{E}_{\A_{1,1}}(\K_1,\K_2)
=\frac{1}{2}\int_{0}^{+\infty} e^{-(\sqrt{\K_1}+\sqrt{\K_2})(e^{\frac{l}{2}}+e^{-\frac{l}{2}})} l d l.\]
Its Mellin transform 
\be
\begin{aligned}
&\mathcal{M}(\mathcal{E}_{\A_{1,1}}; s_1,s_2)=\frac{1}{2} \int_{0}^{+\infty}\int_{0}^{+\infty}\int_{0}^{+\infty} e^{-(\sqrt{\K_1}+\sqrt{\K_2})(e^{\frac{l}{2}}+e^{-\frac{l}{2}})} l d l \cdot \K_1^{s_1} d \log\K_1 \cdot \K_2^{s_2} d \log \K_2
\\&= 2\Gamma(2s_1)\Gamma(2s_2) \int_0^{+\infty}(e^{\frac{l}{2}}+e^{-\frac{l}{2}})^{-2 s_1-2 s_2} l   dl.
\end{aligned}
\ee

\bl 
\label{lemma:Mn}
For any non-negative integer $n$,
\be \la{Mn}
\rM_n(s):=\int_0^{+\infty}(e^{\frac{x}{2}}+e^{-\frac{x}{2}})^{-2s} x^{n}   dx=\frac{n!}{2\pi i\Gamma(2s)}\int_{-i \infty}^{+i \infty} \frac{\Gamma(t)\Gamma(2s-t)\Gamma(s-t)^{n+1}}{\Gamma(1+s-t)^{n+1}} dt.
\ee
The  integration  is over a vertical line  separating the poles of $\Gamma(2s-t)$ and $\Gamma(s-t)$ from those of $\Gamma(t)$.
\el
\bpr
Since 
\[(e^{\frac{x}{2}}+e^{-\frac{x}{2}})^{-2s}=((1+e^{-x})e^{\frac{x}{2}})^{-2s}=e^{-xs}\sum_{k=0}^{+\infty}\frac{(2s)_k}{k!}(-e^{-x})^k,\]
we have
\be \la{15}
\begin{aligned}
&\rM_n(s)= \int_0^{+\infty} e^{-xs}\sum_{k=0}^{+\infty}\frac{(2s)_k}{k!}(-e^{-x})^k  x^{n}   dx
\\&= \sum_{k=0}^{+\infty} \frac{(2s)_k (-1)^k}{k!} \int_0^{+\infty} e^{-x(s+k)} x^{n}   dx
\\&= \sum_{k=0}^{+\infty} \frac{(2s)_k (-1)^k}{k!} \frac{n!}{(s+k)^{n+1}}
\\&= n!\sum_{k=0}^{+\infty} \frac{(2s)_k (-1)^k}{k!(s+k)^{n+1}} 
\\&= \frac{n!}{s^{n+1}}\sum_{k=0}^{+\infty} \frac{(2s)_k (s)_k^{n+1}}{(s+1)_k^{n+1}}\frac{(-1)^k}{k!}
\\&= \frac{n!}{s^{n+1}} \cdot {}_{n+2}F_{n+1}(\{2s,s,\cdots,s\}, \{s+1,\cdots,s+1\},-1)
\\&= n! \left(\frac{\Gamma(s)}{\Gamma(s+1)}\right)^{n+1} \frac{\Gamma(s+1)^{n+1}}{\Gamma(s)^{n+1}\Gamma(2s)} \cdot \frac{1}{2\pi i} \int_{-i \infty}^{+i \infty} \frac{\Gamma(t)\Gamma(2s-t)\Gamma(s-t)^{n+1}}{\Gamma(1+s-t)^{n+1}} dt
\\&=\frac{n!}{2\pi i \Gamma(2s)} \int_{-i \infty}^{+i \infty} \frac{\Gamma(t)\Gamma(2s-t)\Gamma(s-t)^{n+1}}{\Gamma(1+s-t)^{n+1}} dt.
\end{aligned}
\ee
\epr

By Lemma \ref{lemma:Mn}, we get
\[\mathcal{M}(\mathcal{E}_{\A_{1,1}}; s_1,s_2)=2\Gamma(2s_1)\Gamma(2s_2) \rM_1=\frac{\Gamma(2s_1)\Gamma(2s_2)}{\pi i \Gamma(2s_1+2s_2)}  \int_{-i \infty}^{+i \infty} \frac{\Gamma(t)\Gamma(2s_1+2s_2-t)\Gamma(s_1+s_2-t)^{2}}{\Gamma(1+s_1+s_2-t)^{2}} dt, \]
where the integration  is along a vertical line  separating the poles of $\Gamma(2s_1+2s_2-t)$ and $\Gamma(s_1+s_2-t)$ from those of $\Gamma(t)$.
\epr

\bl   The Taylor series in $z$ of the exponential volume $\mathcal{ML}(\mathcal{E}_{\D_1^*}; s,-z)$ is given by \label{cor:de}
\be
\begin{aligned}
&\mathcal{ML}(\mathcal{E}_{\D_1^*}; s,-z)
=\sum_{n=0}^{+\infty} \frac{z^n}{\pi i } \int_{-i \infty}^{+i \infty} \frac{\Gamma(t)\Gamma(2s-t)\Gamma(s-t)^{n+1}}{\Gamma(1+s-t)^{n+1}} dt.
\end{aligned}
\ee
\el
\bpr
By the second line of Equation \eqref{eq:d1comp}, we obtain
\be
\begin{aligned}
&\mathcal{ML}(\mathcal{E}_{\D_1^*}; s,-z)=2\Gamma(2s)\int_0^{+\infty} (e^{\frac{l}{2}}+e^{-\frac{l}{2}})^{-2s} e^{lz} dl
\\&=2\Gamma(2s)\sum_{k=0}^{+\infty} \int_0^{+\infty} (e^{\frac{l}{2}}+e^{-\frac{l}{2}})^{-2s} l^k dl \frac{ z^k}{k!}
\\&=2\Gamma(2s)\sum_{k=0}^{+\infty} \rM_k(s) \frac{ z^k}{k!}.
\end{aligned}
\ee
It remains to use formula (\ref{Mn}) from Lemma \ref{lemma:Mn}. \epr

\subsection{One-cusped pair of pants  $S_{0,3;1}$}
\begin{figure}[ht]
	\centering
	\includegraphics[scale=0.3]{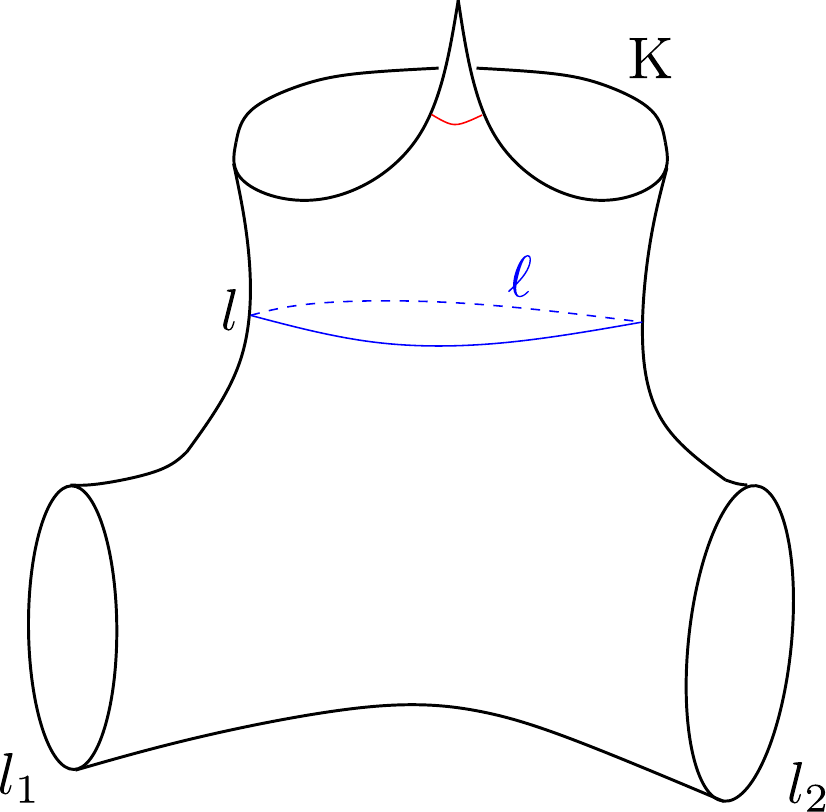}
	\small
	\caption{The pair of pants with one cusp.}
	\label{figure:s031}
\end{figure}
Consider an ideal hyperbolic structure on a pair of pants $S_{0,3;1}$ with one cusp, see Figure \ref{figure:s031}. The neck geodesic $\ell$ cut it into a crown $\D_1^*$ with one cusp and a pair of pants $S_{0,3}$ with boundary lengths $l$, $l_1$,  $l_2$. There is a unique $\K$-coordinate $\K$ on the boundary of $\D_1^*$.

\bp
\be
\begin{aligned}
&\mathcal{ML}(\mathcal{E}_{S_{0,3;1}}; s,z_1,z_2)=  \frac{1}{\pi i z_1 z_2 } \int_{-i \infty}^{+i \infty} \frac{\Gamma(t)\Gamma(2s-t)\Gamma(s-t)^{2}}{\Gamma(1+s-t)^{2}} dt,
\end{aligned}
\ee
The integral  is along a vertical  line separating the poles of $\Gamma(2s-t)$ and $\Gamma(s-t)$ from those of $\Gamma(t)$.\ep

\bpr Recall the exponential volumes
\be \la{pig}
\begin{split}
&\mathcal{E}_{\D_1^*}(\K,l)=e^{-\sqrt{\K}(e^{\frac{l}{2}}+e^{-\frac{l}{2}})},
\\&
\mathcal{E}_{S_{0,3}}(l,l_1,l_2)=1.\\
\end{split}
\ee
Then by the neck recursion formula \cite[Theorem 1.7]{GS24}, we have
\be
\mathcal{E}_{S_{0,3;1}}(\K,l_1,l_2)=\int_0^{+\infty} e^{-\sqrt{\K}(e^{\frac{l}{2}}+e^{-\frac{l}{2}})} l d l.
\ee
Then its Mellin--Laplace transform is
\be
\begin{aligned}
&\mathcal{ML}(\mathcal{E}_{S_{0,3;1}};s,z_1,z_2)=\frac{1}{z_1 z_2}\int_0^{+\infty}\int_{0}^{+\infty}e^{-\sqrt{\K}(e^{\frac{l}{2}}+e^{-\frac{l}{2}})} l dl \K^{s} d\log \K
\\&= \frac{2\Gamma(2s)}{z_1 z_2}  \int_0^{+\infty}(e^{\frac{l}{2}}+e^{-\frac{l}{2}})^{-2s} l  dl
\\&\stackrel{(\ref{Mn})}{=}  \frac{1}{\pi i z_1 z_2 } \int_{-i \infty}^{+i \infty} \frac{\Gamma(t)\Gamma(2s-t)\Gamma(s-t)^{2}}{\Gamma(1+s-t)^{2}} dt.
\end{aligned}
\ee
\epr
\subsection{Once-punctured torus with a single cusp $S_{1,1;1}$}
\begin{figure}[ht]
	\centering
	\includegraphics[scale=0.3]{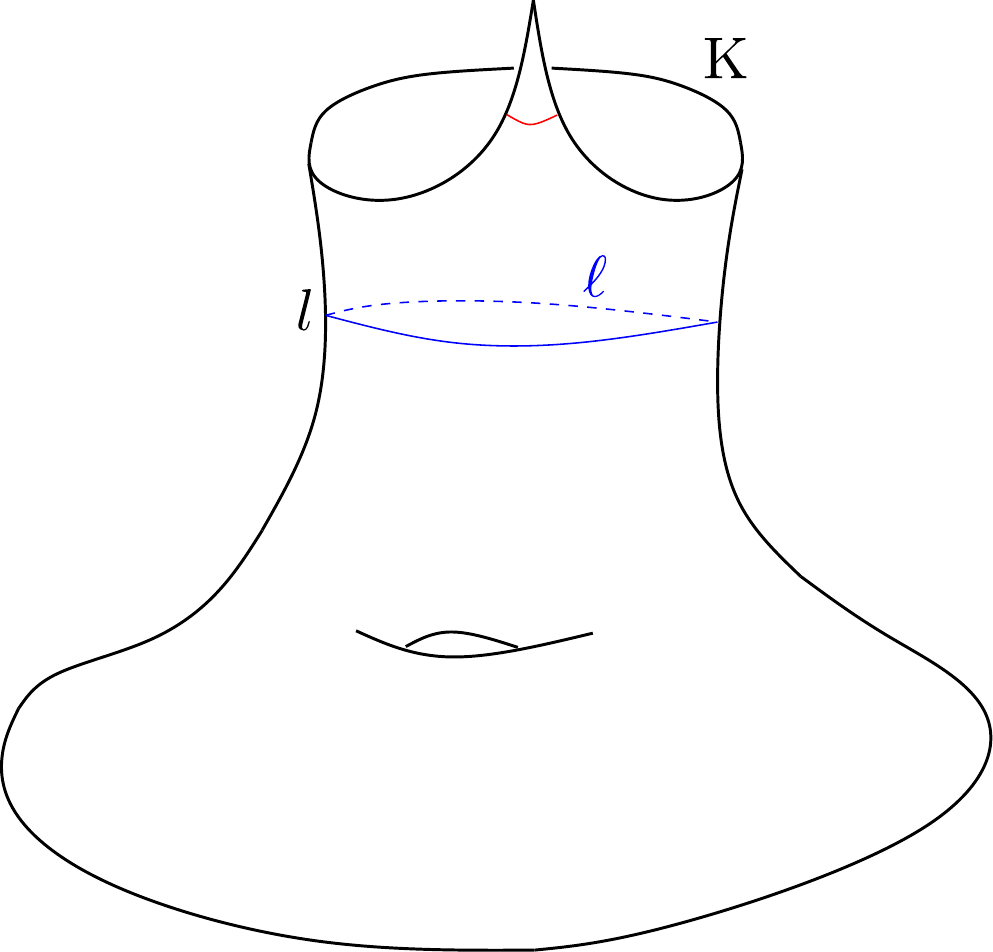}
	\small
	\caption{The once punctured torus with one cusp.}
	\label{figure:s111}
\end{figure}
Take an ideal hyperbolic structure on a once punctured torus $S_{1,1;1}$ with one cusp. The neck geodesic $\ell$ cuts it into a one-cusped crown $\D_1^*$ and a once punctured torus $S_{1,1}$ with boundary lengths $l$, see Figure \ref{figure:s111}. There is a unique $\K$-coordinate $\K$ on the boundary of $\D_1^*$. 
\bp
\be \la{33}
\begin{aligned}
&\mathcal{M}(\mathcal{E}_{S_{1,1;1}};s) =  \frac{\pi}{6i} \int_{-i \infty}^{+i \infty} \frac{\Gamma(t)\Gamma(2s-t)\Gamma(s-t)^{2}}{\Gamma(1+s-t)^{2}} dt  +\frac{1}{4\pi i}  \int_{-i \infty}^{+i \infty} \frac{\Gamma(t)\Gamma(2s-t)\Gamma(s-t)^{4}}{\Gamma(1+s-t)^{4}} dt.
\end{aligned}
\ee Integration contours are  vertical lines  separating poles of $\Gamma(2s-t)$ and $\Gamma(s-t)$ from those of $\Gamma(t)$.\
\ep

\bpr The exponential volumes are $\mathcal{E}_{\D_1^*}(\K,l)$, see (\ref{pig}), and 
\be
\mathcal{E}_{S_{1,1}}(l)=\frac{\pi^2}{6}+\frac{l^2}{24}.
\ee
Then by the neck recursion formula \cite[Theorem 1.7]{GS24}, we have
\be
\mathcal{E}_{S_{1,1;1}}(\K,l_1,l_2)=\frac{\pi^2}{6}\int_0^{+\infty} e^{-\sqrt{\K}(e^{\frac{l}{2}}+e^{-\frac{l}{2}})} l d l+\frac{1}{24}\int_0^{+\infty} e^{-\sqrt{\K}(e^{\frac{l}{2}}+e^{-\frac{l}{2}})} l^3 d l.
\ee
Then its Mellin transform is
\be
\begin{aligned}
&\mathcal{M}(\mathcal{E}_{S_{1,1;1}};s)
\\&=\frac{\pi^2}{6}\int_0^{+\infty} \int_0^{+\infty}  e^{-\sqrt{\K}(e^{\frac{l}{2}}+e^{-\frac{l}{2}})} l d l \K^{s} d\log \K +\frac{1}{24}\int_0^{+\infty}  \int_0^{+\infty} e^{-\sqrt{\K}(e^{\frac{l}{2}}+e^{-\frac{l}{2}})} l^3 d l \K^{s} d\log \K
\\&= \frac{2 \Gamma(2s)\pi^2}{6} \int_0^{+\infty}  (e^{\frac{l}{2}}+e^{-\frac{l}{2}})^{-2s} l d l  +\frac{2 \Gamma(2s)}{24}  \int_0^{+\infty} (e^{\frac{l}{2}}+e^{-\frac{l}{2}})^{-2s} l^3 d l .
\end{aligned}
\ee
Using formula (\ref{Mn}) we get (\ref{33}). \epr

\section{Barnes integrals for  Mellin--Laplace transformed exponential volumes}
As one can see in the examples of Section \ref{sectionex}, the  Mellin--Laplace transforms of exponential volumes of certain decorated surfaces can be expressed 
as  Barnes integrals. Let us prove this fact in general.

\bt
\label{thmint}
For any decorated surface $\bS=S_{g,m;n_1,\cdots,n_r}$, the Mellin--Laplace transform $\mathcal{ML}(\mathcal{E}_\bS;\rm{S},\rm{Z})$ (Definition \ref{defn:LM}) of exponential volume   
has an integral presentation as 
\be
\sum_{0\leq d_1+\cdots d_m\leq 3g-3+m} c_{d_1,\cdots,d_m}  \frac{\F_{d_1,\cdots,d_m}(\rm{S})}{\prod_{i=r+1}^m z_{i-r}^{2d_i+1}}
\ee
where 
\begin{itemize}
\item the sum is over all non-negative integers $d_1,\cdots,d_m$ such that $0\leq d_1+\cdots d_m\leq 3g-3+m$, 
\item the $c_{d_1,\cdots,d_m}$ are non-negative constants,  
\item the $\F_{d_1,\cdots,d_m}(\rm{S})$ can be written as an $k$-multiple integral of products and ratios of gamma functions:
\be
\label{eq:g}
\F_{d_1,\cdots,d_m}(\rm{S})= \left(\frac{1}{2\pi i}\right)^k \int_{-i \infty}^{+i \infty}...\int_{-i \infty}^{+i \infty} \frac{\Gamma \cdots \Gamma}{\Gamma \cdots \Gamma}.
\ee
\end{itemize}
\et
We say a function $g(\rm{S},\rm{Z})$ has a {\em Barnes integral presentation} if $g$ can be written as a linear summation of $h_i(\rm{Z}) f_i(\rm{S})$ where $h_i(\rm{Z})$ is a monomial and $f_i(\rm{S})$ is a multiple integral of products and ratios of gamma functions as in Equation (\ref{eq:g}). The above theorem says that 
$$
\mbox{the Mellin--Laplace transform $\mathcal{ML}(\mathcal{E}_\bS;\rm{S},\rm{Z})$ has a {\em Barnes integral presentation}.}
$$
\bpr
By the computations in Section \ref{sectionex}, Theorem \ref{thmint} is true for $\Delta$, $\D_1^*$, $\A_{1,1}$. Suppose Theorem \ref{thmint} is true for $\bS$, let $\bS'=\bS \cup \Delta$ obtained by gluing along an ideal edge with the $\K$-coordinate $\K$. The Mellin--Laplace transforms are $\mathcal{ML}(\mathcal{E}_{\bS};{\rm{S}},\rm{Z})$, $\mathcal{M}(\mathcal{E}_{\Delta};s,s_1,s_2)$ and $\mathcal{ML}(\mathcal{E}_{\bS'};{\rm{S}}\backslash \{s\},s_1,s_2,\rm{Z})$. Then by Parseval's formula
\be
\label{eq:par}
\begin{aligned}
&\mathcal{ML}(\mathcal{E}_{\bS'};{\rm{S}}\backslash \{s\},s_1,s_2,\rm{Z})
\\&=\int_0^{+\infty} \mathcal{ML}(\mathcal{E}_{\bS}(\K);{\rm{S}}\backslash \{s\},{\rm{Z}})\cdot \mathcal{M}(\mathcal{E}_{\Delta}(\K);s_1,s_2) d\log \K
\\&= \frac{1}{2\pi i} \int_{-i \infty}^{+i \infty}  \mathcal{ML}(\mathcal{E}_{\bS};{\rm{S}}\backslash \{s\},s,{\rm{Z}})\cdot \mathcal{M}(\mathcal{E}_{\Delta};-s,s_1,s_2) ds.
\end{aligned}
\ee
Since Theorem \ref{thmint} is true for $\mathcal{ML}(\mathcal{E}_{\bS};\rm{S},\rm{Z})$ and $\mathcal{M}(\mathcal{E}_{\Delta};-s,s_1,s_2)$ is products of gamma functions, Theorem \ref{thmint} is true for $\bS'$. 
\begin{remark}
Similarly, if $\bS'=\bS \cup \Delta$ obtained by gluing along two ideal edges with the $\K$-coordinates $\K_a$, $\K_b$, we can use Parseval's formula as in Equation (\ref{eq:par}) twice to show that the Mellin--Laplace transform for $\bS'$ has a Barnes integral representation.
\end{remark}
\begin{figure}[ht]
	\centering
	\includegraphics[scale=0.3]{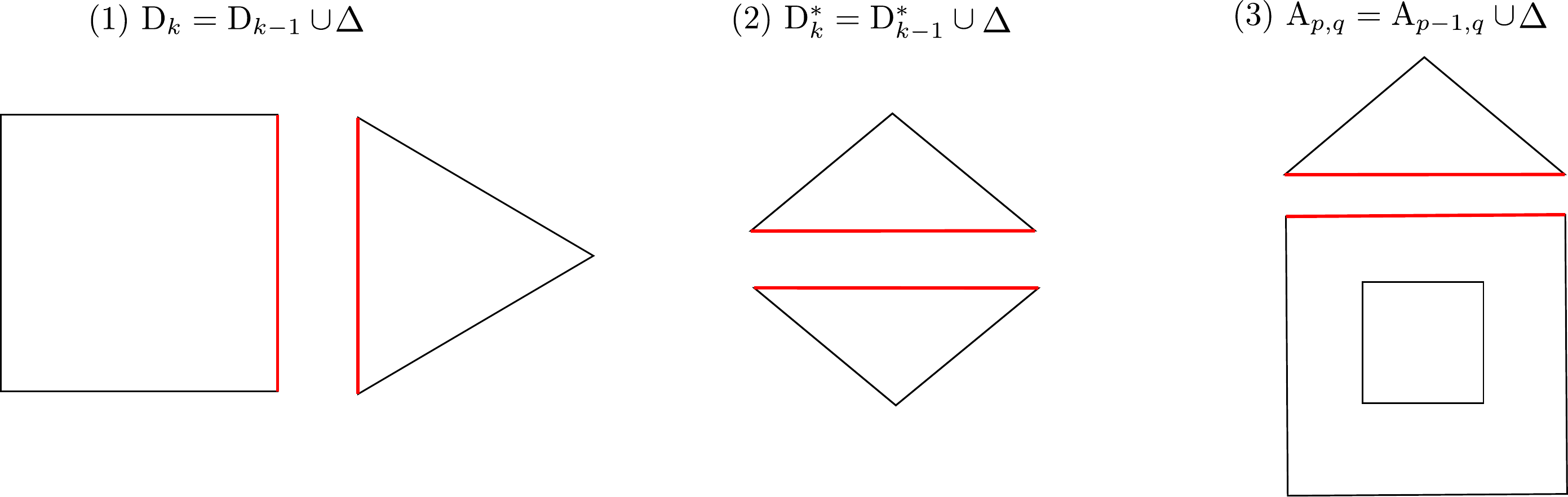}
	\small
	\caption{Gluing triangles.}
	\label{figure:glued3}
\end{figure}
As in Figure \ref{figure:glued3}, since 
$$
\D_k=\D_{k-1} \cup \Delta, \ \ \ \D_k^*=\D_{k-1}^* \cup \Delta, \ \ \ \A_{p,q}=\A_{p-1,q}\cup \Delta=\A_{p,q-1}\cup \Delta,
$$
 Theorem \ref{thmint} is true for $\D_k$, $\D_k^*$ and $\A_{p,q}$.

For any other decorated surface $\bS=S_{g,m;n_1,\cdots,n_r}$, we have $2-2g-m<0$. Since
\[S_{g,m;n_1,\cdots,n_r}=S_{g,m}\cup \D_{n_1}^*\cup \cdots \cup \D_{n_r}^*,\]
by cutting off ideal triangles, it is enough to prove Theorem \ref{thmint} for $\bS=S_{g,m}\cup \cup_{i=1}^r \D_{1}^*$. Given an ideal hyperbolic structure on the decorated surface $\bS=S_{g,m}\cup \cup_{i=1}^r \D_{1}^*$, the boundary parameters are $\K$-coordinates $\K_1,\cdots,\K_{\bf r}$ and geodesic boundary circle lengths $l_{r+1},\cdots,l_m$. After the Mellin--Laplace transform of exponential volume, the $\K$-coordinates $\K_1,\cdots,\K_{\bf r}$ change into $s_1,\cdots,s_r$, and the lengths $l_{r+1},\cdots,l_m$ change into $z_1,\cdots,z_{m-r}$. Then by neck recursion formula in \cite[Theorem 1.7]{GS24}, we obtain
\be
\begin{aligned}
\label{eq:sgmd1r}
&\mathcal{ML}(\mathcal{E}_{\bS};s_1,\cdots,s_r,z_1,\cdots,z_{m-r})
\\&=\int_0^{+\infty}\cdots \int_0^{+\infty} \mathcal{ML}(\mathcal{E}_{S_{g,m}}(l_1,\cdots,l_r);z_1,\cdots,z_{m-r}) \prod_{i=1}^r \mathcal{M}(\mathcal{E}_{\D_1^*}(l_i);s_i) l_i d l_i.
\end{aligned}
\ee
Here
\be
\label{eq:d1l}
\begin{aligned}
&\mathcal{M}(\mathcal{E}_{\D_1^*}(l_i);s_i)=\int_{0}^{+\infty} e^{-\sqrt{\K_i}(e^{\frac{l_i}{2}}+e^{-\frac{l_i}{2}})} \K_i^{s_i} d \log \K_i
\\&=2\Gamma(2s_i) (e^{\frac{l_i}{2}}+e^{-\frac{l_i}{2}})^{-2s_i}.
\end{aligned}
\ee
By Mirzakhani's theorem \cite[Theorem 1.1]{Mir07b}:
\be
\label{eq:M07}
\mathcal{E}_{S_{g,m}}(l_1,\cdots,l_m) = \sum_{0\leq d_1+ \cdots + d_m \leq d} {\cal V}_{d_1, \cdots, d_m}l_1^{2d_1} \ldots  l_m^{2d_m},
\ee
Here $d_i\geq 0$ are integers and ${\cal V}_{d_1, \cdots, d_m} \in \pi^{6g-6+2m-d_1-\cdots-d_m}{\Bbb Q_{>0}}$. Thus 
\be
\label{eq:ml}
\mathcal{ML}(\mathcal{E}_{S_{g,m}}(l_1,\cdots,l_r);z_1,\cdots,z_{m-r})=\sum_{0\leq d_1+ \cdots + d_m \leq d} {\cal V}_{d_1, \cdots, d_m}\prod_{i=1}^r l_i^{2d_i} \prod_{i=r+1}^m \frac{(2d_i)!}{z_{i-r}^{2d_i+1}}.
\ee
Then plug Equation (\ref{eq:d1l}) and Equation (\ref{eq:ml}) into Equation (\ref{eq:sgmd1r}), following Lemma \ref{lemma:Mn}, we get 
\be
\begin{aligned}
\label{eq:sgmd1r}
&\mathcal{ML}(\mathcal{E}_{\bS};s_1,\cdots,s_r,z_1,\cdots,z_{m-r})
\\&=\int \sum_{0\leq d_1+ \cdots + d_m \leq d} {\cal V}_{d_1, \cdots, d_m}\prod_{i=1}^r l_i^{2d_i} \prod_{i=r+1}^m \frac{(2d_i)!}{z_{i-r}^{2d_i+1}} \prod_{i=1}^r 2\Gamma(2s_i) (e^{\frac{l_i}{2}}+e^{-\frac{l_i}{2}})^{-2s_i} l_i d l_i
\\&=\sum_{0\leq d_1+ \cdots + d_m \leq d}  2^r {\cal V}_{d_1, \cdots, d_m}\prod_{i=r+1}^m \frac{(2d_i)!}{z_{i-r}^{2d_i+1}} \prod_{i=1}^r \Gamma(2s_i) \prod_{i=1}^r \int_0^{+\infty} (e^{\frac{l_i}{2}}+e^{-\frac{l_i}{2}})^{-2s_i} l_i^{2d_i+1} d l_i
\\&\overset{\ref{lemma:Mn}}{=}\sum_{0\leq d_1+ \cdots + d_m \leq d}   \frac{2^r {\cal V}_{d_1, \cdots, d_m} \prod_{i=1}^r (2d_i+1)!\prod_{i=r+1}^m (2d_i)!}{\prod_{i=r+1}^m z_{i-r}^{2d_i+1}}  \prod_{i=1}^r
\frac{1}{2\pi i}\int_{-i \infty}^{+i \infty} \frac{\Gamma(t)\Gamma(2s_i-t)\Gamma(s_i-t)^{2d_i+2}}{\Gamma(1+s_i-t)^{2 d_i+2}} dt.
\end{aligned}
\ee
Here the integration contour  for the $i-$th integral in the product  is a vertical line separating the poles of $\Gamma(2s_i-t)$ and $\Gamma(s_i-t)$ from those of $\Gamma(t)$.
So Theorem \ref{thmint} is true for the  $\mathcal{ML}(\mathcal{E}_{\bS};s_1,\cdots,s_r,z_1,\cdots,z_{m-r})$. We conclude that Theorem \ref{thmint} is true for any decorated surface $\bS$.
\epr

 \section{Colored surfaces and weighted exponential volumes} \la{SEC4}

\subsection{Weighted exponential volumes of moduli spaces} \la{SEC4.1}
 
 Recall  from Definitions \ref{DEF2.7} $\&$  \ref{D1.10}  the colored surface $\bS_{\rm{c}}=(\bS; {\cal V}_{\bf r})$ and the weighted colored surface  
 $$
 \bS_{\rm{c}; \omega}=(\bS_{\rm c}; \omega). 
 $$
 
 The weight  is a complex valued  function on the set of blue boundary intervals 
$
\omega: {\cal I}_{\bf b} \lra \C.
$ 

As we explained in Section \ref{sec1.2.2}, it induces a function $t$ on the set of blue cusps: 
 \be \la{v}
 t_b:= \omega_e+\omega_{e'}, \ \ \ \forall b\in {\cal V}_{\bf b}.
\ee
Therefore a {weighted colored surface}  $(\bS_{\rm c}; \omega)$ 
gives rise to the following: 

\begin{itemize}

\item 
The {\em weight factor} $\mathcal{Q}_{{\bS}_{{\rm c}}}(\omega)$  on the Teichm\"uller space ${\cal T}_{\bS_{\rm{c}}}$,  given by
\be \la{G}
\mathcal{Q}_{\bS_{{\rm{c}}}}(\omega):=\prod_{b \in {\cal V}_{\bf b}} \Gamma(t_b)\frac{\prod_{e\in \mathcal{I}_{\bf b}}\K_{e}^{\omega_e}}{ \prod_{b\in \mathcal{V}_{\bf b}} W_b^{t_b}}.
\ee
\item 
 The {\em weighted exponential volume form} on the Teichm\"uller space ${\cal T}_{\bS_{\rm{c}}}$  given by
$$
{\Bbb E}_{\bS_{{\rm{c}}}}(\omega):= {\cal Q}_{\bS_{{\rm{c}}}}(\omega) e^{-W_{\bS_{\rm{c}}}}\Omega_{\bS_{\rm{c}}}.
$$
\end{itemize}

\bp \la{L5.1} Let  $(\bS_{\rm{c}}; \omega)$ be a weighted colored surface. Then

i) The weight factor ${\cal Q}_{\bS_{{\rm{c}}}}(\omega)$ does not depend on the choice of horocycles $h_b$ at the blue vertices $b$ on $\bS_{\rm{c}}$. So it is a well defined  function on the Teichm\"uller space ${\cal T}_{\bS_{\rm{c}}}$.

ii)  The weighted exponential volume form  on the space ${\cal T}_{\bS_{\rm{c}}}$  descends to the moduli space ${\cal M}_{\bS_{\rm{c}}}$. 
\ep

\bpr i) Let us  think about the horocycles $h_p$ as of vectors $v_p$ in the related ${\rm PGL}_2-$local system on $\bS$, see \cite{GS24} or Section \ref{Sect5.3.3}.\footnote{One can easily translate the arguments below to the hyperbolic geometry language, using horocycles rather than vectors.}  
There is a  map forgetting the  horocycles at blue marked points:
\be \la{pi}
\pi:\mathcal{M}_{\bS}\rightarrow  \mathcal{M}_{\bS_{\rm{c}}}.
\ee
Given any $X\in \mathcal{M}_{\bS_{\rm{c}}}$, the fiber $\pi^{-1}(X)$ is the principal homogeneous space for the  fiberwise action $\rho$ of the multiplicative group $(\R_{>0}^\times)^k$, defined as follows: 
an element $\alpha = (\alpha_1, ..., \alpha_k)\in (\R_{>0}^\times)^k$ acts  by 
\[\rho_\alpha: v_i\lms \alpha_i^{-\frac{1}{2}} v_i, \ \ \ i=1, ..., k,\]
for the decorations at blue marked points. 
Then we have 
\be \la{NOT}
\begin{aligned}
&W_{i, \alpha} := \rho_\alpha^* W_i = \alpha_i W_i, \ \ \ \rho_\alpha^*W_{\alpha}=\sum_{i=1}^k \alpha_iW_{i},
\\& \K_{ij;\alpha}:= \rho^*_\alpha\K_{ij} =\alpha_i \alpha_j \K_{ij},
\\& \lambda_{ij;\alpha}:= \rho_\alpha^*\lambda_{ij}  =(\alpha_i \alpha_j)^{-\frac{1}{2}} \lambda_{ij}. \\
\end{aligned}
\ee  
Let us rescale a vector $v_b$ at  blue marked point $b$ by $v_b\lms \alpha^{-1/2} v_b$. Then  the weight function  $\mathcal{Q}_{\bS_{\rm{c}}}(\omega)$ is rescaled by 
$
\alpha^{\omega_e + \omega_{e'} - t_b} =1.
$

ii) The function $\mathcal{Q}_{\bS_{\rm{c}}}(\omega)$ is a ${\rm Mod}(\bS)$-invariant function.
\epr

\paragraph{Proof of Lemma \ref{tLMEVL}.}
\label{rem:deform}
The deformations of $W_{i, \alpha}$ induce defomations of $\K_{ij;\alpha}$. Then Lemma \ref{tLMEVL} follows from \eqref{NOT}. 
Here is the simplest but typical  example of such a calculation. Consider the basic triangle  as in Figure \ref{figure:D3}.
The potential 
$
W=t_1\B_1+t_2\B_2+t_3 \B_3. 
$ 
The exponential volume
\[\mathcal{E}_{\Delta}(\K_1,\K_2,\K_3)= 
 e^{-(t_1\B_1+t_2\B_2+t_3\B_3)}.\]
Then the Mellin transform of exponential volume is calculated via $\Gamma-$functions as follows:
\be
\label{eq:D3x}
\begin{aligned}
&\mathcal{M}(\mathcal{E}_{\Delta}; t_1, t_2, t_3 | s_1,s_2,s_3)
\\&=\int_0^{+\infty}\int_0^{+\infty}\int_0^{+\infty}e^{-(t_1\B_1+t_2\B_2+t_3\B_3)} (\B_3 \B_1)^{s_1} (\B_1 \B_2)^{s_2} (\B_2 \B_3)^{s_3}  2 d \log \B_1 d \log \B_2 d \log \B_3
\\&=2 \int_0^{+\infty} e^{-t_1B_1} \B_1^{s_1+s_2} d \log \B_1 \int_0^{+\infty} e^{-t_2\B_2} \B_2^{s_2+s_3} d \log \B_2 \int_0^{+\infty} e^{-t_3\B_3} \B_3^{s_3+s_1} d \log \B_3
\\&=2 t_1^{-(s_1+s_2)} t_2^{-(s_2+s_3)} t_3^{-(s_3+s_1)} \Gamma(s_1+s_2) \Gamma(s_2+s_3) \Gamma(s_3+s_1).
\end{aligned}
\ee

\subsection{Weighted exponential volumes  via  Mellin transformed exponential volumes}

 Given a weighted colored surface  $(\bS_{\rm{c}}; \omega)$ with the  $\K$-parameters $\K_{\bf r}$ at the  red boundary intervals, recall the weighted exponential volume of the moduli space ${\cal M}_{\bS_{\rm{c}}}$ 
\be
{\cal E}_{\bS_{\rm{c}}}(\omega;\K_{\bf r}, \L):=\int_{{\cal M}_{\bS_{\rm{c}}}(\K_{\bf r}, \L)}{\Bbb E}_{\bS_{\rm{c}}}(\omega;\K_{\bf r}, \L). 
\ee

We enumerate by $\{1, ..., k\}$ the set ${\cal V}_{\bf b}$ of the blue marked boundary points and denote by $t_1, ..., t_k$ the values of the function $t_*$ at these points. The invariant differential form along the  fibers of the map $\pi$ in (\ref{pi}) is $d\log \alpha_1\wedge \ldots \wedge d\log \alpha_k$.  Note that  the number of the blue vertices may differ from the number of blue edges.

Let $\mathcal{M}_{\rm blue}(\cdot)$ be the Mellin transform of the exponential volume function ${\cal E}_\bS(\K, \L)$ along the $\K-$variables $\K_1, ..., \K_b$ at the blue edges. The resulting function is denoted by $\mathcal{M}_{\rm blue}(\mathcal{E}_{{\bS}}(\K_{\bf r}, \L); \omega_1, ..., \omega_b)$. 

\bt
\label{thm:fc}
Given a weighted colored surface  $(\bS_{\rm{c}}; \omega)$, the weighted exponential volume of the related moduli space ${\cal M}_{\bS_{\rm{c}}}(\K_{\bf r}, \L)$   is equal to the Mellin transform of the exponential volume of the moduli space 
${\cal M}_{\bS}(\K, \L)$ for the underlying decorated surface $\bS$, evaluated at the weight parameters $\omega_1, ..., \omega_b$:
\be
\label{eq:fc}
{\cal E}_{\bS_{\rm{c}}}(\omega_1, ..., \omega_b; \K_{\bf r}, \L)=\mathcal{M}_{\rm blue}(\mathcal{E}_{{\bS}}(\K_{\bf r}, \L); \omega_1, ..., \omega_b).
\ee
\et
\bpr
 Recall the notation (\ref{NOT}) for $\K_{e; \alpha}$. The following lemma plays the key role. 

\bl
\label{lem:fb}
\be \la{LEM}
\mathcal{Q}_{\bS_{\rm{c}}}(\omega) = 
\int_{(\mathbb{R}_{>0})^k}\prod_{e\in \mathcal{I}_{\bf b}} \K_{e;\alpha}^{\omega_e} \exp(-W_\alpha) \prod_{i=1}^{k}d\log \alpha_i.
\ee
\el
\bpr
We have
\be
\begin{aligned}
&\int_{(\mathbb{R}_{>0})^k} \prod_{e\in \mathcal{I}_{\bf b}} \K_{e;\alpha}^{\omega_e} \exp(-W_\alpha) \prod_{i=1}^{k}d\log \alpha_i
\\&=\int_{(\mathbb{R}_{>0})^k} \prod_{i=1}^k \alpha_i^{t_i} \prod_{e\in \mathcal{I}_{\bf b}} \K_{e}^{\omega_e} \exp(-\sum_{i=1}^k \alpha_i W_i) \prod_{i=1}^{k} \frac{d\alpha_i}{\alpha_i}
\\&= \prod_{e\in \mathcal{I}_{\bf b}} \K_{e}^{\omega_e}   \prod_{i=1}^{k} \int_{\mathbb{R}_{>0}} \alpha_i^{t_i} \exp(-\alpha_i W_i) \frac{d\alpha_i}{\alpha_i}
\\&=\frac{\prod_{i=1}^{k} \Gamma(t_i) \prod_{e\in \mathcal{I}_{\bf b}} \K_{e}^{\omega_e}}{ \prod_{i=1}^{k} W_i^{t_i}}
\\&= \mathcal{Q}_{\bS_{\rm{c}}}(\omega).
\end{aligned}
\ee
 \epr

Since the differential form along the  fiber is $d\log \alpha_1\wedge \ldots \wedge d\log \alpha_k$,  applying Lemma \ref{lem:fb} we have
\be \la{HF}
\int_{{\cal M}_{\bS_{\rm{c}}}(\K_{\bf r}, \L)}{\Bbb E}_{\bS_{\rm{c}}}(\omega;\K_{\bf r}, \L)=\int_{{\cal M}_{\bS}(\K_{\bf r}, \L)} \left(\prod_{e\in \mathcal{I}_{\bf b}} \K_{e}^{\omega_e}\right) {\Bbb E}_{\bS}(\K_{\bf r}, \L).
\ee
\epr
 
\begin{remark} \la{REM4.4}
{\rm By Theorem \ref{thm:fc}, the integral ${\cal E}_{\bS_{\rm{c}}}(\omega; \K_{\bf r}, \L)$ depends only on  the colors of boundary edges, assuming that the weights $\{\omega_e\}$ are the same.  One can have different colorings of the cusps providing the same colors of the edges, and therefore the same values of the weighted exponential volumes. 
Here is a typical example. 
Given a decorated surface $\bS$, 
let us  pick any subset  ${\cal V}_{\bf r}$ of the red cusps in  the set  ${\cal V}_\bS$, which is {\em sparse} in the following sense: 
$$
\mbox{\it No boundary interval on $\bS$ connects two red cusps, that is the cusps from the set ${\cal V}_{\bf r}$.}
$$
For example, the empty subset ${\cal V}_{\bf r} = \emptyset$ is sparse. All boundary intervals on any sparsely colored surface are blue. Therefore the weighted exponential volumes of 
 the moduli spaces related to  sparsely colored surfaces are the same, provided the weights $\{\omega_e\}$ are the same. Yet the moduli spaces are different. 
Even their dimensions are  different 
if the numbers of red cusps differ.}

\end{remark}
 
Applying the  Mellin--Laplace transform to  both sides of  \eqref{eq:fc}, we obtain
\bc 
\label{cor:mwcs}
Given a weighted colored surface $(\bS_{\rm{c}};\omega)$, the Mellin-Laplace transformed weighted exponential volume of the moduli space ${\cal M}_{\bS_{\rm c}}$ is equal to the 
Mellin-Laplace transformed  exponential volume of the moduli space ${\cal M}_{\bS}$, with the Mellin parameter $s_e$  for each blue boundary interval ${\rm I}_e$  on $\bS_{\rm c}$ specialized to $s_e:= \omega_e$: 
\be
\mathcal{ML}\left({\cal E}_{\bS_{\rm{c}}}(\{\omega_e\}_{e\in \mathcal{I}_{\bf b}});\{s_e\}_{e\in \mathcal{I}_{\bf r}}, z_1,...,z_{m-r} \right)=\mathcal{ML}\left(\mathcal{E}_{\bS};\{\omega_e\}_{e\in \mathcal{I}_{\bf b}},\{s_e\}_{e\in \mathcal{I}_{\bf r}},z_1,...,z_{m-r}\right).
\ee
\ec
Combining this with  Theorem \ref{thmint}, we can represent Mellin--Laplace transformed weighted exponential volumes by Barnes integrals. To recover the weighted exponential volumes ${\cal E}_{\bS_{\rm{c}}}(\omega; \K_{\bf r}, \L)$, we apply the inverse Mellin/Laplace  transform to the function $\mathcal{ML}\left(\mathcal{E}_{\bS};\{s_e\}_{e\in \mathcal{I}_{\bf b}},\{s_e\}_{e\in \mathcal{I}_{\bf r}},z_1,...,z_{m-r}\right)$.

\bc  For any weighted colored surface $(\bS_{\rm{c}}; \omega)$, the weighted exponential volume ${\cal E}_{\bS_{\rm{c}}}(\omega;\K_{\bf r}, \L)$ is finite. 
\ec

\subsection{Chekhov's volumes = Mellin transformed exponential volumes  at $s_e=\frac{1}{2}$}
In \cite{Ch24}, L. Chekhov  defined volumes of the moduli spaces ${\cal M}_\bS(\L)$ of cusped hyperbolic surfaces of the topological type described by a decorated surface $\bS$, with a given set $\L$ of 
 lengths of the neck geodesics, which we call {\it Chekhov's volumes}. 
 They are functions in $\L$   which depend  only on $\bS$. 

On the other hand, our exponential volumes of the moduli spaces ${\cal M}_\bS(\K, \L)$ of the ideal hyperbolic structures \cite{GS24} are functions ${\cal E}_\bS(\K,\L)$. We recall that the very notion of the ideal hyperbolic structure on $\bS$ includes a choice of a horocycle at each cusp. 

The next Proposition shows that for any decorated surface $\bS$,  Chekhov's volumes coincide with the Mellin transformed exponential volumes ${\cal E}_\bS(\K, \L)$ from \cite{GS24}, evaluated at $s_e= \frac{1}{2}$, where 
$e \in {\cal I}_\bS$. 
\bp
Let  $\bS_{\rm{c}}$ be  a decorated surface $\bS$   with all vertices colored blue. Set  $\omega_e=\frac{1}{2}$. Then 
\be
\mbox{\rm Chekhov's volume of ${\cal M}_\bS$}(\L) = 
\mathcal{M}\left(\mathcal{E}_{\bS}(\L);\left\{\frac{1}{2}\right\}_{e\in \mathcal{I}_\bS}\right).
\ee
\ep

\bpr If all vertices on $\bS$ are colored blue, then all edges must be blue by our convention. Moreover, if the weights $\omega_e$ at the edges $e$ 
are  equal to $\omega_e=\frac{1}{2}$,   the weighted exponential volume form 
 coincides with the volume form considered  in \cite{Ch24}. Therefore 
\be
\mbox{\rm Chekhov's volume of ${\cal M}_\bS$}(\L) = {\cal E}_{\bS_{\rm{c}}}\left(\left\{\frac{1}{2}\right\}_{e\in \mathcal{I}_\bS};\L\right).
\ee
On the other hand, by Theorem \ref{thm:fc} we have
\be
{\cal E}_{\bS_{\rm{c}}}\left(\left\{\frac{1}{2}\right\}_{e\in \mathcal{I}_\bS};\L\right) = 
\mathcal{M}\left(\mathcal{E}_{\bS}(\L);\left\{\frac{1}{2}\right\}_{e\in \mathcal{I}_\bS}\right).
\ee
Proposition is proved. \epr


\subsection{Examples}

\subsubsection{Colored triangle}
\begin{figure}[ht]
	\centering
	\includegraphics[scale=0.5]{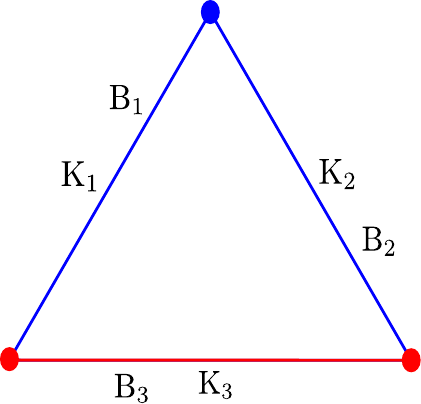}
	\small
	\caption{A colored triangle $\Delta_{\rm{c}}$.}
	\label{figure:tri}
\end{figure}
 
Let us compute the weighted exponential volume for the triangle $\Delta_{\rm{c}}$, colored as Figure \ref{figure:tri}, with the  weights $\omega=(\omega_1, \omega_2)$ and $t_1=\omega_1+\omega_2$. 
The potentials are 
\be
\begin{aligned}
&W_{\Delta_{\rm{c}}}=\B_2+\B_3=\sqrt{\frac{\K_2 \K_3}{\K_1}}+\sqrt{\frac{\K_3 \K_1}{\K_2}}.\\
&W_1 = \B_1 = \sqrt{\frac{\K_1 \K_2}{\K_3}}.\\
\end{aligned}
\ee
Recall \eqref{BK}.
Then the weighted exponential volume is
\be
\begin{aligned}
&{\cal E}_{\Delta_{\rm{c}}}(\omega_1, \omega_2, \K_3)
\\&=\Gamma(\omega_1+\omega_2) \int_{\mathcal{M}_{\Delta_{\rm{c}}}(\K_3)} \frac{\K^{\omega_1}_1\K^{\omega_2}_2}{W^{\omega_1+\omega_2}_1} e^{-\sqrt{\frac{\K_2 \K_3}{\K_1}}-\sqrt{\frac{\K_3 \K_1}{\K_2}}} d\log (\K_1/\K_2) 
\\&\stackrel{(\ref{thm:fc})}{=}\mathcal{M}\left(\mathcal{E}_{\Delta}(\K_3);\omega_1,\omega_2\right)
\\&=\int_0^\infty \int_0^\infty \K^{\omega_1}_1\K^{\omega_2}_2 e^{-\B_1-\B_2-\B_3} 2 d\log \B_1  d\log \B_2
\\&\stackrel{(\ref{BK})}{=} \int_0^\infty \int_0^\infty \left(\frac{\B_1 \K_3}{\B_2}\right)^{\omega_1} (\B_1 \B_2)^{\omega_2} e^{-\B_1-\B_2-\K_3/\B_2} 2 d\log \B_1  d\log \B_2
\\&=2\K_3^{\omega_1} \int_0^\infty \B_1^{\omega_1+\omega_2}e^{-\B_1}d\log\B_1 \int_0^\infty  \B_2^{\omega_2-\omega_1} e^{-\B_2-\K_3/\B_2}  d\log \B_2
\\&=2\Gamma(\omega_1+\omega_2) \K_3^{\frac{\omega_1+\omega_2}{2}} J_{\omega_1-\omega_2}(\K_3). 
\end{aligned}
\ee 
In particular,  for the weights $\omega_1=\omega_2=\frac{1}{2}$ we get
\be
\begin{aligned}
&{\cal E}_{\Delta_{\rm{c}}}(\frac{1}{2},\frac{1}{2} ; \K_3)
 \ = \ \mathcal{M}\left(\mathcal{E}_{\Delta}(\K_3);\frac{1}{2},\frac{1}{2}\right)
\ = \ 2\K_3 J_0(\K_3).
\end{aligned}
\ee

\subsubsection{Colored crown with one cusp}
For the case $\bS=\D_1^*$, we parametrise the moduli space as in Section \ref{subsection:d1}, and the color of the vertex and the color of the edge are both blue, and $t=\omega$. Then the potential
\[W=\B+\B\L=\sqrt{\K}(e^{\frac{l}{2}}+e^{-\frac{l}{2}}).\]
Then by the second line of Equation \eqref{eq:d1comp}, we get
\be
\begin{aligned}
&{\cal E}_{\D^*_{1,\rm{c}}}(\omega;\L)
\\&\stackrel{(\ref{thm:fc})}{=} \mathcal{M}\left(\mathcal{E}_{\D_1^*}(\L);\omega\right)
\\&= \int_0^{+\infty}e^{-(\sqrt{\K}(e^{\frac{l}{2}}+e^{-\frac{l}{2}}))} \K^{\omega}   d\log \K 
\\&= 2\Gamma(2\omega)  \left(e^{\frac{l}{2}}+e^{-\frac{l}{2}}\right)^{-2\omega}.
\end{aligned}
\ee

\subsubsection{Colored quadrilaterals}
\begin{figure}[ht]
	\centering
	\includegraphics[scale=0.35]{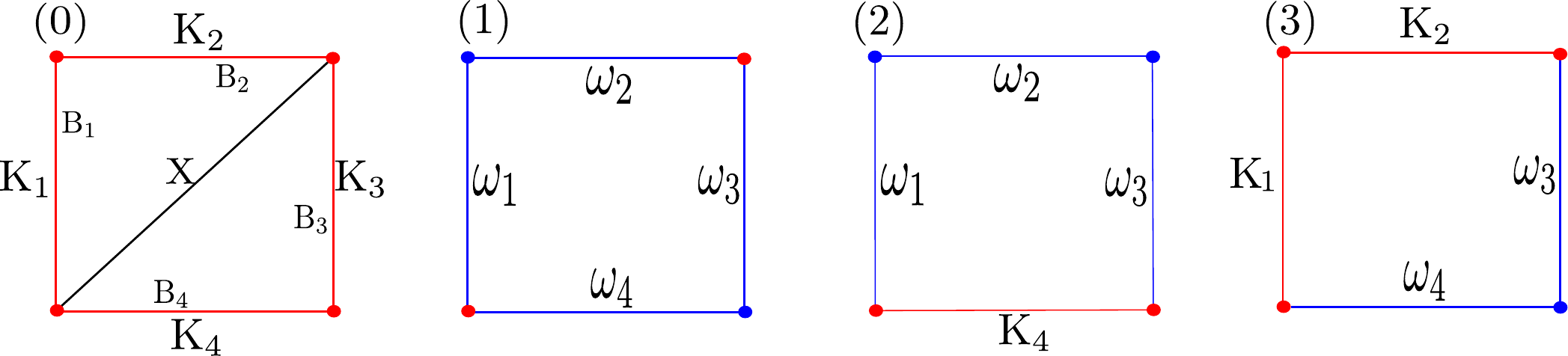}
	\small
	\caption{Different colored quadrilaterals.}
	\label{figure:square}
\end{figure}

0. For the case $\bS=\square$ in Figure \ref{figure:square}(0), the cluster Poisson coordinates are $\X,\B_1,\B_2,\B_3,\B_4$, while the boundary $\K$-parameters are $\K_1,\K_2,\K_3,\K_4$. We have
\[\B_1 \B_2=\K_2,\;\; \B_2 \X \B_3=\K_3,\;\; \B_3 \B_4= \K_4,\;\; \B_4 X \B_1=\K_1.\]
The potentials
\[W_1=\B_1,\;\;\;W_2=\B_2+\B_2 \X,\;\;\;W_3=\B_3,\;\;\;W_4=\B_4+\B_4 \X.\]

1. For Figure \ref{figure:square}(1), we have the potential
\[W_{\square_{\rm{c}_1}}=W_2+W_4.\]
The weight $\omega=(\omega_1,\omega_2,\omega_3,\omega_4)$ with
$t_1=\omega_1+\omega_2,\;\;\;t_3=\omega_3+\omega_4$. 
Then we get
\be
\begin{aligned}
&{\cal E}_{\square_{\rm{c}_1}}(\omega_1, ..., \omega_4)
\\&=\Gamma(\omega_1+\omega_2)\Gamma(\omega_3+\omega_4)\int_{\mathcal{M}_{\square_{\rm{c}_1}}} \frac{\K_1^{\omega_1} \K_2^{\omega_2} \K_3^{\omega_3} \K_4^{\omega_4}}{W_1^{\omega_1+\omega_2} W_3^{\omega_3+\omega_4}} \exp(-W_2-W_4) d \rm{Vol}
\\&\stackrel{(\ref{thm:fc})}{=} \mathcal{M}\left(\mathcal{E}_{\square};\omega_1,\omega_2,\omega_3,\omega_4\right)
\\&\stackrel{(\ref{eq:square4})}{=}    2\frac{\prod_{1 \leq i<j\leq 4}\Gamma(\omega_i+\omega_j)}{\Gamma(\omega_1+\omega_2+\omega_3+\omega_4)}.
\end{aligned}
\ee 
 
Note that we  get the same weighted exponential volume if the vertices $1$ and $3$ are colored red, or if just one of the vertices is colored red, or if there is no red vertices. 
Indeed, the answer depends only on the weights at the blue intervals, which are the same in all these cases. 
 The dimension of the relevant moduli space  is equal to 1+ the number of red vertices.  So the weighted exponential volume forms are different, and even the dimensions, 
  can be different, while the integrals are the same. This demonstrates the general phenomenon explained in Remark \ref{REM4.4}. 
\vskip 1mm

2. For Figure \ref{figure:square}(2), we have the potential
\[W_{\square_{\rm{c}_2}}=W_3+W_4.\]
The weight $\omega=(\omega_1,\omega_2,\omega_3)$ with
$t_1=\omega_1+\omega_2, \;\;\; t_2=\omega_2+\omega_3$. 
Then we get
\be
\begin{aligned}
&{\cal E}_{\square_{\rm{c}_2}}(\omega_1, \omega_2, \omega_3, \K_4)
\\&=\Gamma(\omega_1+\omega_2)\Gamma(\omega_2+\omega_3)\int_{\mathcal{M}_{\square_{\rm{c}_2}}(\K_4)} \frac{\K_1^{\omega_1} \K_2^{\omega_2} \K_3^{\omega_3}}{W_1^{\omega_1+\omega_2} W_2^{\omega_2+\omega_3}} \exp(-W_3-W_4) d \rm{Vol}
\\&\stackrel{(\ref{thm:fc})}{=}\mathcal{M}_{\rm blue}\left(\mathcal{E}_{\square}(\K_4);\omega_1,\omega_2,\omega_3\right).
\end{aligned}
\ee 
Let $\mathcal{M}^{-1}(f(\omega_4);\K_4)$ be the inverse Mellin transform of $f$ transforming the parameter $\omega_4$ into $\K_4$. Then by Equation \eqref{eq:square4}, the last term of the last equation is
\be
\begin{aligned}
\\&=\mathcal{M}^{-1}\left(\frac{2\prod_{1\leq i<j\leq 4}\Gamma(\omega_i+\omega_j)}{\Gamma(\omega_1+\omega_2+\omega_3+\omega_4)}; \K_4\right)
\\&=2\prod_{1 \leq i < j \leq 3}\Gamma(\omega_i+\omega_j)\ \mathcal{M}^{-1}\left(\frac{\Gamma(\omega_1+\omega_4)\Gamma(\omega_2+\omega_4)\Gamma(\omega_3+\omega_4)}{\Gamma(\omega_1+\omega_2+\omega_3+\omega_4)}; \K_4\right)
\\&=\frac{\prod_{1 \leq i < j \leq 3}\Gamma(\omega_i+\omega_j)\ }{\pi i} \int_{-i\infty}^{+i\infty}\frac{\Gamma(\omega_1+\omega_4)\Gamma(\omega_2+\omega_4)\Gamma(\omega_3+\omega_4)}{\Gamma(\omega_1+\omega_2+\omega_3+\omega_4)}\K_4^{-\omega_4} d\omega_4
\\&=2\prod_{1 \leq i < j \leq 3}\Gamma(\omega_i+\omega_j)\  \MeijerG{3}{0}{1}{3}{\omega_1+\omega_2+\omega_3}{\omega_1,\omega_2,\omega_3}{\K_4}.
\end{aligned}
\ee

3. For Figure \ref{figure:square}(3), we have the potential
\[W_{\square_{\rm{c}_3}}=W_1+W_2+W_4.\]
The weights are $\omega=(\omega_3,\omega_4)$ with
$t_3=\omega_3+\omega_4$. 
Then we get
\be
\begin{aligned}
&{\cal E}_{\square_{\rm{c}_3;\omega}}( \K_1,\K_2, \omega_3, \omega_4)
\\&=\Gamma(\omega_3+\omega_4)\int_{\mathcal{M}_{\square_{\rm{c}_3}}(\K_1,\K_2) } \frac{\K_3^{\omega_3} \K_4^{\omega_4}}{W_3^{\omega_3+\omega_4}} \exp(-W_1-W_2-W_4) d \rm{Vol}
\\&\stackrel{(\ref{thm:fc})}{=}\mathcal{M}_{(3,4)}\left(\mathcal{E}_{\square}(\K_1,\K_2);\omega_3, \omega_4\right)
\\&= \mathcal{M}_{(1,2)}^{-1}\left(\frac{2\prod_{1\leq i<j\leq 4}\Gamma(\omega_i+\omega_j)}{\Gamma(\omega_1+\omega_2+\omega_3+\omega_4)}; \K_1,\K_2\right)
\\&=\frac{2}{(2\pi i)^2}\int_{-i\infty}^{+i\infty} \int_{-i\infty}^{+i\infty}\frac{\prod_{1\leq i<j\leq 4}\Gamma(\omega_i+\omega_j) }{\Gamma(\omega_1+\omega_2+\omega_3+\omega_4)}\K_1^{-\omega_1}\K_2^{-\omega_2} d\omega_1 d\omega_2.
\end{aligned}
\ee 
In last line  stands a Meijer G-function with two parameters. Here ${\cal M}_{(3,4)}$ is the Mellin transform in the variables $(\K_3, \K_4) \to (\omega_3, \omega_4)$, and 
${\cal M}^{-1}_{(1,2)}$ is the inverse Mellin transform in the variables $ (\omega_1, \omega_2)\to (\K_1, \K_2) $. \\

The  last  two examples show the main steps of the calculation of the  weighted exponential volumes of moduli spaces related to colored decorated surfaces using Barnes integrals. Let us summarize them, as illustrated on Figure \ref{figure:square4}.

\begin{figure}[ht]
	\centering
	\includegraphics[scale=0.8]{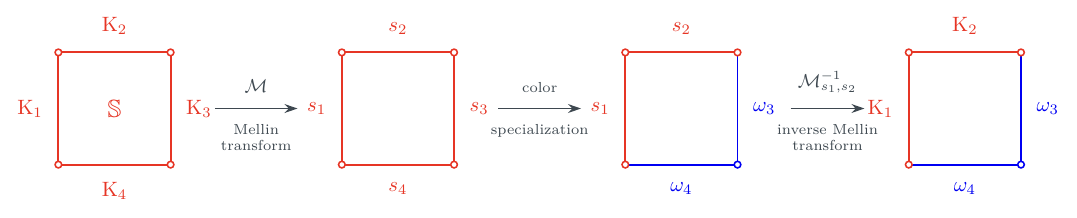}
	\small
	\caption{Schematic calculation of the weighted exponential volumes using Barnes integrals.}
	\label{figure:square4}
\end{figure}

\begin{enumerate}

\item Start with  the exponential volume  of a decorated surface $\bS$, see the left picture of Figure \ref{figure:square4}: 
$$
{\cal E}_{\bS}(\K_1, ..., \K_n, \L).
$$

 \item Do the Mellin transform $\K_i\to s_i$,  getting a function  in the Mellin variables 
$(s_1, ..., s_n)$ and $\L$:
 $$
 {\cal M}({\cal E}_{\bS}; s_1, ..., s_n, \L).
 $$
 By one of the main results of this paper, it is given by a Barnes $\Gamma-$integral.

\item Color  in blue some marked points of $\bS$. This way we get   a colored surface $\bS_{\rm c}$ with $d$ red  and $b$ blue boundary intervals,  $d+b=n$. 
 Specialize the Mellin variables $(s_1, ..., s_b)$ at the blue boundary intervals to the weights $(\omega_1, ..., \omega_b)$ -  the weights on  the colored surface $\bS_{\rm c}$. 
This way we get a function 
$$
{\cal M}({\cal E}_{\bS}; \omega_{1}, ..., \omega_{b}; s_{b+1}, ..., s_{b+d}, \L)
$$

\item Apply the inverse Mellin transform to  the red Mellin variables $(s_{b+1}, ..., s_{b+d})$. It transforms them back to  the $\K-$variables $(\K_{b+1}, ..., \K_{b+d})$. 
 The resulting function is 
$$
{\cal M}^{-1}_{\rm red}{\cal M}({\cal E}_{\bS}; \omega_{1}, ..., \omega_{b}; \K_{b+1}, ..., \K_{b+d}, \L).
$$
 
\end{enumerate}

The inverse Mellin transform applied to a Barnes $\Gamma-$integral can be viewed as a generalised  Meijer-$\rm G$ function, with $d-$fold integration. 
\vskip 1mm

One can also
do this in just two steps:
\begin{figure}[ht]
	\centering
	\includegraphics[scale=0.8]{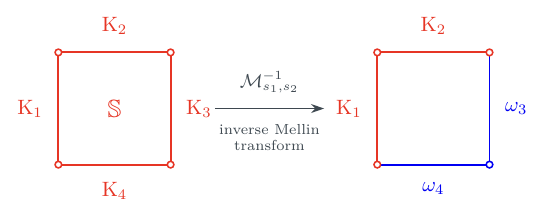}
	\small
	\caption{Schematic calculation of the weighted exponential volumes}
	\label{figure:square4b}
\end{figure}

\begin{enumerate}

\item Start with  the exponential volume  $
{\cal E}_{\bS}(\K_1, ..., \K_n, \L)
$ for a decorated surface $\bS$. 
 Color  in blue some marked points of $\bS$,  getting   a colored surface $\bS_{\rm c}$ with $d$ red  and $b$ blue boundary intervals,  $d+b=n$. 

\item Do the Mellin transform $\K_i\to \omega_i$ for the $\K-$variables at the blue boundary intervals: 
 $$
{\cal M}_{\rm blue}({\cal E}_{\bS}; \omega_{1}, ..., \omega_{b}; \K_{b+1}, ..., \K_{b+d}, \L).
 $$
 
 Evidently, we have 
 \be
 \begin{split}
&{\cal M}^{-1}_{\rm red}{\cal M}({\cal E}_{\bS}; \omega_{1}, ..., \omega_{b}; \K_{b+1}, ..., \K_{b+d}, \L) \\
 &= {\cal M}_{\rm blue}({\cal E}_{\bS}; \omega_{1}, ..., \omega_{b}; \K_{b+1}, ..., \K_{b+d}, \L). \\
  &\stackrel{}{=} {\cal E}_{\bS}(\omega_{1}, ..., \omega_{b}; \K_{b+1}, ..., \K_{b+d}, \L). \\
  \end{split}
 \ee
 
\end{enumerate}

\subsection{Teichm\"uller spaces for colored surfaces  via cluster varieties} \la{Sect5.3.3}

Let $G=\operatorname{PGL}_2$. Denote by $\B$  a maximal Borel subgroup,  by $\U$ the maximal unipotent subgroup   of $G$, by  $\mathcal{B}:=\mathbb{P}^1$ 
the flag variety, and by $\mathcal{A}:=G/U = (\A^2-0)/\pm 1$ the decorated  flag variety. 

 Let $\mathcal{L}$ be a $\operatorname{PGL}_2$-local system on $\bS$. 
A {\em framing/decoration} on $\mathcal{L}$ at a point $p$ on $\bS$ is a flat section of the associated local system $\mathcal{L}_{\mathbb{P}^1}:= {\cal L}/\B$/$\mathcal{L}_{\mathcal{A}}:={\cal L}/\U$ near $p$.

\bd Let $\bS_{\rm{c}}=(\bS; {\cal V}_{r})$ be a colored surface. The space $\mathcal{P}_{\bS_{\rm{c}}}$ parametrises triples $(\mathcal{L},\alpha,\beta)$ where $\mathcal{L}$ is a $\operatorname{PGL}_2$-local system on $\bS$ with a decoration $\alpha$ at  red marked points, and a framing $\beta$ at the blue marked points and  puncture. 
\ed
 
The space $\mathcal{P}_{\bS_{\rm{c}}}$ carries canonical ${\rm Mod}(\bS)-$equivariant  cluster Poisson structure \cite{GSh19}. Its cluster are parametrised by the triangulations of the associated decorated surface $\bS'$, 
 defined as follows. For each red point on $\bS_{\rm c}$ we introduce an extra marked boundary point on $\bS'$. So if $\bS$ has $r$ red points and $b$ blue points, then $\bS'$ has $2r+b$ marked boundary points, see \cite[Proposition 2.13]{GS24}. 
The $\mathcal{X}$-coordinates are defined in \cite[Section 2.1.1]{GS24}.  
So the same way as in \cite[Definition 2.8]{GS24}, the $\mathcal{P}_{\bS_{\rm{c}}}$ and its subspace $\mathcal{P}_{\bS_{\rm{c}}}(\K, \L)$ carry the ${\rm Mod}(\bS)$-invariant cluster volume forms 
\[\Omega_{\bS_{\rm{c}}},\;\;\;\Omega_{\bS_{\rm{c}}}(\K,\L).\]

The  space $\mathcal{P}_{\bS_{\rm{c}}}(\mathbb{R}_{>0})$ coincides with the enhanced Teichm\"uller space ${\cal T}_{\bS'}$, parametrising the ideal hyperbolic structures on the decorated surface $\bS'$ plus  a choice of an orientation of the boundary for each  hole.  The Teichm\"uller space ${\cal T}_{{\bS}_c}$ is a subspace: 
\be
{\cal T}_{{\bS}_c}\subset \mathcal{P}_{\bS_{\rm{c}}}(\mathbb{R}_{>0}) = {\cal X}_{\bS'}(\R_{>0}).
\ee 
The cutting and gluing maps in \cite[Section 2.2]{GS24} for $\mathcal{P}_{\bS}$ which works the same way for $\mathcal{P}_{\bS_{\rm{c}}}$, except that we could only cutting and gluing along an edge with two red vertices.

\subsubsection{Cluster Poisson structure on the space of configurations $(\A_p, \F_1,...,\F_\N, \A_q)$.}
\label{section:cpAF}

Let us elaborate an example of the red-based half-disc $\D^-_{\N+2,\rm{c}}$, that is a colored surface with two adjacent red points $p,q$ and the other $\N$ points blue, see Figure \ref{figure:interval1}.    We will use this in Sections \ref{Sect9},  \ref{Sect10}.

Let us introduce coordinates on the space of configurations $(\A_p, \F_1,...,\F_\N, \A_q)$. Here $\A_p, \A_q$ 
are vectors in a real two dimensional symplectic vector space $(\V_2, \omega_2)$ with the symplectic form $\omega_2$, and $\F_i$ are lines, that is one-dimensional subspaces  in $\V_2$. The positivity means that  
$\omega_2(\A_p, \A_q)>0$,  the lines $\F_1, ..., \F_\N$ are inside of the convex cone spanned by the vectors $\A_p, \A_q$, and follow clockwise one after the other. 
To define the local potential  at $p$ for the triangle $(\F_{i+1}, \A_p, \F_{i})$, we choose a non-zero vector $v_i\in \F_i$ for $i=1,...,\N$. Let $v_{\N+1}:=\A_q$. Set
\be \la{WL}
  W(\F_{i+1}, \A_p, \F_{i}) := \frac{\omega_2(v_i, v_{i+1})}{\omega_2(\A_p, v_i)\omega_2(\A_p, v_{i+1})}, \ \ \ \ \lambda:= \omega_2(\A_p, \A_q).
\ee
 The $W$ does not depend on the choice of the vector $v_i$.     
\bl 
\la{LPO}
The positivity of the configuration $(\A_p, \F_1,...,\F_\N, \A_q)$ is equivalent to $$
W(\F_{i+1}, \A_p, \F_{i})>0 \ \ \forall i=1, ..., \N \ \mbox{and} \ \lambda>0.
$$
\el

Let us consider the ideal triangulation $\tau_p$ such that all the edges start from $p$. We define 
\be \la{BI}
\B_i:=W(\F_{i+1}, \A_p, \F_{i}), \ \ \ \ i=1,...,\N.
\ee
 Then the cluster Poisson coordinates for the triangulation $\tau_p$ are 
 $$
 \B_\N, \ \ \K_{pq}/\B_\N, \ \ \X_i=\B_i/\B_{i+1}\ \ \forall i=1,...,\N-1.
 $$
  Thus 
\be
\Omega_{\D_{\N+2,\rm{c}}}=2 \bigwedge_{i=1}^\N d\log \B_i \wedge d\log \K_{pq}.
\ee

 \subsubsection{The weighted exponential volume on ${\rm Conf}_{+}(\A_p, \F_1,...,\F_\N, \A_q)$} 
\begin{figure}[ht]
	\centering
	\includegraphics[scale=0.7]{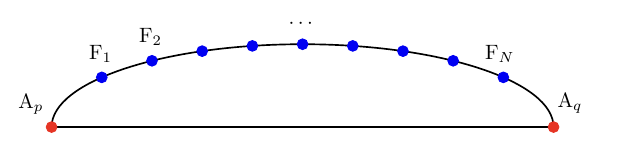}
	\small
	\caption{A red-based half-disc $\D^-_{\N+2, {\rm c}}$ is  a disc with two adjacent red points.}
	\label{figure:interval1}
\end{figure}

 The Teichm\"uller space ${\cal T}_{\D^-_{\N+2,\rm{c}}}$ coincides with the space of positive configurations $(\A_p, \F_1,...,\F_\N, \A_q)$, as well as with the moduli space  ${\cal M}_{\D^-_{\N+2,\rm{c}}}$: 
 $$
 {\cal M}_{\D^-_{\N+2,\rm{c}}} = {\cal T}_{\D^-_{\N+2,\rm{c}}} = {\rm Conf}_{+}(\A_p, \F_1,...,\F_\N, \A_q).
 $$

 Given an ideal hyperbolic structure and a weight $\omega$ on $\D^-_{\N+2,\rm{c}}$, recall  from Definition \ref{LD.11} the weight factor ${\cal Q}_{\D^-_{\N+2,\rm{c}}}(\omega)$, 
  the potential of $\D^-_{\N+2,\rm{c}}$, given by the sum of local potentials at the two  red points:
\be
W_{\D^-_{\N+2,\rm{c}}}=W_p+W_q, 
\ee
and the weighted exponential volume form.  
\bd \label{defn:weightDh}
For any $\kappa\in \C$ and $\beta\in \R$, the { $(\kappa,\beta)$-weighted exponential volume form} for the red-based half-disc $\D^-_{\N+2,\rm{c}}$ as on Figure \ref{figure:interval1} is 
\be \la{EVFK} 
{\Bbb E}_{\D^-_{\N+2,\rm{c}}}^{(\kappa,\beta)}(\omega):={\cal Q}_{\D^-_{\N+2,\rm{c}}}^\kappa(\omega) e^{-\beta W_{D_{\N+2,\rm{c}}}} \Omega_{\D^-_{\N+2,\rm{c}}}
\ee
\ed
The parameter $\kappa$  changes the weight from $\omega$ to $\kappa\cdot\omega$, and the parameter $\beta$ rescales the potential. Complex powers make sense since thanks to the positivity of 
configurations $(\A_p, \F_1,...,\F_\N, \A_q)$ and Lemma \ref{LPO} the factor ${\cal Q}_{\D^-_{\N+2,\rm{c}}}(\omega)$ is   positive, and so 
$$
{\Bbb E}_{\D^-_{\N+2,\rm{c}}}^{(\kappa,\beta)}(\omega)={\rm exp}\Bigl(\kappa \log  {\cal Q}_{\D^-_{\N+2,\rm{c}}}(\omega) -\beta W_{D^-_{\N+2,\rm{c}}}\Bigr)\Omega_{\D^-_{\N+2,\rm{c}}}.$$

Let us glue two such half-discs $\D^-_{\N_1+2,\rm{c}_1}(\omega_1)$, $\D^-_{\N_2+2,\rm{c}_2}(\omega_2)$ and a triangle $\Delta=pqu$ along the two pairs of the 
red edges $pq$ and $qu$ as in Figure \ref{figure:interval3}.
\begin{figure}[ht]
	\centering
	\includegraphics[scale=0.7]{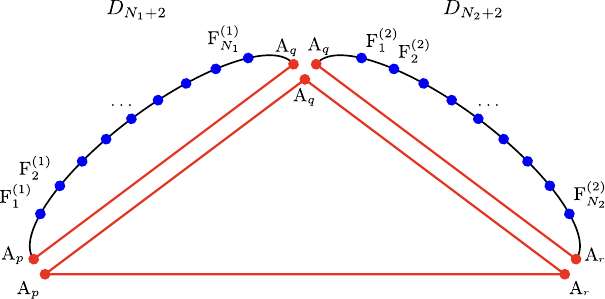}
	\small
	\caption{Glueing two red-based half-discs discs and a red triangle.}
	\label{figure:interval3}
\end{figure}
Denote by  $(\D^-_{\N_1+\N_2+3,\rm{c}};\omega)$ the obtained weighted colored surface where the weight $\omega$ and the color $\rm{c}$ are inherited from $\omega_1$, $\omega_2$ and $\rm{c}_1$, $\rm{c}_2$. 
Let 
\[{\Bbb E}^{(\kappa,\beta)}_{\D^-_{\N_1+2,\rm{c}_1}}(\omega_1), \ \ \;{\Bbb E}^{(\kappa,\beta)}_{\D^-_{\N_2+2,\rm{c}_2}}(\omega_2),\ \ \; {\Bbb E}^{(\kappa,\beta)}_{\D^-_{\N_1+\N_2+3,\rm{c}}}(\omega)\]
be their $(\kappa,\beta)$-weighted exponential volume forms. Since the triangle ${\Bbb E}_{\Delta}$ is colored red, its $(\kappa,\beta)$-weighted exponential volume form is 
\[{\Bbb E}^{(\kappa,\beta)}_{\Delta}= e^{- \beta \left(\sqrt{\frac{\K_{qr}\K_{pr}}{\K_{pq}}}+ \sqrt{\frac{\K_{pq}\K_{pr}}{\K_{qr}}} +\sqrt{\frac{\K_{qr}\K_{pq}}{\K_{pr}}}\right)}d\log \K_{pq}\wedge d\log \K_{qr} \wedge d\log \K_{pr}.\] 
The cutting and gluing property for the exponential volume forms \cite{GS24} immediately implies a similar property for the  weighted exponential volume form ${\Bbb E}_{\bS_{\rm{c}}}(\omega)$. So in our case we arrive to 
\bp
\label{prop:cge}
\be
  {\Bbb E}^{(\kappa,\beta)}_{\D^-_{\N_1+2,\rm{c}_1}}(\omega_1;\K_{pq})\wedge{\Bbb E}^{(\kappa,\beta)}_{\D^-_{\N_2+2,\rm{c}_2}}(\omega_2;\K_{qr})\wedge {\Bbb E}^{(\kappa,\beta)}_{\Delta}={\Bbb E}^{(\kappa,\beta)}_{\D^-_{\N_1+\N_2+3,\rm{c}}}(\omega).
\ee
\ep

\section{Recursions} \la{Section6}
By \cite[Theorem 5.3]{GS24}, we obtain the recursion formula for exponential volumes by cutting off trouser legs and triangles. 
The formula is a sum of  finitely many integrals. We briefly recall the story in Section \ref{Sec5.1}. 

In Section \ref{Sec5.2} we show that the Mellin--Laplace transforms of each integral in the recursion formula has a Barnes integral presentation. This way we get recursions for the Barnes integrals for Mellin-Laplace transformed exponential volumes. They provide  linear relations among the Barnes integrals.

In Section \ref{SEC10}  we prove a recursion  for the weighted exponential volumes of the  moduli spaces  ${\cal M}_{\bS_{\rm{c}}}$. It suggested the shape of  recursions  
for the  volumes of probabilistic moduli spaces  ${\cal M}_{\bS_\varphi}$ in Section \ref{Sect9}.

\subsection{Recursion for the exponential volumes  \cite{GS24}} \la{Sec5.1}
\label{subsection:re}
Let us recall the recursion formula for exponential volumes  in \cite[Theorem 5.3]{GS24}.

Consider an ideal hyperbolic structure on a decorated surface $\bS$, choose a cusp $p$ on a crown of $\bS$. 
As in Figure \ref{figure:cuttriangle}(1), let $\T$ be an embedded trouser leg containing the cusp $p$ with the boundary circle $\ell_\T$ length $l_\T$ and boundary $\K$-coordinate $\K_\T$, the potential is 
\[W_\T=\sqrt{\K_\T}\left(e^{\frac{l_T}{2}}+e^{-\frac{l_T}{2}} \right).\]
If $\ell_\T\not\subset \partial \bS$, the McShane gap term is 
\[2\sqrt{\K_\T}e^{-\frac{l_T}{2}}.\]
If $\ell_\T\subset \partial \bS$, the McShane gap term is 
\[W_\T=\sqrt{\K_\T}\left(e^{\frac{l_T}{2}}+e^{-\frac{l_T}{2}} \right).\]
Let ${\cal H}_p^{\circ}$ (${\cal H}_p^{\partial}$ resp.) be the collection of homotopy classes of embedded trouser legs containing the cusp $p$ such that the boundary circle is not (is resp.) the boundary of $\bS$, and ${[\cal H}_p^{\circ}]$ ($[{\cal H}_p^{\partial}]$ resp.) its finitely many mapping class group orbit.

As in Figure \ref{figure:cuttriangle}(2), let $\tau$ be an ideal triangle containing the cusp $p$. We say that $\tau$ is {\em $p$-narrowest} if there is not any other ideal triangle contained in $\tau$ around the region of cusp $p$. Suppose the $\K$-coordinates are $\K_1$, $\K_2$, $\K_3$ with $K_1$ along $\partial \bS$ opposite to $p$, and the $\B$-coordinates $\B_1$, $\B_2$, $\B_3$ as in Figure \ref{figure:cuttriangle}.
The potential is 
\be
W_\tau(\K_1, \K_2, \K_3)=\B_1+\B_2+\B_3= \sqrt{\frac{\K_1 \K_2}{\K_3}}+ \sqrt{\frac{\K_2 \K_3}{\K_1}}+ \sqrt{\frac{\K_3 \K_1}{\K_2}}.
\ee
The McShane gap term is 
\[\sqrt{\frac{\K_2 \K_3}{\K_1}}.\]
Let $\Delta_p^\circ$ ($\Delta_p^\partial$ resp.) be the collection of homotopy classes of embedded $p$-narrowest ideal triangles containing the cusp $p$ such that exactly one (two resp.) ideal edges of the ideal triangle is on the boundary of $\bS$, and $[\Delta_p^\circ]$ ($[\Delta_p^\partial]$ resp.) its finitely many mapping class group orbit.
 \begin{figure}[ht]
	\centering
	\includegraphics[scale=0.5]{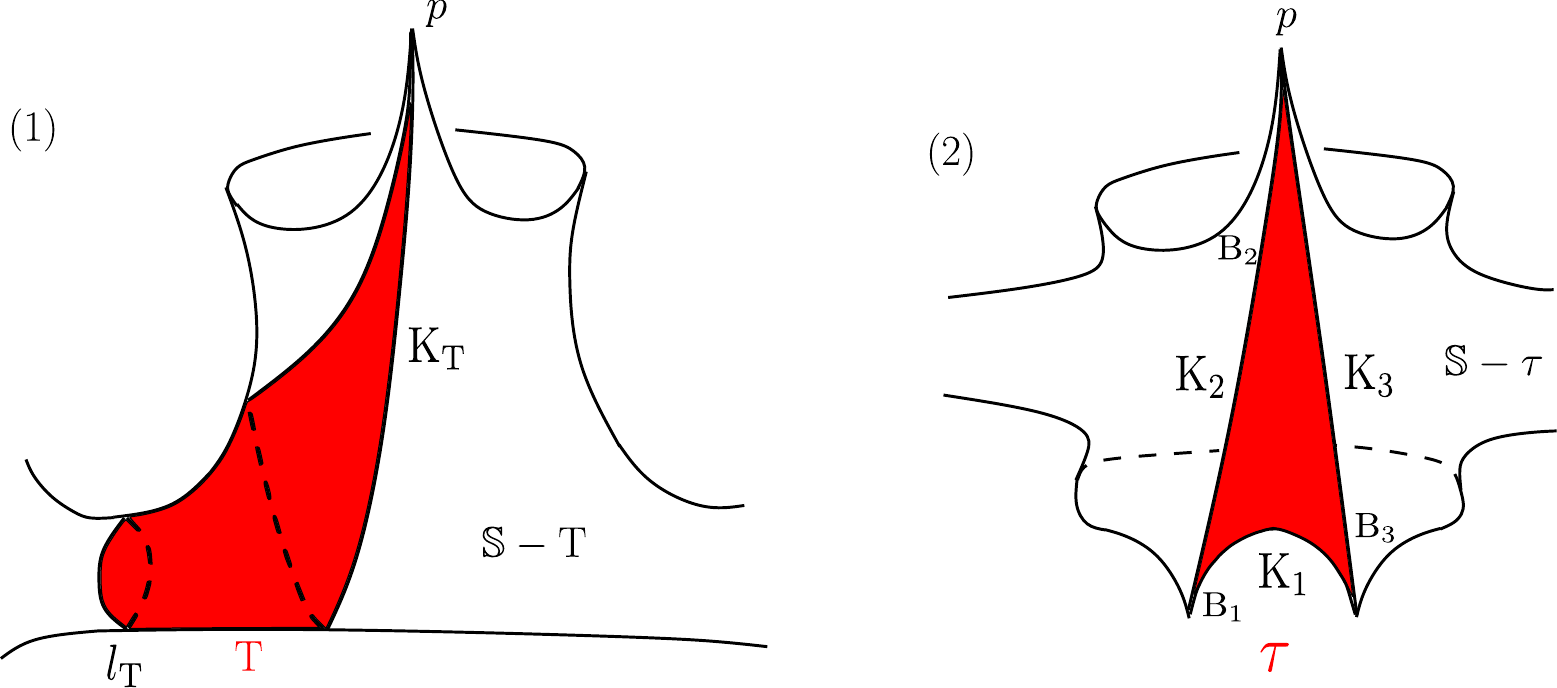}
	\small
	\caption{Figure (1) for cutting the trouser leg. Figure (2) for cutting the $p$-narrowest ideal triangle.}
	\label{figure:cuttriangle}
\end{figure}

\begin{theorem}
\label{thm:re}
{\rm \cite[Theorem 5.3]{GS24}} Let $\mu_\T$ be the number of $S_{1,1}$ obtained by cutting off $\T$ from $\bS$.
Given any decorated surface $\bS$ and a cusp $p$ on a crown of $\bS$, we have the following recursion formula with finitely many terms:
\begin{equation}
\label{eq:re}
\begin{aligned}
& \int _{{\cal M}_\bS({\K}, \L)}   W_p\  {\Bbb E}_\bS 
\\&= \sum_{[\T]\in  [{\cal H}_p^\circ]}   
 \int_0^{\infty}  \int_0^{\infty}   2^{1-\mu_\T}  \sqrt{\K_\T}e^{-\frac{l_\T}{2}} exp\left(-\sqrt{\K_\T}\left(e^{\frac{l_\T}{2}}+e^{-\frac{l_\T}{2}} \right)\right) \mathcal{E}_{\bS-\T}(\K,\L, \K_\T,l_\T) d\log \K_\T \cdot l_\T d l_\T+
\\& \sum_{[\T]\in  [{\cal H}_p^\partial]}   
 \int_0^{\infty}      \sqrt{\K_\T}\left(e^{\frac{l_\T}{2}}+e^{-\frac{l_\T}{2}} \right) exp\left(-\sqrt{\K_\T}\left(e^{\frac{l_\T}{2}}+e^{-\frac{l_\T}{2}} \right)\right) \mathcal{E}_{\bS-\T}(\K,\L,\K_\T) d\log \K_\T+
\\&\sum_{[\tau]\in  [\Delta_p^\circ]}   
 \int_0^{\infty}  \int_0^{\infty}    \sqrt{\frac{\K_2 \K_3}{\K_1}} {\it exp} \left(-W_\tau(\K_1, \K_2, \K_3)\right) 
  \mathcal{E}_{\bS-\tau}(\K\backslash \{\K_1\},\L,\K_2,\K_3) d\log \K_2 \cdot d\log \K_3+\\&\sum_{[\tau] \in  [\Delta_p^\partial]}   
 \int_0^{\infty}      \sqrt{\frac{\K_2 \K_3}{\K_1}} exp\left(-W_\tau(\K_1, \K_2, \K_3)\right) \mathcal{E}_{\bS-\tau}(\K\backslash \{\K_1\},\L,\K_2) d\log \K_2.
  \end{aligned}
\end{equation}
\end{theorem}

\subsection{Recursion for Mellin--Laplace transformed exponential volumes} \la{Sec5.2}
Let us Mellin--Laplace transform each integral of Equation (\ref{eq:re}). 

Below we  use the following notation for the Mellin--Laplace transform $\mathcal{ML}(\F)$ of a function $\F$,  separating by the semicolon $;$ the function from its arguments:
 \be
\mathcal{ML}(\F; \rm{S},\rm{Z}):= \mathcal{ML}(\F)( {\rm{S},\rm{Z}}).
\ee

\bt
\label{thm:girre}

Given any decorated surface $\bS$ and a cusp $p$ on a crown of $\bS$, we have the following recursion formula  for the  Barnes integrals for the   Mellin--Laplace transformed  exponetial volumes:
\begin{equation} \la{RFBI}
\begin{aligned}
& \mathcal{ML}(\mathcal{EW}_\bS; \rm{S}, \rm{Z}) 
\\&=  \sum_{[\T]\in  [{\cal H}_p^\circ]}   
 \int_{-i\infty}^{+i\infty}  \int_{-i\infty}^{+i\infty}    \frac{2^{1-\mu_\T}}{(2\pi i)^2} \mathcal{ML}\left(\sqrt{\K_\T}e^{-\frac{l_\T}{2}} \mathcal{E}_{\T}; s,z  \right) \mathcal{ML}\left(l_\T \mathcal{E}_{S-\T};{\rm{S}, \rm{Z}},-s,-z  \right) ds d z+
\\&  \sum_{[\T]\in  [{\cal H}_p^\partial]}   
   \int_{-i\infty}^{+i\infty}   \frac{1}{2\pi i} \mathcal{ML}\left(\sqrt{\K_\T}\left(e^{-\frac{l_\T}{2}}+e^{\frac{l_\T}{2}}\right) \mathcal{E}_{\T}; s,z \right) \mathcal{ML}\left(\mathcal{E}_{S-\T};{\rm{S}},{\rm{Z}},-s  \right) ds +
\\&\sum_{[\tau]\in  [\Delta_p^\circ]}   
\int_{-i\infty}^{+i\infty}  \int_{-i\infty}^{+i\infty}   \frac{1}{(2\pi i)^2}  \mathcal{M}\left(\sqrt{\frac{\K_2 \K_3}{\K_1}} \mathcal{E}_{\tau}; s_1,s_2,s_3  \right) \mathcal{ML}\left(\mathcal{E}_{S-\tau}; {\rm{S}\backslash \{s_1\},\rm{Z},}-s_2,-s_3  \right) d s_2 d s_3+
\\&\sum_{[\tau] \in  [\Delta_p^\partial]}   
\int_{-i\infty}^{+i\infty}    \frac{1}{2\pi i}  \mathcal{M}\left(\sqrt{\frac{\K_2 \K_3}{\K_1}} \mathcal{E}_{\tau}; s_1,s_2,s_3  \right) \mathcal{ML}\left(\mathcal{E}_{S-\tau};\rm{S}\backslash \{s_1\},\T,-s_2  \right) d s_2 ,
  \end{aligned}
\end{equation}
Here the first term is given by \eqref{eq:EW} and the terms of the recursion match the ones in (\ref{thm:re}).
\et
\bpr
Let us show that the Mellin--Laplace transform of each integral on the right in (\ref{eq:re}) with respect to boundary parameters $(\K,\L)$ has a Barnes integral presentation, 
which  induces the recursion (\ref{RFBI}).

Let us show that for the integral in first line  of Equation \eqref{eq:re}, which is the left hand side of the equality, 
\be \label{eq:EW}
\mathcal{ML}(\mathcal{EW}_\bS; \rm{S},\rm{Z}):=\mathcal{ML}\left(\int _{{\cal M}_\bS({\K}, \L)}  \ W_p\  {\Bbb E}_\bS\right)\left(\rm{S},\rm{Z}\right)
\ee
 has a Barnes integral presentation by induction. 
 
 1. Let  $\bS=\Delta$, and let us fix the notation as in Figure \ref{figure:D3}. Let us $t$-deform the vertex shared by edges for $\K_2$ and $\K_3$.
 By Lemma \ref{tLMEVL}  the Mellin transform of the $t-$deformed exponential volume is given by
 \be \la{A22}
  \mathcal{M}(\mathcal{E}_{\Delta, t}; s_1,s_2,s_3) =  t^{s_2+s_3} \mathcal{M}(\mathcal{E}_{\Delta}; s_1,s_2,s_3) .
    \ee
    On the other hand, by the very definition of the $t-$deformed exponential volume, 
    \be
\label{eq:WD3z}
\begin{aligned}
& \mathcal{M}(\mathcal{EW}_{\Delta}; s_1,s_2,s_3)={ \frac{d}{dt}}\left|\right._{t=1} \mathcal{M}(\mathcal{E}_{\Delta, t}; s_1,s_2,s_3).\\
\end{aligned}
\ee We have, using the $t-$deformed exponential volumes
\be
\label{eq:WD3}
\begin{aligned}
& \mathcal{M}(\mathcal{EW}_{\Delta}; s_1,s_2,s_3)\\
&\stackrel{(\ref{eq:WD3z})}{=}{ \frac{d}{dt}}\left|\right._{t=1} \mathcal{M}(\mathcal{E}_{\Delta, t}; s_1,s_2,s_3)\\
&\stackrel{(\ref{A22})}{=}(s_2+s_3)\mathcal{M}(\mathcal{E}_{\Delta, t}; s_1,s_2,s_3)
\\&\stackrel{(\ref{eq:D3})}{=}2 \Gamma(s_1+s_2) \Gamma(s_2+s_3+1) \Gamma(s_3+s_1).
\end{aligned}
\ee

2. Let $\bS=\D_1^*$. We fix the notation as in Figure \ref{figure:D1}.  By Lemma \ref{tLMEVL}  the Mellin transform of the $t-$deformed exponential volume is given by 
 \be \la{A22a}
  \mathcal{M}(\mathcal{E}_{\D_1^*, t}; s,z) =  t^{2s} \mathcal{M}(\mathcal{E}_{\D_1^*}; s,z) .
    \ee
   On the other hand, by the very definition of the $t-$deformed exponential volume, 
    \be
\label{eq:WD3zz}
\begin{aligned}
& \mathcal{M}(\mathcal{EW}_{\D_1^*}; s, z)={ \frac{d}{dt}}\left|\right._{t=1} \mathcal{M}(\mathcal{E}_{\D_1^*, t}; s, z).\\
\end{aligned}
\ee 
    We have, using the $t-$deformed exponential volumes
\be
\label{eq:WD1}
\begin{aligned}
& \mathcal{ML}(\mathcal{EW}_{\D_1^*}; s,z)\\
&\stackrel{(\ref{eq:WD3zz})}{=}{ \frac{d}{dt}}\left|\right._{t=1}\mathcal{ML}(\mathcal{E}_{\D_{1, t}^*}; s,z)\\
&\stackrel{(\ref{A22a})}{=}2s \mathcal{ML}(\mathcal{EW}_{\D_{1}^*}; s,z)
\\&= \frac{1}{\pi i}  \int_{-i \infty}^{+i \infty} \frac{\Gamma(2s+1)\Gamma(t)\Gamma(2s-t)\Gamma(s+z-t)}{\Gamma(2s)\Gamma(s+z+1-t)} d t,
\end{aligned}
\ee
where the contour of integration separates the poles of $\Gamma(2s-t)$ and $\Gamma(s+z-t)$ from those of $\Gamma(t)$.

 \vskip 1mm

 By Equation \eqref{eq:WD3}, we can use Parseval's formula as in Equation \eqref{eq:par} to reduce the general decorated surface $\bS$ to $\bS=S_{g,m}\cup \cup_{i=1}^r \D_1^*$.
We have
\[\mathcal{ML}(\mathcal{EW}_{\D_1^*}(\L); s)=2\Gamma(2s+1) (e^{\frac{l}{2}}+e^{-\frac{l}{2}})^{-2s}.\]
Then similarly to  \eqref{eq:sgmd1r} we get
\be
\begin{aligned}
&\mathcal{ML}(\mathcal{EW}_{\bS}; s_1,\cdots,s_r,z_1,\cdots,z_{m-r})
\\&=\int \sum_{0\leq d_1+ \cdots + d_m \leq d} {\cal V}_{d_1, \cdots, d_m}\prod_{i=1}^r l_i^{2d_i} \prod_{i=r+1}^m \frac{(2d_i)!}{z_{i-r}^{2d_i+1}} \prod_{i=1}^r 2\Gamma(2s_i+1) (e^{\frac{l_i}{2}}+e^{-\frac{l_i}{2}})^{-2s_i} l_i d l_i
\\&=\sum_{0\leq d_1+ \cdots + d_m \leq d}  2^r {\cal V}_{d_1, \cdots, d_m}\prod_{i=r+1}^m \frac{(2d_i)!}{z_{i-r}^{2d_i+1}} \prod_{i=1}^r \Gamma(2s_i+1) \prod_{i=1}^r \int_0^{+\infty} (e^{\frac{l_i}{2}}+e^{-\frac{l_i}{2}})^{-2s_i} l_i^{2d_i+1} d l_i
\\&\overset{\ref{lemma:Mn}}{=}\sum_{0\leq d_1+ \cdots + d_m \leq d}   \frac{2^r {\cal V}_{d_1, \cdots, d_m} \prod_{i=1}^r (2d_i+1)!\prod_{i=r+1}^m (2d_i)!} 
{\prod_{i=r+1}^m z_{i-r}^{2d_i+1}}   \\&\times \prod_{i=1}^r
 \frac{1}{2\pi i}\int_{-i \infty}^{+i \infty} \frac{\Gamma(2s_i+1)\Gamma(t)\Gamma(2s_i-t)\Gamma(s_i-t)^{2d_i+2}}{\Gamma(2s_i)\Gamma(1+s_i-t)^{2 d_i+2}} dt. 
\end{aligned}
\ee
Here for each $i$, the integration contour is a vertical  line  separating the poles of $\Gamma(2s_i-t)$ and $\Gamma(s_i-t)$ from those of $\Gamma(t)$.

To prove that the integral in the third line of  \eqref{eq:re} has a Barnes integral presentation, by using Parseval's formula as in  \eqref{eq:par}, it is enough to prove that 
\[\mathcal{ML}\left(\sqrt{\K_\T}\left(e^{\frac{l_\T}{2}}+e^{-\frac{l_\T}{2}} \right) exp\left(-\sqrt{\K_\T}\left(e^{\frac{l_\T}{2}}+e^{-\frac{l_\T}{2}} \right)\right); s,z\right)\]
has a Barnes integral presentation, which is shown in \eqref{eq:WD1}.

To prove that the integral in the fourth or fifth line of  \eqref{eq:re} has a Barnes integral presentation, by using Parseval's formula as in  \eqref{eq:par}, it is enough to prove that 
\[\mathcal{M}\left(\sqrt{\frac{\K_2 \K_3}{\K_1}} exp\left(-\sqrt{\frac{\K_1 \K_2}{\K_3}}- \sqrt{\frac{\K_2 \K_3}{\K_1}}- \sqrt{\frac{\K_3 \K_1}{\K_2}}\right); s_1,s_2,s_3\right)\]
has a Barnes integral presentation, which is shown in Equation \eqref{eq:WD3}.

To prove that the integral in the second line of Equation \eqref{eq:re} has a Barnes integral presentation, by using Parseval's formula as in  \eqref{eq:par}, it is enough to prove that  
\[\mathcal{ML}\left(\sqrt{\K_\T}e^{-\frac{l_\T}{2}}  exp\left(-\sqrt{\K_\T}\left(e^{\frac{l_\T}{2}}+e^{-\frac{l_\T}{2}} \right)\right); s,z\right)\]
and $l_{\T}\mathcal{ML}(\mathcal{E}_{\bS-\T}(\K,\L, \K_\T,l_{\T}); {\rm{S}},{\rm{Z}},s,z)$ have Barnes integral presentations.
By Theorem \ref{thmint}, we have $\mathcal{ML}(\mathcal{E}_\bS;\rm{S},\rm{Z})$ has a Barnes integral presentation for any decorated surface $\bS$. Since the Laplace transform of $l^k$ is $k!/z^{k+1}$, for any product $\mathbb{L}=\prod_{i=1}^{m-r} l_{i+r}^{m_i}$ of boundary circle lengths, we have $\mathcal{ML}(\mathbb{L}\cdot\mathcal{E}_\bS(\K,\L); \rm{S},\rm{Z})$ has a Barnes integral presentation. Thus $l_{\T}\mathcal{ML}(\mathcal{E}_{\bS-\T}(\K,\L, \K_\T,l_{\T}); {\rm{S}},{\rm{Z}},s,z)$ has a Barnes integral presentation. Moreover
\be
\label{eq:WDQ}
\begin{aligned}
&\int_0^{+\infty}\int_0^{+\infty} \sqrt{\K}e^{-\frac{l}{2}} e^{-(\sqrt{\K}(e^{\frac{l}{2}}+e^{-\frac{l}{2}}))} \K^{s-1} e^{-lz}  d\K dl
\\&= 2\Gamma(2s+1) \int_0^{+\infty}(e^{\frac{l}{2}}+e^{-\frac{l}{2}})^{-2s-1} e^{-\frac{l}{2}} e^{-lz}   dl
\\&= 2\Gamma(2s+1) \int_0^{+\infty}\left( e^{\frac{l}{2}} (1+e^{-l})\right)^{-2s-1} e^{-\frac{l}{2}} e^{-lz}   dl
\\&= 2\Gamma(2s+1) \int_0^{+\infty}\sum_{k=0}^{+\infty} \left( e^{-l(s+z+1)} \frac{(2s+1)_k}{k!}(-e^{-l})^k \right)  dl
\\&= 2\Gamma(2s+1) \sum_{k=0}^{+\infty} \left(\frac{(2s+1)_k}{s+z+k+1}\cdot \frac{(-1)^k}{k!}\right)
\\&= \frac{2\Gamma(2s+1)}{s+z+1} {}_2F_1(\{2s+1,s+z+1\},\{s+z+2\},-1)
\\&= \frac{2\Gamma(2s+1)}{s+z+1} \frac{\Gamma(s+z+2)}{\Gamma(2s+1)\Gamma(s+z+1)}\frac{1}{2\pi i} \int_{-i \infty}^{+i \infty} \frac{\Gamma(x)\Gamma(2s+1-x)\Gamma(s+z+1-x)}{\Gamma(s+z+2-x)} d x
\\&= \frac{1}{\pi i} \int_{-i \infty}^{+i \infty} \frac{\Gamma(x)\Gamma(2s+1-x)\Gamma(s+z+1-x)}{\Gamma(s+z+2-x)} d x.
\end{aligned}
\ee
Here the  integration contour separates the poles of $\Gamma(2s+1-x)$ and $\Gamma(s+z+1-x)$ from those of $\Gamma(x)$.
Hence each integral in the second line of Equation \eqref{eq:re} has a Barnes integral representation.
We conclude that  the Mellin--Laplace transform of each integral in Equation \eqref{eq:re} has a Barnes integral representation.
\epr
\begin{remark}
Theorem \ref{thm:girre} provides 
 a family of linear relations among Barnes integrals.  These relations can be also considered as the relations among the generalized hypergeometric functions ${}_p F_q$. 
\end{remark}

\subsection{Recursion for the weighted exponential volumes} \la{SEC10}

Let $(\bS_{\rm{c}};\omega)$ be a weighted  colored surface. In Section \ref{SEC10} we derive a recursion formula for the weighted exponential volume of the  moduli space  ${\cal M}_{\bS_{\rm{c}}}$. 
To run the recursion, let us pick  a red cusp $p$. 

Let us cut off from the surface $\bS$ a triangle $\tau$ or a trouser leg  $\T$ with a cusp at $p$, as in Figure \ref{figure:cuttriangle}.  
The colors and the weights for the surfaces $\bS-\tau$ and $\tau$, as well as the surfaces   $\bS-\T$ and $\T$, are inherited from the ones of  $(\bS_{\rm{c}}; \omega)$. The recursion formula will have several terms. 
Each of them matches one of the four terms on the right hand side of formula \eqref{eq:re} for the recursion for $\int_{{\cal M}_{\bS_{\rm{c}}}(\K_{\bf r}, \L)} W_p {\Bbb E}_{\bS_{\rm{c}}}(\omega;\K_{\bf r}, \L)$. 
Let us elaborate these terms.  

($\T$) Consider  a trouser leg $\T$ with the cusp $p$. Then the   $\K$-parameter of $\T$ is $\K_\T$, and the geodesic 
length of the circle is $\ell_\T$.  Suppose that the cusp $p$ lies on the crown end $\D_k^*$ as in Figure \ref{figure:cuttrouser}. After cutting off $\T$ from $(\bS_{\rm{c}}; \omega)$, the cusp $p$ splits into $p$ and $p'$,  
 and the crown $\D_k^*$ changes into $\D_{k+1}^*$. Since the cusp $p$ is red, the points $p$, $p'$ and the edge for $\K_\T$ must be red. The color $\rm{c}'$ and the weight $\omega'$ of the surface $\bS-\T$ is 
 inherited from  $(\bS_{\rm{c}}; \omega)$. So the weighted exponential volume $\mathcal{E}_{(\bS-\T)_{c'}}(\omega';\K_{\bf r},\L, \K_\T,l_\T)$ is defined. Recall the notation $[{\cal H}_p^\circ]$, $[{\cal H}_p^\partial]$, $[\Delta_p^\circ]$, $[\Delta_p^\partial]$ from Section \ref{subsection:re}. Then the first two sums of the recursion are
\be\la{S1}
\sum_{[\T]\in  [{\cal H}_p^\circ]}   
 \int_0^{\infty}  \int_0^{\infty}   2^{1-\mu_\T} \sqrt{\K_\T}e^{-\frac{l_\T}{2}} \exp\left(-2\sqrt{\K_\T} \ {\rm ch}\left({\frac{l_\T}{2}}\right) \right) \mathcal{E}_{(\bS-\T)_{c'}}(\omega';\K_{\bf r},\L, \K_\T,l_\T) d\log \K_\T \cdot l_\T d l_\T
 \ee
and
\be \la{S2} \sum_{[\T]\in  [{\cal H}_p^\partial]}   
 \int_0^{\infty}     \sqrt{\K_\T}\  {\rm ch}\left({\frac{l_\T}{2}} \right)\exp\left(-2 \sqrt{\K_\T}\ {\rm ch}\left({\frac{l_\T}{2}} \right)\right) \mathcal{E}_{(\bS-\T)_{c'}}(\omega';\K_{\bf r},\L,\K_\T) d\log \K_\T.
 \ee
\begin{figure}[ht]
	\centering
	\includegraphics[scale=0.22]{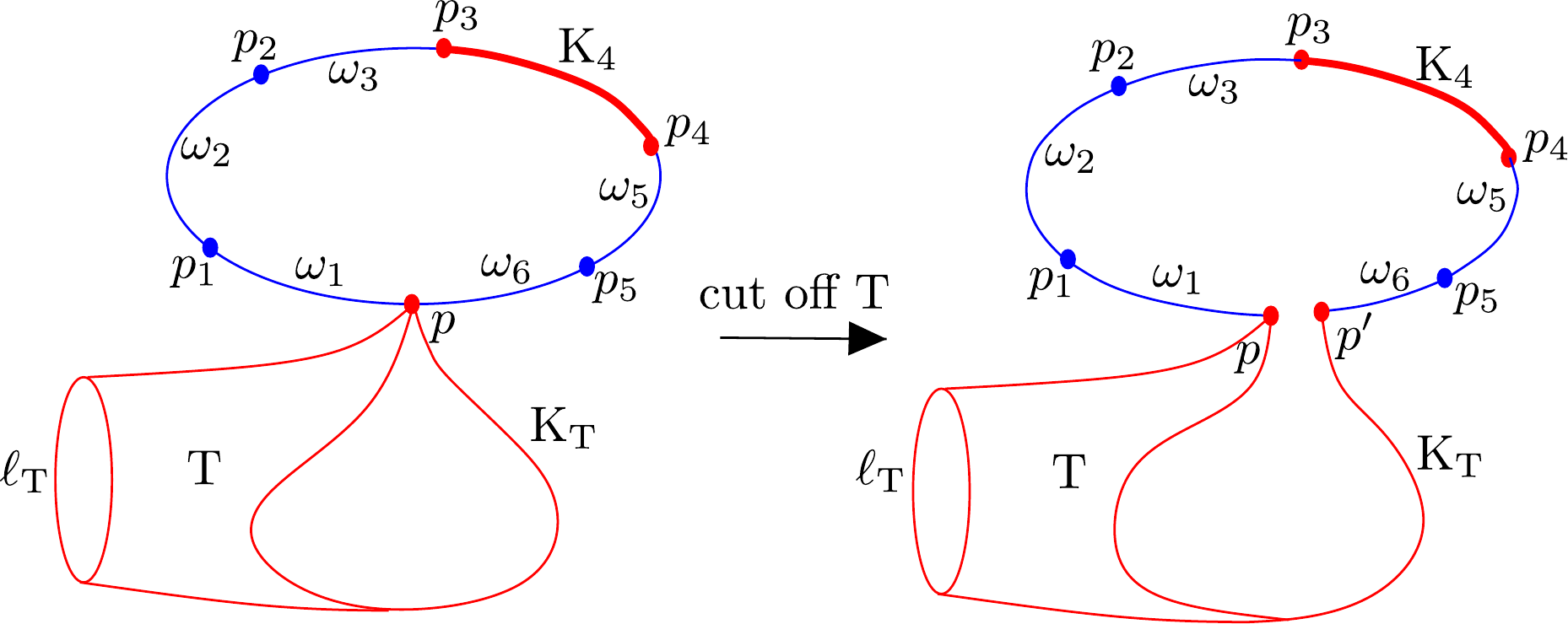}
	\small
	\caption{Cut off a trouser leg.}
	\label{figure:cuttrouser}
\end{figure}
These two terms are the analogs of the first two terms of the recursion \eqref{eq:re}. 
\vskip 2mm

($\tau$) The other terms arise from cutting our a triangle $\tau$ with the cusp $p$, and accounted as follows.
\vskip 1mm

  Suppose that the cusp $p$ lies on the crown end $\D_k^*$. 
Suppose that the other two vertices $x$, $y$ of the triangle $\tau$ 
lie on any crown $\D_a^*$ with $xy \subset \partial \bS$. Denote the $\K$-parameters of the triangle $\tau$ by $\K_1$, $\K_2$, $\K_3$. Then there are two cases:
\begin{enumerate}
\item When we cut off $\tau\in [\Delta_p^\circ]$ from $(\bS_{\rm{c}}; \omega)$ as in Figure \ref{figure:cuttri}, the red cusp $p$ splits into red $p$ and red $p'$, and the two crowns $\D_k^*$ and $\D_a^*$ combine into the crown $\D_{k+a+1}^*$;
\begin{figure}[ht]
	\centering
	\includegraphics[scale=0.2]{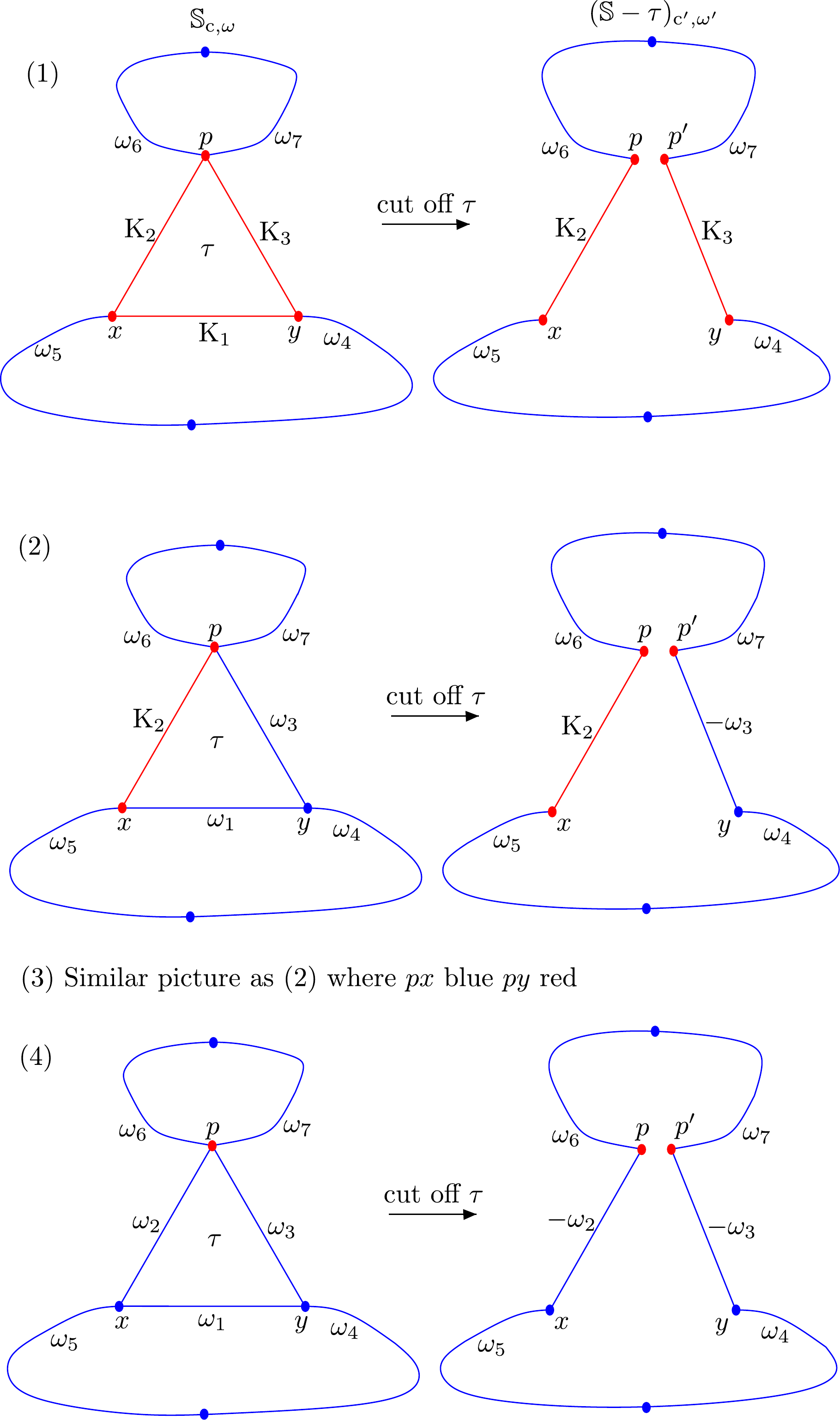}
	\small
	\caption{Cut off an interior triangle.}
	\label{figure:cuttri}
\end{figure}
\item For any triangle $\tau\in [\Delta_p^\partial]$, there is another edge on the boundary of $\bS$. We always assume that $py\subset \partial \bS$. When we cut off $\tau$ from $(\bS_{\rm{c}}; \omega)$ as in Figure \ref{figure:cuttri2}, the crown $\D_k^*$ changes into $\D_{k-1}^*$.
\begin{figure}[ht]
	\centering
	\includegraphics[scale=0.2]{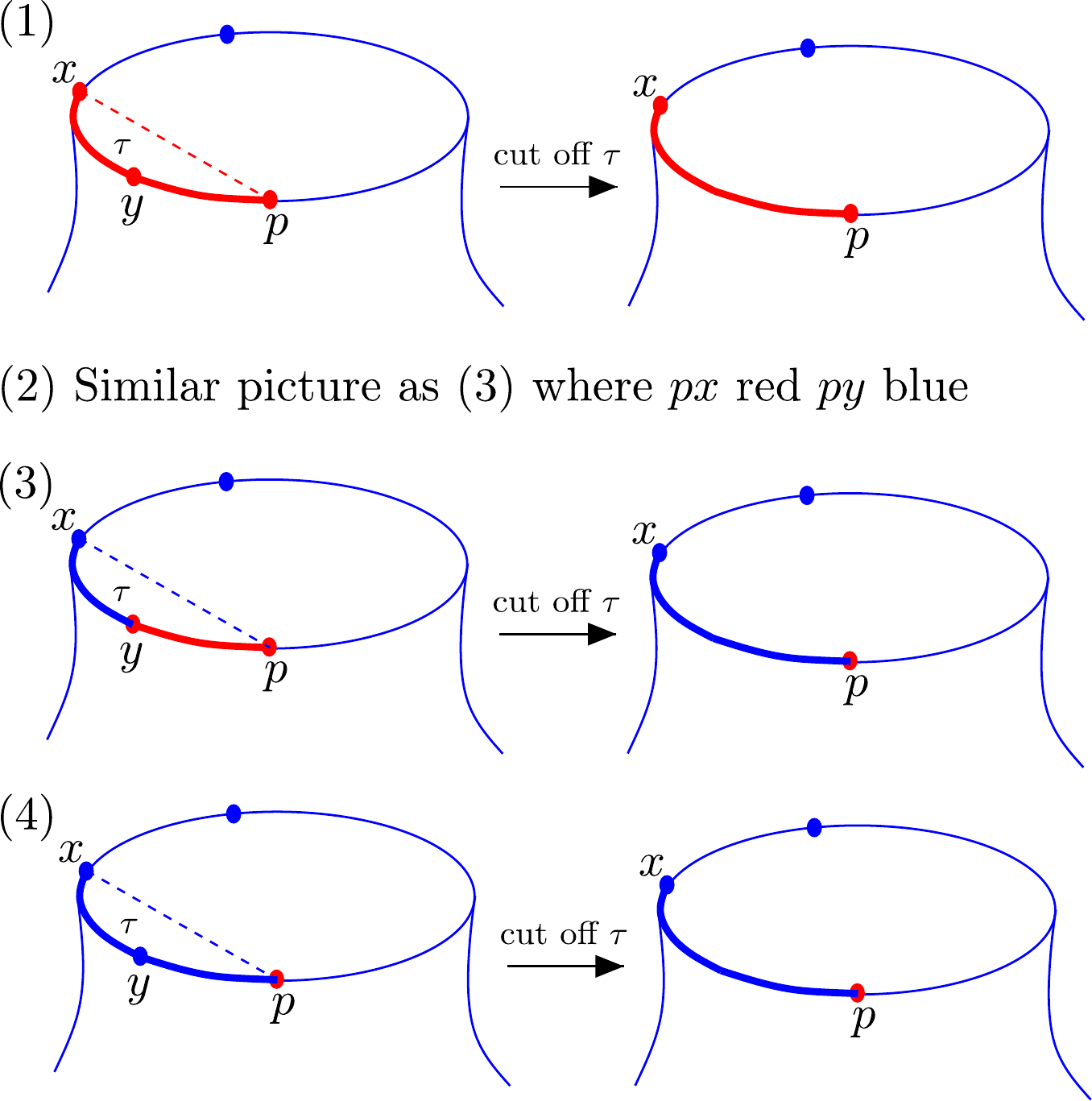}
	\small
	\caption{Cut off a boundary triangle.}
	\label{figure:cuttri2}
\end{figure}
\end{enumerate}
In both two cases, the color $\rm{c}'$ and weight $\omega'$ of the surfaces $\bS-\tau$ is inherited from  $(\bS_{\rm{c}}; \omega)$. There are four possibilities for the  colors of $x$ and $y$:
\begin{enumerate}
\item If both $x$ and $y$ are red as in Figure \ref{figure:cuttri}(1), there is the  weighted exponential volume 
$$
\mathcal{E}_{(\bS-\tau)_{c'}}(\omega';\K_{\bf r}\backslash \{\K_1\},\K_2,\K_3, \L).
$$
Below we will use the shorthand
\be \la{FF}
f(\K_1,\K_2,\K_3):= \sqrt{\frac{\K_2 \K_3}{\K_1}} \exp\left(-\sqrt{\frac{\K_1 \K_2}{\K_3}}- \sqrt{\frac{\K_2 \K_3}{\K_1}}- \sqrt{\frac{\K_3 \K_1}{\K_2}}\right).
\ee
Then the related two terms  in the recursion formula, which are similar to  the third and the fourth terms of  the recursion \eqref{eq:re},  are:
\be \la{S3}
\sum_{[\tau]\in  [\Delta_p^\circ], {\rm x,y \;red}}   
 \int_0^{\infty}  \int_0^{\infty}    f(\K_1,\K_2,\K_3)  \mathcal{E}_{(\bS-\tau)_{c'}}(\omega';\K_{\bf r}\backslash \{\K_1\},\K_2,\K_3,\L) d\log \K_2 \cdot d\log \K_3, 
 \ee
and for Figure \ref{figure:cuttri2}(1)
\be \la{S4}
\sum_{[\tau] \in  [\Delta_p^\partial], {\rm x,y \;red}}   
 \int_0^{\infty}     f(\K_1,\K_2,\K_3) \mathcal{E}_{(\bS-\tau)_{c'}}(\omega';\K_{\bf r}\backslash \{\K_1,\K_3\},\K_2,\L) d\log \K_2.
 \ee

\item Let us suppose that one of $x$ and $y$ is red. 

Consider first the case when  $x$ is red and $y$ is blue as in Figure \ref{figure:cuttri}(2). Then $px$ is red, $py$ and $xy$ are blue. When we cut off $\tau \in [\Delta_p^\circ]$, the integral is almost the same as the last case, except $\K_1$, $\K_3$ are Mellin transformed into $\omega_1$, $\omega_3$. We integrate from $0$ to $+\infty$ for $\K_2$ as before, but the previous integral from $0$ to $+\infty$ for $\K_3$ is changed into integral from $-i\infty$ to $+i\infty$ for $\omega_3$ after Plancherel formula. The weight $py$ for $(\bS_{\rm{c}}; \omega)$ is $0$ since $py$ is not the boundary edge there. Thus the weight for $py$ on $\tau$ is $\omega_3$ while the weight for $p'y$ on $\bS-\tau$ is $-\omega_3$. 

Then the related two terms  in the recursion formula, which are similar to  the third and the fourth terms of  the recursion \eqref{eq:re},  are:
\be \la{S5} 
\sum_{[\tau]\in  [\Delta_p^\circ], {\rm x \;red, \;y \;blue}}   
 \frac{1}{2\pi i} \int_0^{\infty}  \int_{-i\infty}^{+i\infty}    \mathcal{M}(f(\K_2); \omega_1, \omega_3)  \mathcal{E}_{(\bS-\tau)_{c'}}(\omega\backslash \{\omega_1\}\cup\{-\omega_3\};\K_{\bf r},\K_2,\L) d\log \K_2 \cdot d \omega_3, 
 \ee
and for Figure \ref{figure:cuttri}(2)
\be \la{S6}
\sum_{[\tau] \in  [\Delta_p^\partial], {\rm x \;red, \;y \;blue}}   
 \int_0^{\infty}      \mathcal{M}(f(\K_2); \omega_1, \omega_3)  \mathcal{E}_{(\bS-\tau)_{c'}}(\omega\backslash \{\omega_1,\omega_3\};\K_{\bf r},\K_2,\L) d\log \K_2.
 \ee

By the same argument, when $x$ is blue and $y$ is red as in Figure \ref{figure:cuttri}(3) and Figure \ref{figure:cuttri2}(3),  the related two terms  in the recursion formula, which are similar to  the third and the fourth terms of  the recursion \eqref{eq:re},  are:
\be \la{S7} 
\sum_{[\tau]\in  [\Delta_p^\circ], {\rm x \;blue, \;y \;red}}   
\frac{1}{2\pi i} \int_{-i\infty}^{+i\infty} \int_0^{\infty}    \mathcal{M}(f(\K_3); \omega_1, \omega_2)  \mathcal{E}_{(\bS-\tau)_{c'}}(\omega\backslash \{\omega_1\}\cup\{-\omega_2\};\K_{\bf r},\K_3,\L)  d \omega_2 \cdot d\log \K_3, 
 \ee
and 
\be \la{S8}
\sum_{[\tau]\in  [\Delta_p^\partial], {\rm x \;blue, \;y \;red}}   
 \frac{1}{2\pi i} \int_{-i\infty}^{+i\infty}    \mathcal{M}(f(\K_3); \omega_1, \omega_2)  \mathcal{E}_{(\bS-\tau)_{c'}}(\omega\backslash \{\omega_1\}\cup\{-\omega_2\};\K_{\bf r}\backslash \{\K_3\},\L)  d \omega_2.
  \ee
\item Let us suppose now that both $x$ and $y$ are blue as in Figure \ref{figure:cuttri}(4). 

When we cut off $\tau \in [\Delta_p^\circ]$, the $px$, $py$ and $xy$ are all blue. Then $\K_1$, $\K_2$ and $\K_3$ are Mellin transformed into $\omega_1$, $\omega_2$, $\omega_3$. The weights of $px$ and $py$ for $(\bS_{\rm{c}}; \omega)$ are $0$ since $px$ and $py$ are not the boundary edges there. Thus the weights for $px$, $py$ on $\tau$ are $\omega_2$, $\omega_3$ respectively while the weight for $px$, $p'y$ on $\bS-\tau$ are $-\omega_2$, $-\omega_3$. 

Therefore the related two terms  in the recursion formula, matching the third and the fourth terms of  the recursion \eqref{eq:re},  are:
\be \la{S9}
\sum_{[\tau]\in  [\Delta_p^\circ],{\rm x,y \;blue}}   
\frac{1}{(2\pi i)^2} \int_{-i\infty}^{+i\infty}  \int_{-i\infty}^{+i\infty}    \mathcal{M}(f; \omega_1,\omega_2, \omega_3)  \mathcal{E}_{(\bS-\tau)_{c'}}(\omega\backslash \{\omega_1\}\cup\{-\omega_2,-\omega_3\};\K_{\bf r},\L) d \omega_2 \cdot d \omega_3, 
 \ee
and for Figure \ref{figure:cuttri2}(4)
\be \la{S10}
\sum_{[\tau]\in  [\Delta_p^\circ],{\rm x,y \;blue}}   
\frac{1}{2\pi i} \int_{-i\infty}^{+i\infty}   \mathcal{M}(f; \omega_1,\omega_2, \omega_3)  \mathcal{E}_{(\bS-\tau)_{c'}}(\omega\backslash \{\omega_1,\omega_3\}\cup\{-\omega_2\};\K_{\bf r},\L) d \omega_2.
 \ee
\end{enumerate}

We are now ready to state one of the main results of the paper: the recursion formula for the weighted exponential volumes of the moduli spaces related to the weighted colored  surfaces. 
\begin{theorem}
\label{thm:colre}
Given a weighted decorated surface $(\bS_{{\rm c}};\omega)$ and a red cusp $p$, we have the recursion formula written $\int_{{\cal M}_{\bS_{\rm{c}}}(\K_{\bf r}, \L)} W_p {\Bbb E}_{\bS_{\rm{c}}}(\omega;\K_{\bf r}, \L)  $ as the sum of Equations \eqref{S1} - \eqref{S10}. Precisely,
\be
\begin{split}
&\int_{{\cal M}_{\bS_{\rm{c}}}(\K_{\bf r}, \L)} W_p {\Bbb E}_{\bS_{\rm{c}}}(\omega;\K_{\bf r}, \L)  \\
=& \eqref{S1} + \eqref{S2} + \eqref{S3} + \eqref{S4} + \eqref{S5} + \eqref{S6} + \eqref{S7} + \eqref{S8} + \eqref{S9} + 
\eqref{S10}.\\
\end{split}
\ee
\end{theorem}
\bpr
 Follows immediately from Theorem \ref{thm:girre}, Theorem \ref{thm:fc} and the above discussion. 
\epr

\section{Probabilistic moduli spaces} \la{SEC6}

   \subsection{Measures on elementary probabilistic moduli spaces}

\subsubsection{The  measure $\Theta^{\K, \T}_{\sigma^2}$ for the interval.}

Let $\mathbb{W}^\T_{\sigma^2}$ be the {\em Wiener probability measure} on 
$$
{\rm C}_{0}[0,\T]:=\{\xi \in {\rm C}[0,\T]\;|\; \xi(0)=0\}
$$
with the {\it Wiener process} $\mathrm{W}_t(\xi)=\xi(t)$ with zero mean, and covariance
$
\mathbb{E}_{\mathbb{W}^\T_{\sigma^2}}[\mathrm{W}_s \mathrm{W}_t]=\sigma^2 \T \min\{s,t\}$.

Let
\be
\label{eq:aT}
{\rm C}_{0,a}[0,\T]:=\{\xi  \in {\rm C}([0,\T]\rightarrow [0,a])\;|\; \ \xi (0)=0,\xi (\T)=a\}.
\ee

\bd {\rm \cite[Definition 14]{Lo24}}. 
The unnormalised Brownian bridge measure $\mathcal{B}^{a,\T}_{\sigma^2}$ is a measure on ${\rm C}_{0}[0,\T]$ supported on ${\rm C}_{0,a}[0,\T]$ 
such that
\be
\label{eq:WB}
\mathbb{W}^\T_{\sigma^2}(d\xi)=\int_{\mathbb{R}}   \mathcal{B}^{a,\T}_{\sigma^2}(d\xi) da.
\ee
\ed
Consistently with (\ref{eq:WB}) and the fact that the Wiener measure is a probability measure, the total mass of $\mathcal{B}^{a,\T}_{\sigma^2}$ is the Gaussian function
\be
\label{eq:intB}
\int_{{\rm C}_{0,a}[0,\T]} d\mathcal{B}^{a,\T}_{\sigma^2}=\frac{e^{-\frac{a^2}{2\T \sigma^2}}}{\sqrt{2\pi \T }\sigma}:=g_{\sqrt{T}\sigma}(a). 
\ee
Disintegrating the Brownian bridge measure $\mathcal{B}^{a,\T}_{\sigma^2}$, we introduce the measure $\widetilde{\Theta}^{a,b}_{\sigma^2,\T}$.
\bd {\rm \cite[Definition 14]{Lo24}}.  There is a unique measure $\widetilde{\Theta}^{a,b}_{\sigma^2,\T}$ supported on  
$$
\xi(\T)=a, \ \ \ \int_0^\T e^{\xi(s)}ds=b,
$$
such that
\be \la{122}
\mathcal{B}^{a,\T}_{\sigma^2}(d\xi)=\int_0^\infty  \widetilde{\Theta}^{a,b}_{\sigma^2,\T}(d\xi) db.
\ee
\ed

Let
\be \la{124}
\begin{split}
\mathrm{C}(\lambda, \T):=&\{\xi  \in {\rm C}_{0,0}[0,\T]\;|\; \int_0^\T e^{\xi(t)}dt=\lambda/\pi\}\\
=&\{\xi  \in {\rm C}[0,\T]\;|\; \ \ \ \xi (0)=0, \ \xi (\T)=0, \  \int_0^\T e^{\xi(t)}dt=\lambda/\pi\}.\\
\end{split}
\ee

\bd
The { measure $\Theta^{\lambda, \T}_{\sigma^2}$ for an interval} is a measure on $\mathrm{C}(\lambda, \T)$ given by 
\[\Theta^{\lambda, \T}_{\sigma^2}:=\frac{2\lambda}{\pi}\exp\left(\frac{4\pi}{\sigma^2\lambda}\right) \widetilde{\Theta}^{0,\frac{\lambda}{\pi}}_{\sigma^2,\T}.\]
\ed

Note that formula (\ref{122}) just means that the  measure $\Theta^{\lambda, \T}_{\sigma^2}$ satisfies 
\be
\label{eq:btheta}
\mathcal{B}^{0,\T}_{\sigma^2}(d\xi)=\int_0^\infty \frac{d \lambda}{2\lambda} \exp\left(-\frac{4\pi}{\sigma^2\lambda}\right)\Theta^{\lambda, \T}_{\sigma^2}(d\xi).
\ee

By \cite{Yo80, MY05} the total mass of the measure $ \Theta^{\lambda, \T}_{\sigma^2} $ is finite and calculated  via Yor's integral $\theta$:
\be
\label{eq:HWYx}
\begin{aligned}
&\int_{\mathrm{C}(\lambda, \T)} \Theta^{\lambda, \T}_{\sigma^2}=
\frac{1}{2\pi^2} \int_0^\infty \exp\left(-\frac{\tau^2 \sigma^2 \T}{8}\right) J_{-i \tau}\left(\frac{4\pi^2\lambda^2}{\sigma^4 }\right) \tau \sinh(\tau\pi) d\tau= :\theta\left(\frac{4\pi \lambda}{\sigma^2 },\frac{\sigma^2 \T}{4}\right). 
\end{aligned}
\ee

\subsubsection{The measure for the $\varphi-$colored red-based half-disc $\D^{-}_\varphi$.}

Recall the set 
 \be
\label{eq:intdiff}
\mathrm{Diff}^1_+(\lambda, \T):=\{\varphi  \in \mathrm{Diff}^1_+([0,\T]\rightarrow [0,1])\;|\; \varphi(0)=0, \ \varphi(\T)=1, \  \ \varphi '(0)=\varphi '(\T)=\pi/\lambda\}.
\ee
Recall  the {\em Malliavin--Shavgulidze map} \cite{BLW24}:
\be \la{MSM}
 \varphi(t):= \P(\xi)(t):=\frac{\int_0^t e^{\xi(s)} ds}{\int_0^\T e^{\xi(s)} ds}.
\ee

\bl \la{MSH} The  map $\P:\xi \lms \P(\xi)(t)$ provides an isomorphism
\be \la{128}
\begin{split}
\P:& \ \mathrm{C}(\lambda, \T) \stackrel{\sim}{\lra} \mathrm{Diff}_+^1(\lambda, \T).\\
\end{split}
\ee
\el

\bpr If $\xi$ lies  in the support of the  measure $\Theta^{\lambda, \T}_{\sigma^2}$ in  (\ref{124}), 
then we have 
\be 
\label{eq:phieL}
\begin{split}
&\varphi'(0)= \left(\int_0^\T e^{\xi(t)}dt\right)^{-1} = \frac{\pi}{\lambda},\\
&\varphi(t)=\varphi'(0) \int_0^t e^{\xi(t)} dt.\\
\end{split}
\ee
Since 
$\xi(\T)=0$, we get 
$
\varphi'(\T)=\varphi'(0)=\frac{\pi}{\lambda}.$ 
So $
\P(\mathrm{C}(\lambda, \T))\subset \mathrm{Diff}_+^1(\lambda, \T)$.  
The inverse map  is  
\be \la{xiL}
\xi(t):= \log \frac{\varphi'(t)}{\varphi'(0)}.
\ee
Indeed, $\xi(t)$  is a well defined  since $\varphi'(t)>0$,  continuous, and 
by \eqref{eq:phieL}, we recover formula (\ref{MSM}).
\epr

\bd \la{D9.13}
Let $\K:=\lambda^{-2}$. 
The {Schwarzian measure $\Omega^{\K, \T}_{\sigma^2}$ for the interval} is
\be \la{125}
\Omega^{\K, \T}_{\sigma^2} := \P_*\left(\Theta^{\lambda, \T}_{\sigma^2}\right).
\ee
\ed
By (\ref{125}) and (\ref{eq:HWYx}),  the total mass of the Schwarzian  measure $ \Omega^{\K, \T}_{\sigma^2} $ is finite, and given by (\ref{eq:HWYx}).

 \bd \la{D9.13a} The measure for the $\varphi-$colored red-based half-disc $\D^{-}_\varphi$  is given by
\be
{\bf E}_{\D_{\varphi}^{-};\sigma^2}(\T, \lambda^{-2}):= \Omega_{\sigma^2}^{\lambda^{-2}, \T}, \ \ \ \ \lambda:=\omega(\A_1, \A_2), 
\ee
 \ed

\subsubsection{Schwarzian measure and  formal Schwarzian action for the interval.} The formal  Schwarzian action for the interval is defined by 
\be
\mathcal{S}(\varphi):=-\frac{1}{2\sigma^2} \int_0^\T \left(\frac{\varphi''(t)}{\varphi'(t)}\right)^2 dt.
\ee
If the function $\xi(t)$ in (\ref{xiL}) is smooth, then $\xi'(t)=\frac{\varphi''(t)}{\varphi'(t)}$, and the formal Schwarzian action  for the interval can be written 
as 
\be \la{133a}
\mathcal{S}(\varphi)=-\frac{1}{2\sigma^2} \int_0^\T \xi'(t)^2 dt.
\ee
The formal Schwarzian action for the disc case was defined in \cite{SW17}.  
Since the Schwarzian of $\varphi$ is 
\be
\label{eq:schf}
\mathrm{Sch}(\varphi,t)=\left(\frac{\varphi''(t)}{\varphi'(t)}\right)'-\frac{1}{2} \left(\frac{\varphi''(t)}{\varphi'(t)}\right)^2,
\ee
we can write the formal Schwarzian action for the  interval  as 
\be
\mathcal{S}(\varphi)=\frac{1}{\sigma^2}  \int_0^\T \mathrm{Sch}(\varphi,t) dt -\frac{1}{\sigma^2}  \left[\frac{\varphi''(t)}{\varphi'(t)}\right]_0^\T. 
\ee
Then the  Schwarzian measure for the interval is  given formally by $e^{-\mathcal{S}(\varphi)}{\cal D}\varphi$. 
This expression is just a motivation since  both factors are ill defined.  Yet we can define their product as 
  the Schwarzian measure $\Omega^{\K, \T}_{\sigma^2}$, similarly  to the measure defined in \cite{BLW24} for the disc,   turning this into a definition: 
  $$
  e^{-\mathcal{S}(\varphi)}{\cal D}\varphi:=   \Omega^{\K, \T}_{\sigma^2}.
  $$

\subsubsection{Measure for the cylinder ${\rm C}_\varphi$ with a single $\varphi-$colored boundary circle.}

We start from  the space of $(l, \T)-$quasiperiodic continuous functions $\xi(t)$ on $\R$, vanishing at $t=0$:
 \be
 \mathrm{C}_{\rm qp}(\R; l, \T):= \{\xi \in \mathrm{C}(\R)\;|\; \ \ \  \xi(t+\T)=\xi(t)+l, \ \ \ \forall t \in \R, \ \ \ \xi(0)=0\}.
\ee
Next, consider  the space of $(l, \T)-$quasiperiodic  monotone ${\rm C}^1$ functions $\varphi(t)$,   vanishing at $t=0$:
 \be \la{160}
 \mathrm{Diff}_{\rm qp}^1(\R; l, \T)_+:= \{\varphi \in \mathrm{C}^1(\R)\;|\; \ \ \  \varphi(t+\T)=\varphi(t)+l,\ \ \ \varphi'(t)> 0 \ \ \ \forall t \in \R; \ \ \ \varphi(0)=0\}.
\ee
Similarly to Lemma \ref{MSH}, there is an isomorphism between these two spaces defined as follows. Set 
\be \la{161}
\varphi(t):=\mathcal{R}_l(\xi)(t)=\log \left(1+(e^l-1)\frac{\int_0^t e^{\xi(s)} ds }{\int_0^\T e^{\xi(s)} ds } \right).
\ee
  

\bl
The map $\mathcal{R}_l$ provides an isomorphism of spaces
\be
\begin{split}
\mathcal{R}_l:& \ \mathrm{C}_{\rm qp}(\R; l, \T) \stackrel{\sim}{\lra} \mathrm{Diff}_{\rm qp}^1(\R; l, \T)_+\\
&\xi \lms \varphi(t):=\mathcal{R}_l(\xi)(t). \\
\end{split}
\ee
\el
\bpr
Since the function $\xi(t)$ is $(l,\T)-$quasiperiodic,  the function $\varphi(t)$ also $(l, \T)-$quasiperiodic. 
The function $\varphi(t)$ vanishes at $t=0$,  and its derivative is positive and continuous, since 
$$
\varphi'(t)=\frac{(e^l-1) e^{\xi(t)}}{\int_0^\T e^{\xi(s)} ds + (e^l-1)\int_0^t e^{\xi(s)} ds}.
$$
So  
$$
\frac{\varphi'(t)}{\varphi'(0)}=e^{\xi(t)}    \Bigl(1+ (e^l -1)\frac{\int_0^t e^{\xi(s)}ds}{\int_0^\T e^{\xi(s)} ds}\Bigr)^{-1} = e^{\xi(t) - \varphi(t)}. 
$$
Therefore
we arrived at the inverse map:  
\be \la{133}
\xi(t)=\varphi(t)+\log \varphi'(t)- \log \varphi'(0).
\ee
We use $\varphi'(t)>0$ to define $\log\varphi'(t)$. The quasiperiodicity of $\xi(t)$ follows from the one of $\varphi(t)$.
\epr

The space  ${\rm C}_{0,l}[0,\T]$ of continuous functions $\xi(s)$ on the interval $[0,\T]$ such that $\xi(0)=0, \xi(\T)=l$, see (\ref{eq:aT}), 
  is  identified with the space $ \mathrm{C}_{\rm qp}(\R; \T, l)$ of $(l,\T)-$quasiperiodic continuous functions $\xi(s)$ on $\R$ vanishing at $s=0$:
\be \la{CANI}
\eta: {\rm C}_{0,l}[0,\T] \stackrel{\sim}{\lra}  \mathrm{C}_{\rm qp}(\R; l, \T).
 \ee
Precisely, we extend uniquely a  function $\xi(s)$ on the interval $[0,\T]$ to a $(l,\T)-$quasiperiodic function on  $\R$. 
The condition $\xi(\T)=l$ implies that the extension of a continuous function is continuous. 

Denote by $\mathcal{B}_{\sigma^2}(\R; l, \T)$ the  measure on the space $\mathrm{C}_{\rm qp}(\R; l, \T)$ obtained as the image of the Brownian bridge measure $\mathcal{B}^{l,\T}_{\sigma^2}$ on ${\rm C}_{0,l}[0,\T] $ under this isomorphism:
\be \la{MBB}
\mathcal{B}_{\sigma^2}(\R; \l, \T) = \eta_* \mathcal{B}^{l,\T}_{\sigma^2}.
\ee
\bd \la{D9.14}
The measure for a cylinder ${\rm C}_{\varphi}$ with the length $l$ neck loop and a $\varphi-$colored boundary is a measure  on $\mathrm{Diff}^1_{\rm qp}({\R; l, \T})_+$ obtained by applying the map ${\cal R}_l$ to the measure (\ref{MBB}): 
\be \la{DSG}
{\bf E}_{{\rm C}_{\varphi}, \sigma^2}(l, \T):=4 \mathcal{R}_{l_*}\Bigl(\mathcal{B}_{\sigma^2}(\R; l, \T)\Bigr).
\ee
\ed

\begin{remark}
{\rm Suppose $\varphi$ and $\xi$ are smooth,  quasiperiodic, 
and related by (\ref{133}), which implies 
$\xi'(t)=\varphi'(t)+\frac{\varphi''(t)}{\varphi'(t)}.$ 
Then the  Schwarzian action is derived from the  Dirichlet action (\ref{133a}):
\be
\begin{aligned}
\mathcal{S}(\varphi)&=-\frac{1}{2\sigma^2}\int_0^\T \xi'(t)^2 dt
\\ & \ = \ -\frac{1}{2\sigma^2} \int_0^\T \left(\frac{\varphi''(t)}{\varphi'(t)} \right)^2 dt-\frac{1}{2\sigma^2} \int_0^\T \varphi'(t)^2 dt
\\& = \ \frac{1}{\sigma^2} \int_0^\T \mathrm{Sch}(e^\varphi,t) dt.
\end{aligned}
\ee
The last equality follows from \eqref{eq:schf},  the chain rule $\mathrm{Sch}(e^\varphi,t)=\mathrm{Sch}(\varphi,t)-\frac{1}{2}\varphi'(t)^2$, and  the quasiperiodicity of $\varphi$, which implies that $\varphi^{(k)}(0)=\varphi^{(k)}(\T)$ for $k\geq 1$.}

\end{remark}

\subsubsection{Measure for a single-cusped $\varphi-$colored crown $\D^*_{\varphi, 1}$}

Recall the  measure   for the $\varphi-$decorated red-based half disc $\D_{\varphi}^{-}$  from Definition \ref{D9.13a}:
\be
{\bf E}_{\D_{\varphi}^{-};\sigma^2}( \lambda, \T):= \Omega_{\sigma^2}^{\lambda^{-2}, \T}.
\ee

Recall that the exponential volume form for the triangle $\tau$ with the local potentials at the red points $\A_1, \A_2$ is given by the  function 
$$
{\cal E}_\tau(\lambda^{-2},e^l):= \exp({-\lambda^{-1}(e^{l/2}+e^{-l/2})}),  
$$
where 
$$
l:= 2\log \frac{\omega(\A_1, \widetilde \A)}{\omega(\A_2, \widetilde \A)}, \ \ \ \lambda:= \omega(\A_1, \A_2).
$$

\bd \la{D9.17}  The  measure ${\bf  E}_{\D_1^*; \sigma^2}(\T, e^l)$ for  the  $\varphi-$colored crown with a single red cusp and the length $l$ neck geodesic is given by
\be \la{RHS}
\begin{split}
{\bf E}_{\D^*_{1}; \sigma^2}(\T, e^l)\stackrel{\mbox{{\rm Figure }\ref{figure:glutrie}}}{:=} 
&2\int {\Omega}^{\lambda^{-2}, \T}_{\sigma^2} \cdot \exp({-\lambda^{-1}(e^{l/2}+e^{-l/2})})d\log \lambda.\\
\end{split}
\ee
\ed

Definition \ref{D9.17} is motivated by the following geometric argument, illustrated by Figure \ref{figure:glutrie}.
\begin{figure}[ht]
	\centering
	\includegraphics[scale=0.7]{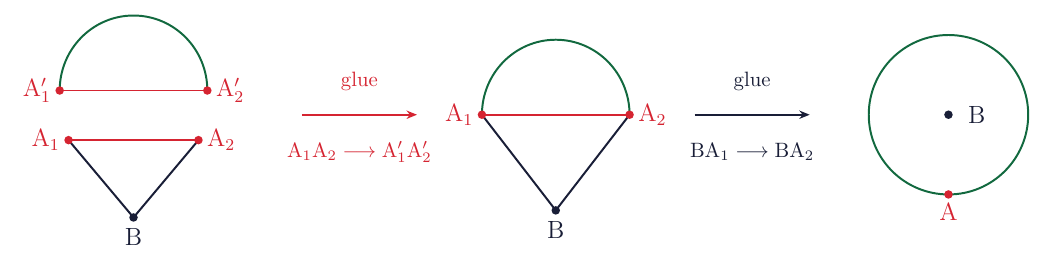}
	\small
	\caption{Glueing the $\varphi-$decorated red-based half-disc to a red-based triangle, and then gluing the two blue sides of the obtained surface, resulting a one-cusped $\varphi-$decorated crown.}
	\label{figure:glutrie}
\end{figure}

Take a triangle $\tau = (\B, \A_1, \A_2)$ with the red base $(\A_1, \A_2)$ and a blue vertex $\B$. Glue it to the red-based half-disc, 
$\varphi-$decorated along the green boundary semi-circle $(A_1', \A_2')$, see Figure \ref{figure:glutrie}. The resulting $\varphi-$decorated surface  is shown in the 
middle of Figure \ref{figure:glutrie}. Gluing of the red diameter $(\A_1, \A_2)$ to the red base $(\A'_1, \A'_2)$ requires  $\omega(\A_1, \A_2) = \omega(\A'_1, \A'_2)$. It amounts to the integral of the product of the exponential volume forms of the triangle and  of 
the red based half-discs over $d\log \lambda$, producing the exponential volume form for the surface in the middle of Figure \ref{figure:glutrie}. Finally, 
glueing  the   boundary intervals $\B\A_1$ and $\B\A_2$ on the surface  in the middle, we get 
a once-punctured disc  $\D^*_{1, \varphi}$ with the $\varphi-$colored boundary,   a single red cusp $\A$, and the puncture $\B$. This amounts to identify $2\log \frac{\omega(\A_1, \widetilde \A)}{\omega(\A_2, \widetilde \A)}$ with the neck geodesic length $l$.

 Alternatively, and equivalently to the above construction, Definition \ref{D9.17} amounts to applying the cutting and gluing formula for the red circle starting in $\A$ and containing the puncture $\B$.

 \subsubsection{Reduction formula for the  measure for  $n-$cusped  $\varphi-$colored crown $\A_{0,n}$.} \la{Sec9.4.5}

\bp \la{D9.18} The  measure ${\bf E}_{\D_{\varphi, n}^*}( \T_1, ..., \T_n, e^l)$ for  the $\varphi-$colored crown with $n$ red cusps and the length $l$ neck geodesic is reduced to the $n=1$ case by 
\be \la{RHS}
{\bf E}_{\D_{\varphi, n}^*}(\T_1, ..., \T_n,e^l):= {\bf E}_{\D_{\varphi, 1}^*}(\T_1+\ldots +\T_n, e^l). 
\ee
 Here we omit $\sigma^2$ from the notation. 
 \ep
 
In other words, the reduction formula tells that 
removing a cusp from a $\varphi-$decorated crown with $n\geq 2$ red cusps amounts to taking the sum of the adjacent time parameters, e.g.: 
$$
{\bf E}_{\D_{\varphi, n}^*}( \T_1,\T_2,  ..., \T_n,e^l) = {\bf E}(_{D_{\varphi, n-1}^*}( \T_1+\T_2,  ..., \T_n, e^l).
 $$ 
 
 \bpr 
 
 Given a triangle $\Delta=pqr$ as in Figure \ref{figure:semicircle} with the $\lambda-$lengths $\lambda_1$, $\lambda_2$, $\lambda$ for $pr$, $pq$, $qr$,  the potential is  $W_\Delta=\frac{\lambda}{\lambda_1 \lambda_2}+\frac{\lambda_1}{\lambda \lambda_2}+\frac{\lambda_2}{\lambda_1 \lambda}$. Let $\T=\T_1+\T_2$. Then the Schwarzian measure $\Theta^{\lambda, \T}_{\sigma^2}$ has the gluing property 
 similar to  Proposition \ref{prop:cge}, which  
   follows from \cite[Prop. 15 and Lemma 18]{Lo24}:
\be
\label{eq:tde}
\Theta^{\lambda, \T}_{\sigma^2}= \int_0^\infty \int_0^\infty e^{-\frac{2\pi}{\sigma^2}W_\Delta} d\Theta^{\lambda_1, \T_1}_{\sigma^2} d\Theta^{\lambda_2, \T_2}_{\sigma^2} d\log \lambda_1 d\log\lambda_2. 
\ee

 \begin{figure}[ht]
	\centering
	\includegraphics[scale=0.7]{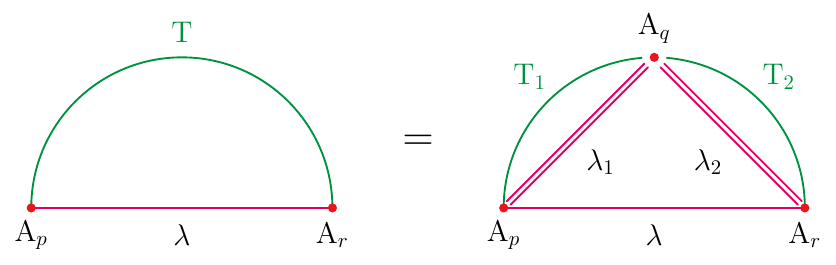}
	\small 
	\caption{Reducung two discs glued to a triangle to red-based half disc.}
	\label{figure:semicircle}
\end{figure}
 
 Then we apply the gluing argument illustrated on Figure \ref{figure:glutrie}, with the   the red-based half-disc with the base $\A_1'\A_2'$ on the left replaced by the half-disc with the same base and $n-1$ red cusps. 
Here we use the following  factorization formula for the  exponential volume for any disc $\D_m$. 

 Consider any  triangulation $\tau$ of the disc $\D_m$. Then we can express the exponential volume $ {\cal E}_{{\rm D}_{m}}$ for the disc $\D_m$ via integrals of  the exponential volumes  ${\cal E}_{\Delta}$ for the triangles of the triangulation. 
The result is independent of the triangulation. For example, for the two triangulations   of  $\square$  we have 
\be \la{sec9.4.5a}
\begin{split}
 {\cal E}_\square(\K_1, \K_2, \K_3, \K_4) = &\int  {\cal E}_\Delta(\K_1, \K_2, \K)  {\cal E}_\Delta(\K, \K_3, \K_4) d\log \K\\
 =&\int  {\cal E}_\Delta(\K_2, \K_3, \K')  {\cal E}_{\Delta}(\K', \K_4, \K_1) d\log \K'.
 \\
\end{split}
 \ee
Applying the factorization formula we get the result. 
 \epr

\subsection{The  measure on  the probabilistic moduli space  ${\cal M}_{\bS_\varphi}$}


There are the following three types of projections associated with a $\varphi-$colored surface $\bS_\varphi$:

\begin{enumerate}

\item The canonical projection
$$
\pi_{\bS_\bullet}: \mathcal{M}_{\bS_\varphi}(\T_\bullet, \K_{\bf r}, \L_{\bf r})\lra {\cal M}_{\bS_\bullet}(\K_\bullet, \L_\bullet).
$$
The moduli space ${\cal M}_{\bS_\bullet}(\K_\bullet, \L_\bullet)$ carries the exponential volume form with the potential  scaled by $\beta=2\pi/\sigma^2$ and denoted by:
$$
{\Bbb E}^{\beta}_{\bS_\bullet}(\K_\bullet, \L_\bullet). 
$$
The  pull back of this exponential volume by the map $\pi_{\bS_\bullet}$ gives us a differential form
\be \la{m1}
\pi^*_{\bS_\bullet}\Bigl({\Bbb E}^\beta_{\bS_\bullet}(\K_\bullet, \L_\bullet)\Bigr).
\ee

\item For each $\varphi-$colored edge $e$ of $\bS$, there is a canonical projection
$$
\pi_{e}: \mathcal{M}_{\bS_\varphi}(\T_\bullet, \K_{\bf r}, \L_{\bf r})\lra \mathrm{Diff}_+^1(\lambda_e,\T_e).
$$
So we can pull back the measure from Definition \ref{D9.13a}, multiplied by $d \log \K_e$:
\be \la{m3}
\pi^*_e\left({\bf E}_{\D_{\varphi}^{-};\sigma^2}(\K_e, \T_e)\wedge d \log \K_e\right).
\ee

\item For each $\varphi-$colored boundary component $S^1_\gamma$  
there is a canonical projection
$$
\pi_{\gamma}: \mathcal{M}_{\bS_\varphi}(\T_\bullet, \K_{\bf r}, \L_{\bf r})\lra  \mathrm{Diff}_{\rm qp}^1(\R;  l_\gamma, \T_\gamma)_+\times \R.
$$
Recall Definition \ref{D9.14}. We pull back $\mbox{(measure (\ref{DSG})} \wedge d l_\gamma)  \wedge d\theta_\gamma$, where $\theta_\gamma$ is the twist parameter along $\gamma$, by the map $\pi_\gamma^*$:
\be \la{m2}
\pi^*_\gamma\Bigl({\bf E}_{{\rm C}_\varphi, \sigma^2}(l_\gamma, \T_\gamma) \wedge d l_\gamma \Bigr) \wedge d\theta_\gamma.
\ee
\end{enumerate}


\bd  \la{D6.20}  The measure ${\bf E}_{\bS_\varphi}$ on $\mathcal{M}_{\bS_\varphi}(\T_\bullet, \K_{\bf r}, \L_{\bf r})$ is  the product of  the measures (\ref{m1}) - (\ref{m2}):
\be
\begin{split}
&{\bf E}_{\bS_\varphi}(\T_\bullet, \K_{\bf r}, \L_{\bf r})\\
:=
& \pi^*_{\bS_\bullet}\left({\Bbb E}^\beta_{\bS_\bullet}(\K_\bullet, \L_\bullet)\right) \wedge \bigwedge_{e\in {\cal I}_\varphi} \pi^*_e\left({\bf E}_{\D_{\varphi}^{-};\sigma^2}(\K_e, \T_e) \wedge d \log \K_e\right)\wedge \bigwedge_{\gamma\in {\cal J}_\varphi} \pi^*_\gamma\Bigl({\bf E}_{{\rm C}_\varphi, \sigma^2}(l_\gamma, \T_\gamma) \wedge dl_\gamma\Bigr)\wedge d\theta_\gamma.\\
\end{split}
\ee
Then the volume of the {measure ${\bf E}_{\bS_\varphi}$}
\be
\mathcal{E}_{\bS_\varphi}(\T_\bullet, \K_{\bf r}, \L_{\bf r}) := \int_{\mathcal{M}_{\bS_\varphi}(\T_\bullet, \K_{\bf r}, \L_{\bf r})} {\bf E}_{\bS_\varphi}(\T_\bullet, \K_{\bf r}, \L_{\bf r}).
\ee
\ed

\bp The  measure ${\bf E}_{\bS_\varphi}$ for a $\varphi-$decorated surface $\bS_\varphi$ containing a crown ${\rm C}$ with $n$ cusps and  time parameters 
$\T_1, ..., \T_n$ at the boundary intervals on the crown is equal to the  one for the $\varphi-$decorated surface $\bS'_\varphi$ obtained from $\bS_{\varphi}$ by replacing the  crown ${\rm C}$ by a
single-cusped  crown 
${\rm C}'$ with the time parameter 
$\T_1+ \ldots +  \T_n$. So,  suppressing the rest of the parameters which remain intact: 
$$
{\bf E}_{\bS_\varphi}(\T_1,\T_2,  \ldots , \T_n) = {\bf E}_{\bS_{\varphi'}}( \T_1+\ldots +\T_n).
 $$ 
\ep

\bpr Follows from Definition \ref{D6.20} and Proposition \ref{D9.18} by using the neck recursion formula. 
\epr 

\subsection{Volumes of probabilistic moduli spaces  ${\cal M}_{\bS_\varphi}(\T_\bullet, \K_{\bf r}, \L_{\bf r})$ are finite}

\bt   \la{Th6.14} The volume of the  measure ${\bf E}_{\bS_\varphi}$ on $\mathcal{M}_{\bS_\varphi}(\T_\bullet, \K_{\bf r}, \L_{\bf r})$ is finite.
\et

\bpr

Let $\bS_\varphi$ be a $\varphi-$decorated surface 
 obtained by gluing   the punctured disc $\D_\N^*$  with $\N$  red boundary intervals and the  red-based half-discs $\D_{e,\varphi_e}$ along some of the boundary intervals ${\rm I}_e$:
$$
\bS_\varphi=\D^*_\N \cup \bigcup_{e\in \mathcal{I}_\varphi} \D_{e,\varphi_e}, \ \ \ \ k:= |\mathcal{I}_\varphi|\leq \N.
$$

 Recall 
 $
 e^l=\L$, $\lambda_i=\K_i^{-\frac{1}{2}}$, and $\beta=\frac{2\pi}{\sigma^2}.
 $ 
  Let $f(t):=e^{-\frac{2}{t}} \sqrt{\pi t }$. Recall  the Gaussian 
 \be \la{Ga}
 g_a(t)=\frac{1}{\sqrt{2\pi a}}e^{-\frac{t^2}{2a}}.
 \ee. 
 Let 
 $
 \tau_e:=\sigma^2 \T_e/4$ and  $x:=\sum_{e\in \mathcal{I}_\varphi} \tau_e +\sum_{i\notin \mathcal{I}_\varphi} \lambda_i/(2\beta).
 $
\bp
\label{lem:crowneq}
The  volume of the moduli space  $\mathcal{M}_{\bS_\varphi}(\T_\bullet, \K_{\bf r}, \L_{\bf r})$ satisfy the  inequality
\[\mathcal{E}_{\bS_\varphi}(\T_\bullet, \K_{\bf r}, \L_{\bf r})\leq 2^{k+1}  \prod_{i\notin \mathcal{I}_\varphi}  f(\lambda_i/\beta) g_{x}\left(l/2\right).\]
\ep
\bpr
For any $e\in \mathcal{I}_\varphi$, recall Yor's integral  \eqref{eq:HWYx}. We have 
\be \la{FW}
\mathcal{E}_{\D_{e,\varphi_e}}(\T_e,\K_e)=\theta\left(2\beta/\lambda_e, \tau_e\right).
\ee
Given an ideal triangulation of $\D^*_{\N}$ as in Figure \ref{figure:D13}, we use the notations as in Proposition \ref{Meij}. Then 
\be \la{FW}
\prod_{i=1}^\N \X_i=\frac{\prod_{i=1}^\N \K_i}{\left(\prod_{i=1}^\N \B_i\right)^2}=e^l,\;\;\; W_{\D^*_\N}= \sum_{i=1}^\N \left(\B_i +\frac{\K_i}{\B_i}\right).
\ee
Let $s_i:= \log (\B_i \lambda_i)$. Then we can rewrite formulas (\ref{FW}) as follows:
\be \la{HS}
\sum_{i=1}^\N s_i=-l/2,\;\;\;W_{\D^*_\N}=\sum_{i=1}^\N \frac{2\cosh s_i}{\lambda_i}. 
\ee
By \cite[Corollary 2.21]{GS24}, the cluster volume form for the punctured disc $\D_\N^*$, restricted to the hypersurface  given by the first equation in (\ref{HS}),  is
\[2 d \log \X_1...d \log \X_{\N-1} d \log \B_1... d \log \B_{\N}=2^{\N+1} d \log \lambda_1...d \log \lambda_\N d s_1... d s_{\N-1}.\]
Then using for each edge $e\in {\cal I}_\varphi$   the gluing formula we get 
\be
\begin{aligned}
\label{eq:compdnst}
\mathcal{E}_{\bS_\varphi}(\T_\bullet, \K_{\bf r}, \L_{\bf r})&=2^{k+1} \int \exp\left(-\beta\sum_{i=1}^\N \frac{2\cosh s_i}{\lambda_i}\right) \prod_{e\in \mathcal{I}_\varphi} \Bigl(\theta (2\beta/\lambda_e, \tau_e)
d \log \lambda_e \Bigr)\cdot d s_1... d s_{\N-1}
\\&=2^{k+1} \int \prod_{e\in \mathcal{I}_\varphi} \left(\int\exp\left(-2\beta \frac{\cosh s_e}{\lambda_e}\right)\theta\left(2\beta/\lambda_e, \tau_e\right) 
d \log \lambda_e \right) \\
&\times \prod_{i\notin \mathcal{I}_\varphi} \exp\left(-2\beta \frac{\cosh s_i}{\lambda_i}\right)  d s_1... d s_{\N-1}.
\end{aligned}
\ee

\bl \la{JL} For any $e \in \mathcal{I}_\varphi$, we have
\be
\int_0^\infty \exp\left(-2\beta \frac{\cosh s_e}{\lambda_e}\right)\theta\left(2\beta/\lambda_e, \tau_e\right) d \log \lambda_e =
g_{\tau_e}(s_e).
\ee
\el 

\bpr This follows from  \cite[Formula 4.6]{MY05} and Equation \eqref{eq:intB}.  
\epr

For any $i \notin \mathcal{I}_\varphi$, we have
\be
\exp\left(-2\beta \frac{\cosh s_i}{\lambda_i}\right)\leq \exp\left(-2\beta \frac{1+s_i^2/2}{\lambda_i}\right)
=f(\lambda_i/\beta) g_{\frac{\lambda_i}{2\beta}}(s_i).
\ee
Then using Lemma \ref{JL}, the right hand side of Equation \eqref{eq:compdnst} is
\be
\begin{aligned}
&\leq  2^{k+1}  \prod_{i\notin \mathcal{I}_\varphi}  f(\lambda_i/\beta) \int_{\sum_{i=1}^\N s_i =-l/2}  \prod_{e\in \mathcal{I}_\varphi} g_{\tau_e}(s_e) \prod_{i\notin \mathcal{I}_\varphi} g_{\frac{\lambda_i}{2\beta}}(s_i)  d s_1... d s_{\N-1}
\\&=  2^{k+1}  \prod_{i\notin \mathcal{I}_\varphi}  f(\lambda_i/\beta) g_{x}\left(l/2\right). 
\end{aligned}
\ee
The last equality follows by calculating the  convolution of  Gaussians via the Fourier transform. 
\epr

Let $\mathcal{N}$ be the collection of the neck geodesics along the crown ends. For any $\gamma\in \mathcal{N}$, let $\varphi_{\gamma}$ be the collection of $\phi_e$ on the crown $\D^*_{\gamma}$. 
By cutting along all the neck geodesic circles, we get
\[\bS_\varphi= S_{g,m}\cup \bigcup_{\gamma \in \mathcal{J}_\varphi} \D^*_{\gamma,\varphi_\gamma} \cup \bigcup_{\gamma \in \mathcal{N}} \D^*_{\gamma,\varphi_{\gamma}}.\]
Then by  \eqref{eq:M07}, the  volume $\mathcal{E}_{\bS_\varphi}(\T_\bullet, \K_{\bf r}, \L_{\bf r})$ is calculated by the following sum: 
\be
\label{eq:totmas1}
\sum_{0\leq d_1+ \cdots + d_m \leq d} {\cal V}_{d_1, \cdots, d_m}\prod_{i=1}^m l_i^{2d_i}  \prod_{\gamma \in \mathcal{J}_\varphi}  \mathcal{E}_{\D^*_{\gamma,\varphi_\gamma}}(\T_\gamma,\L_\gamma)  l_\gamma d l_\gamma \prod_{\gamma \in \mathcal{N}}
\mathcal{E}_{\D^*_{\gamma,\varphi_{\gamma}}}(\{\T_{e}\}_{e\in \D^*_{\gamma}},\K_{\gamma},\L_\gamma)  l_\gamma d l_\gamma.
\ee
When $\gamma \in \mathcal{J}_\gamma$, by  \eqref{eq:intB} we have
$\mathcal{E}_{\D^*_{\gamma,\varphi_\gamma}}(\T_\gamma,\L_\gamma)=4g_{\sqrt{T_\gamma}\sigma}(l_\gamma)$.  
Therefore for any $n\geq 0$, we have
\be
\label{eq:Dgamma}
\int_0^\infty l_\gamma^n \mathcal{E}_{\D^*_{\gamma,\varphi_\gamma}}(\T_\gamma,\L_\gamma) d l_\gamma= 4\int_0^\infty l_\gamma^n g_{\sqrt{T_\gamma}\sigma}(l_\gamma)d l_\gamma<+\infty.
\ee
When $\gamma \in \mathcal{N}$, by Lemma \ref{lem:crowneq} we obtain
\be
\label{eq:Ngamma}
\int_0^\infty l_\gamma^n   \mathcal{E}_{\D^*_{\gamma,\varphi_{\gamma}}}(\{\T_{e}\}_{e\in \D^*_{\gamma}},\K_{\gamma},\L_\gamma) d l_\gamma \leq  2^{k+1}  \prod_{i\notin \mathcal{I}_\varphi}  f(\lambda_i/\beta) \int_0^\infty  l_\gamma^n  g_{x}\left(l/2\right) d l_\gamma <+\infty.
\ee
Then plugging  \eqref{eq:Dgamma}-\eqref{eq:Ngamma} into \eqref{eq:totmas1}, we conclude that the  volume is finite.
\epr

\section{Recursions for   probabilistic moduli spaces  ${\cal M}_{\bS_\varphi}$} \la{Sect9}

In Section \ref{Sect9} we think about  Teichm\"uller spaces using the language  recalled in Section \ref{Sect5.3.3}. 
So we talk about flags $\B$ and decorated flags $\A$, which mean points on $\R\P^1$ and vectors in $\R^2-\{0\}$, as   in Section \ref{section:cpAF}. 
\subsection{The triangle case} \la{Sect7.1}

 \begin{figure}[ht]
	\centering
	\includegraphics[scale=0.3]{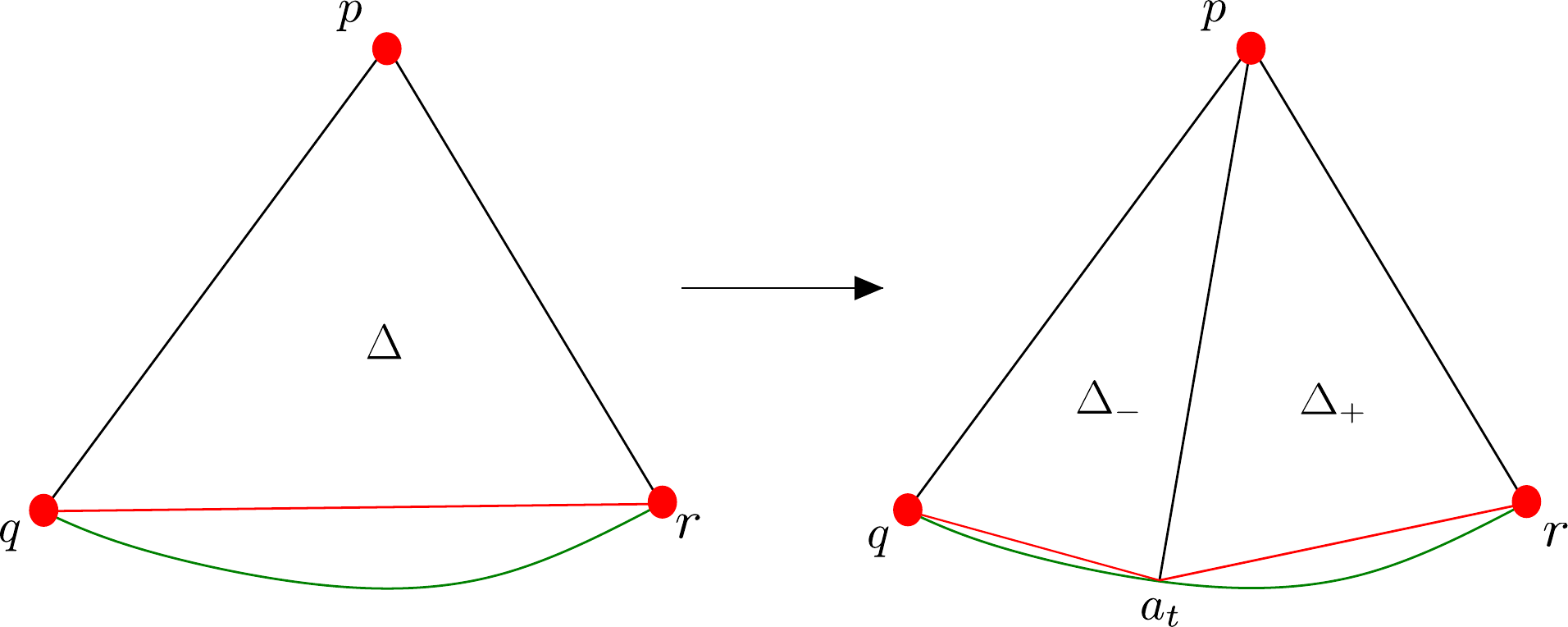}
	\small
	\caption{Cutting a triangle with a $\varphi-$colored   base  along the edge $p a_t$.}
	\label{figure:recursiontr}
\end{figure}

From the geometric perspective of Figure \ref{figure:recursiontr}, in Sections \ref{Sect7.1}-\ref{YK} we pass from the red-based  half-disc ${qr}$ with the green arc on the left, to the concatenation of  similar red-based  
half-discs $qa_t$ and $a_tr$ on the right. 
Then we attach a triangle to each Brownian bridge (= red-based half-disc with green arc), and do a flip of the triangulation of the red rectangle on the right.   Let us implement this in formulas. 
\vskip 1mm

Consider a $\varphi$-colored surface $\Delta_\varphi$ on the left hand side of Figure \ref{figure:recursiontr}. A Schwarzian structure on $\Delta_\varphi$ is parametrised by a triple of decorated flags $(\A_p, \A_q, \A_r)$  and  $\varphi$ along the interval $qr$. 
Recall $\lambda=\K^{-1/2}$. Let 
\[a:=\lambda_{qp}, \;\;\;b:=\lambda_{pr},\;\;\; c:=\lambda_{qr}.\] 
Then the potential of the triangle $pqr$ is
\[W_{\Delta}=\frac{a}{bc}+\frac{b}{ac}+\frac{c}{ab}.
\]

Cut the    $\Delta_\varphi$ along the edge $pa_t$, as on the right hand side of Figure \ref{figure:recursiontr}, getting a decomposition  $$
\Delta_\varphi=\Delta \cup \D^{-}_\varphi.
$$ 
\vskip 1mm
Recall the notation
\be 
\begin{split}
&{\rm I}(s):=\int_0^s e^{\xi(x)} dx. \\
& 
\varphi(s)= \frac{{\rm I}(s)}{{\rm I}(t)}.\\
\end{split}
\ee

Recall from (\ref{eq:intdiff}) that that $\varphi\in \mathrm{Diff}_+^1(c,\T)$  just means that $\varphi'(t)>0$ and continuous and 
$$
\varphi(0)=0, \ \varphi(\T)=1, \ \varphi'(0)=\varphi'(\T)=\frac{\pi}{c}.
$$

Recall  from (\ref{124}) that $\xi\in \mathrm{C}(c,\T)$  just means that $\xi$ is continuous and  
$$
\xi(0)=\xi(\T)=0, \ {\rm I}(\T) =\frac{c}{\pi}.
$$ 
Therefore if $\varphi(s) \in \mathrm{Diff}_+^1(c,\T)$ then 
$
\xi(s)=\log \frac{\varphi'(s)}{\varphi'(0)}\in \mathrm{C}(c,\T).
$ 

Denote by $[\A]$  the flag underlying a decorated flag $\A$. The function $\varphi$ encodes a map  
\be
\begin{split}
\B: [0, \T] \lra &\R\P^1\\
s \lms &\B_s=[(1-\varphi(s)) \A_q + \varphi(s) \A_r].\\
\end{split}
\ee
 
 Consider {\em canonical lift } of the family of flags $\B_s$ to a family of decorated flags  between $\A_q$ and $\A_r$:\footnote{Note that it satisfies the following condition 
$$
\lim_{\epsilon \rightarrow 0} \frac{\left<\A_s, \A_{s+\varepsilon}\right>}{\varepsilon}=\frac{c \varphi'(s)}{g(s)^2}= \frac{\pi e^{\xi(s)}}{g(s)^2}=\pi.
$$}
\be \la{CL}
\A_s:=e^{-\frac{\xi(s)}{2}}\left(  \left( 1-\frac{{\rm I}(s)}{{\rm I}(\T)} \right) \A_q +\frac{{\rm I}(s)}{{\rm I}(\T)} \A_r\right).
\ee
Set $h:=e^{\frac{\xi(t)}{2}}$, where $t$ is the given number. Let us set   then 
\be \la{193a}
\begin{split}
&c_-:=\left<\A_q, \A_t\right>= 
\frac{\pi {\rm I}(t)}{h},\\
&c_+:=\left<\A_t, \A_r\right>= \frac{\pi ({\rm I}(\T)-{\rm I}(t))}{h},\\
& c:= h(c_- + c_+) =\pi {\rm I}(\T),\\
&z:=\left<\A_p, \A_t\right>=(a c_+ + b c_-)/c.\\
\end{split}
\ee
Here  the equalities in the lines 1-2 follow from (\ref{CL}) and ${\rm I}(\T) =\frac{c}{\pi}$, the equality in the line 3 is evident, while the equality in the line 4 is the Pl\"ucker relation.

\bp \la{Prop11.3}
Given any $a,b>0$, $0<t<\T$, let $\mathcal{C}_\T=\{(c,\xi)| c>0, \xi\in {\rm C}(c,\T)\}$,
\be
\begin{split}
&\mathcal{D}_t:=\{(z,c_-,c_+,\xi_-,\xi_+)| z,c_-,c_+>0, \xi_- \in {\rm C}(c_-,t),  \xi_+ \in {\rm C}(c_+,\T-t)\}.
\end{split}
\ee
There is a measurable isomorphism of spaces
\be
\begin{split}
{\rm Cut}_t: &\mathcal{C}_\T\longrightarrow  \mathcal{D}_t.\\
\end{split}
\ee
\ep

\bpr 
Set: 
\be  
\label{eq:Usxi}
\begin{split}
\U_-(s):=\frac{{\rm I}(s)}{{\rm I}(t)}, \ \ \ \ \ \ &\varphi_-(s):= \frac{h \U_-(s)}{1+(h-1) \U_-(s)},   \ \ \ \ \ \   s\in [0,t]. \\
\U_+(s):=\frac{{\rm I}(s+t)-{\rm I}(t)}{{\rm I}(\T)-{\rm I}(t)}, \ \ \ \ \ \ & \varphi_+(s): =\frac{\U_+(s)}{h+ (1-h)\U_+(s)}, \ \ \ \ \ \ \ s\in [0,\T-t].\\
\end{split}
 \ee
 Then we  define the pair $(\xi_-, \xi_+)$ by setting:
 \be \la{191}
\begin{split}
 &\xi_-(s):=\log \frac{\varphi'_-(s)}{\varphi'_-(0)} =\xi(s)-2 \log \Bigl(1+(h-1) \U_-(s)\Bigr), \ \ \ \ \ \  \ \ \  \  \;\forall s\in [0,t].   \\
 &\xi_+(s):=\log \frac{\varphi'_+(s)}{\varphi'_+(0)}  = \xi(s+t)-2\log\Bigl(h+ (1-h)\U_+(s)\Bigr), \ \ \ \ \ \ \ \forall s\in [0,\T-t].\\  
 \end{split}\ee


\bl We have  
 \be \la{991}
\begin{split}
&\xi_-(0)=  \xi_-(t)=0, \ \ \ \ \ \  \ \ \ \ \int_0^t e^{\xi_-(x)}dx = 
\frac{c_-}{\pi}.\\
&\xi_+(0)= \xi_+(\T-t)=0, \ \ \  \ \int_0^{\T-t} e^{\xi_+(x) dx }= 
\frac{c_+}{\pi}.\\
\end{split}
\ee
\el
This just means that
\be
\begin{split}
&\xi_- \in \mathrm{C}(c_-,t), \ \ \ \xi_+ \in \mathrm{C}(c_+,\T-t).\\
&\varphi_- \in \mathrm{Diff}_+^1(c_-,t), \ \ \ \varphi_+ \in \mathrm{Diff}_+^1(c_+,t)\\
\end{split}
\ee

\bpr The left two formulas in each  line of (\ref{991}) are clear from (\ref{191}). Next
$$
\int_0^{t} e^{\xi_-(x) }dx \stackrel{(\ref{191})}{=} \varphi'_-(0)^{-1} \cdot (\varphi_-(t) - \varphi_-(0)) \stackrel{(\ref{eq:Usxi})}{=}  \varphi'_-(0)^{-1}  = \frac{c_-}{\pi}
$$
The second claim is similar.\epr

The canonical lift of the subpath from $0$ to $t$ of flags $\B_s$ between $\A_q$ and $\A_t$ is identical with $\A_s$:
\be
e^{-\frac{\xi(s)}{2}} \left( (1-\varphi(s) ) \A_q +\varphi(s)\A_r\right) = e^{-\frac{\xi_-(s)}{2}}( (1-\varphi_-(s))\A_q+ \varphi_-(s) \A_t).
\ee
The canonical lift of the subpath from $0$ to $\T-t$ of flags $\B_{s+t}$ between $\A_t$ and $\A_r$ delivers   $\A_{s+t}$:
\be
e^{-\frac{\xi(s+t)}{2}}\Bigl( (1-\varphi(s+t))\A_q+ \varphi(s+t) \A_r\Bigr)=e^{-\frac{\xi_+(s)}{2}}\Bigl( (1-\varphi_+(s))\A_t+ \varphi_+(s) \A_r\Bigr).
\ee

This construction provides the  map ${\rm Cut}_t$. 
Let us compute the inverse map ${\rm Cut}_t^{-1}(z,c_-,c_+,\xi_-,\xi_+)=(c,\xi)$. 
The second equation in the first line of  \eqref{eq:Usxi} gives
\be
\label{eq:Us}
\begin{split}
&\U_-(s)=\frac{\varphi_-(s)}{h+(1-h)\varphi_-(s)}, \ \ \ \ \U_+(s)= \frac{h\varphi_+(s)}{1+(h-1)\varphi_+(s)}. \\ 
\end{split}
\ee
Plugging this into \eqref{191} we get
\be
\label{eq:xit}
\begin{split}
\xi(s)&=\xi_-(s)-2\log\left(1+(h^{-1}-1)\varphi_-(s) \right), \ \  s\in [0,t].\\
\xi(t+s)-\xi(t)&=\xi_+(s)-2\log\left(1+(h-1)\varphi_+(s)\right), \ \ \ \ \  s \in [0,\T-t].\\
\end{split}
\ee
Proposition \ref{Prop11.3} proved. 
\epr
 
\vskip 2mm

   This map induces a canonical  map $\mathrm{Cut}_t$ on the probabilistic moduli space. 
\bt
\label{thm:cutt}
\be
\begin{split}
&(\mathrm{Cut}_t)_* \left(\int_{c=0}^\infty \mathbf{E}_{\Delta_\varphi}(\T,a^{-2},b^{-2})\right) \\
=&\frac{1}{2}\int_{z=0}^\infty \int_{c_-=0}^\infty \int_{c_+=0}^\infty \mathbf{E}_{\Delta_{\varphi_-}}(t,a^{-2},z^{-2})\otimes \mathbf{E}_{\Delta_{\varphi_+}}(\T-t,z^{-2},b^{-2}) \frac{dz}{z}.\\
\end{split}
\ee
\et
We prove Theorem \ref{thm:cutt} in Section \ref{YK}.

\bc \la{cor}
For any $0<t<T$, we have
\be
\mathcal{E}_{\Delta_\varphi}(\T,\K_{pq},\K_{pr})=\frac{1}{4}\int  \mathcal{E}_{\Delta_{\varphi_1}}(t,\K_{pq},\K_{p a_t}) \mathcal{E}_{\Delta_{\varphi_2}}(\T-t,\K_{p a_t},\K_{pr}) d \log \K_{p a_t},
\ee
where $q a_t$ is $\varphi_1$-colored and $ a_t r$ is $\varphi_2$-colored.
\ec

We define 
\be \la{184}
\widetilde{\mathcal{E}}_{\Delta_\varphi}(\T,\K_{pq},\K_{pr}):=\int_{\mathcal{M}_{\Delta_\varphi}(\T,\K_{pq},\K_{pr})} \frac{c}{ab} \mathbb{E}_{\Delta_\varphi}(\T,\K_{pq},\K_{pr}).
\ee

\bc
\be
\label{eq:dectr}
\begin{aligned}
&\widetilde{\mathcal{E}}_{\Delta_\varphi}(\T,\K_{pq},\K_{pr})= 
\\&\frac{1}{4}\int_0^\infty \left( \widetilde{\mathcal{E}}_{\Delta_\varphi}(t,\K_{pq},\K_{p a_t}) \mathcal{E}_{\Delta_\varphi}(\T-t,\K_{p a_t},\K_{pr})+\mathcal{E}_{\Delta_\varphi}(t,\K_{pq},\K_{p a_t})\widetilde{\mathcal{E}}_{\Delta_\varphi}(\T-t,\K_{p a_t},\K_{pr})\right) d \log \K_{p a_t}.
\end{aligned}
\ee
\ec
\bpr
The potentials satisfy
\be \la{PId}
\frac{c}{ab}=\frac{c_-}{az}+\frac{c_+}{zb}.
\ee

By Theorem \ref{thm:cutt}, we have
\be
\begin{aligned}
&\widetilde{\mathcal{E}}_{\Delta_\varphi}(\T,\K_{pq},\K_{pr})=\int_{\mathcal{M}_{\Delta_\varphi}(\T,\K_{pq},\K_{pr})} \left(\frac{c_-}{az}+\frac{c_+}{zb}\right) \mathbb{E}_{\Delta_\varphi}(\T,\K_{pq},\K_{pr})
\\&=\frac{1}{2}\int_{z=0}^\infty  \int_{\mathcal{M}_{\D_{\varphi_-}}(t,\K_{pq},\K_{p a_t})} \int_{\mathcal{M}_{\D_{\varphi_+}}(\T-t,\K_{p a_t},\K_{p r})} 
 \left(\frac{c_-}{az}+\frac{c_+}{zb}\right) \mathbb{E}_{\Delta_\varphi}(t,\K_{p q},\K_{p a_t}) \mathbb{E}_{\Delta_\varphi}(\T-t,\K_{p a_t},\K_{p r})  \frac{dz}{z}
\end{aligned}
\ee
Here the equality in the first line follows from the definition (\ref{184}) and  relation (\ref{PId}). The next one follows from Theorem \ref{thm:cutt}. Then we write the integrand of the last formula as
\be
\begin{aligned}
 &\frac{c_-}{az}\cdot \mathbb{E}_{\Delta_\varphi}(t,\K_{p q},\K_{p a_t}) \mathbb{E}_{\Delta_\varphi}(\T-t,\K_{p a_t},\K_{p r})
+ \frac{c_+}{zb}\cdot  \mathbb{E}_{\Delta_\varphi}(t,\K_{p q},\K_{p a_t}) \mathbb{E}_{\Delta_\varphi}(\T-t,\K_{p a_t},\K_{p r}).\\
 \end{aligned}
\ee
 After this  definition (\ref{184})  implies immediately 
formula (\ref{eq:dectr}).
\epr

\subsection{Proof of Theorem \ref{thm:cutt}.} \la{YK}

 Section \ref{YK}  relies on the concatenating formula for the unnormalised Brownian bridge measures in \cite[Proposition 15]{Lo24}  and the change of variables formula for the unnormalised Brownian bridge in \cite[Lemma 18]{Lo24}, proved in \cite[Proposition 21]{BLW24}. 
\vskip 1mm

\bl
\label{lem:changev}
For any $\xi \in  \mathrm{C}(c,\Y)$, let $h=e^{k/2}$. As in \eqref{eq:xit}, we define
\be
R_k(\xi):=\xi-2\log(1+(h^{-1}-1)\P(\xi)):=\psi.
\ee
Recall $\beta=2\pi/\sigma^2$. Then 
\be
(R_k^{-1})_* \mathcal{B}_{\sigma^2}^{0,k; \Y} = \exp\left(-\frac{\beta (h-1)^2}{h c} \right) \mathcal{B}_{\sigma^2}^{0,0; \Y}.
\ee
\el
\bpr
We  use \cite[Lemma 18]{Lo24}.
Recall $\P(\xi)(x)=\int_0^x e^{\xi(s)}ds \cdot ( \int_0^Y e^{\xi(s)}ds)^{-1}$ and $\P^{-1}(\varphi)(s)=\log\frac{\varphi'(s)}{\varphi'(0)}$. Let 
\[g(u)=\frac{h u}{1+(h-1)u},\;\;\; \ \ \ \ g^{-1}(u)=\frac{u}{h+(1-h)u}.\]
Let $\L_{f}$ be left composition. From \eqref{eq:xit}, 
\eqref{eq:Us} and \eqref{eq:Usxi} we obtain
\[R_k=\P^{-1}\circ \L_{g^{-1}} \circ \P.
\]
We have
\[g'(u)=\frac{h}{(1+(h-1)u)^2}, \;\;\; \ \ \ \ \frac{g''(u)}{g'(u)}=-\frac{2(h-1)}{1+(h-1)u}.\]
Then
\[g'(0)=h,\;g'(1)=h^{-1},\; \log g'(1)-\log g'(0)=-k,\; \frac{1}{\sqrt{g'(0)g'(1)}}=1.\]
Sicne $g$ is a M\"obius transformation, we have $\mathrm{Sch}(g)=0$. 
And we have
\[J=\int_0^\Y e^{\psi(x)}dx=\frac{h c}{\pi},\; \ \ \P(\psi)'(0)=\frac{1}{J},\; \ \ \P(\psi)'(Y)=\frac{h^2}{J}.\]
Then
\[\frac{1}{\sigma^2} \left(\frac{g''(0)}{g'(0)}\P(\psi)'(0)-\frac{g''(1)}{g'(1)}\P(\psi)'(\Y) \right)=\frac{2(h-1)^2}{\sigma^2 J}=\frac{\beta (h-1)^2}{h c}.\]
Thus by \cite[Lemma 18]{Lo24}, we get 
\be
(R_k)_* \mathcal{B}_{\sigma^2}^{0,0; \Y} = \exp\left(\frac{\beta (h-1)^2}{h c} \right) \mathcal{B}_{\sigma^2}^{0,k; \Y}.
\ee
Hence
\be
(R_k^{-1})_* \mathcal{B}_{\sigma^2}^{0,k; \Y} = \exp\left(-\frac{\beta (h-1)^2}{h c} \right) \mathcal{B}_{\sigma^2}^{0,0; \Y}.
\ee

\epr

By applying Lemma \ref{lem:changev} to \eqref{eq:btheta}, we obtain
\bc
\be
4 \mathcal{B}_{\sigma^2}^{0,l; \T}=(R_k)_* \left(2\int_0^\infty  \exp(-\beta (e^{l/2}+e^{-l/2})/c) \Theta_{\sigma^2}^{c,\T} dc/c \right).
\ee
\ec
When a boundary circle is colored green, by Proposition \ref{D9.18}, we can reduce the red points to one. By the above corollary, we can get rid of the last single red point.

\bpr[Proof of Theorem \ref{thm:cutt}]
Combining the above equation with \eqref{eq:btheta}, we get
\be
\label{eq:rbY}
(R_k^{-1})_* \mathcal{B}_{\sigma^2}^{0,k; \Y} =\int_0^\infty \frac{dc}{2c} e^{-\beta(h+h^{-1})/c} d\Theta_{\sigma^2}^{c, \Y}.
\ee

The Schwarzian measure on $\Delta_\varphi$ is
\[\mathbb{E}_{\Delta_\varphi}(\T,a^{-2},b^{-2})=2\exp(-\beta W_{\Delta})\frac{dc}{c} d\Theta_{\sigma^2}^{c,\T}.\]
Then by \eqref{eq:btheta}, we get
\[ \int_c \mathbb{E}_{\Delta_\varphi}(\T,a^{-2},b^{-2}) = 4 \exp(-\beta (W_{\Delta}-2/c) ) d \mathcal{B}_{\sigma^2}^{0,0; \T}.\]
By \cite[Proposition 15]{Lo24}, we have
\be
\mathcal{B}_{\sigma^2}^{0,0; \T}=\int_{\mathbb{R}} \rm{Concat}_*\left(\mathcal{B}_{\sigma^2}^{0,k; t}\otimes \mathcal{B}_{\sigma^2}^{k,0; \T-t}\right) dk.
\ee
Applying \eqref{eq:rbY} and the similar formula for the bridge from $k$ to $0$, we get exactly our previous paths $\xi_-$ and $\xi_+$, then the concatenated measure is
\be
\label{eq:tsqm}
\int_{k} \int_{c_-} \int_{c_+}  dk \frac{d c_-}{c_-}  \frac{d c_+}{c_+}\exp[-\beta(W_{\Delta} -2/c+(h+h^{-1})(1/c_-+1/c_+) )] d\Theta_{\sigma^2}^{c_-,t}(\xi_-)  d\Theta_{\sigma^2}^{c_+,\T-t}(\xi_+).
\ee
Since $c=h(c_-+c_+)$, we get
\[
(h+h^{-1})(1/c_-+1/c_+) 
=
\frac{c}{c_- c_+}+\frac{c_-}{c c_+}+\frac{c_+}{c_- c}+2/c.
\]
Thus since the potential of the  square equals to the sum of the potentials for  any of its triangulations,
\[W_{\Delta} +(h+h^{-1})(1/c_-+1/c_+) - 2/c = W_{pq a_t}+W_{p r a_t}.\]
And, see (\ref{193a}), 
\[h=\frac{c_- b+c_+ a}{(c_- + c_+)z}.\] 
Then, since  $dk=2 dz/z$, we get 
\[\eqref{eq:tsqm} =\frac{1}{4} \int_{z} \int_{c_-} \int_{c_+} \frac{2 dz}{z} \left(2 e^{-\beta W_{p q a_t}}\frac{d c_-}{c_-} d \Theta_{\sigma^2}^{c_-,t}\right)\left(2 e^{-\beta W_{p r a_t}}\frac{d c_+}{c_+} d \Theta_{\sigma^2}^{c_+,\T-t}\right).\]
Theorem \ref{thm:cutt} is proved.\\

\epr

\subsection{Recursions for general surfaces}

To run  a  recursion,  the surface $\bS_\varphi$ should have at least one red cusp. 
Given a ${\varphi}-$colored surface $\bS_\varphi$ and a red cusp $p$, we study the recursion for 
$
\int_{\cal M_{\bS_\varphi}(\T_\bullet,\K_{\bf r},\L_{\bf r})} W_p {\Omega}_{\bS_\varphi}(\T_\bullet,\K_{\bf r},\L_{\bf r}).
$
 \vskip 2mm

Let $p'$ and $p''$ be the red cusps adjacent to the cusp $p$. If the boundary interval $pp'$  (respectively $pp''$) carries blue boundary edges, we cut out from $\bS$ the red-based half disc $\D_{pp'}$ with the base $pp'$ (respectively the disc $\D_{pp''}$ with the base $pp''$), and use the gluing formula to express the Schwarzian measure of the  surface $\bS_\varphi$  via the ones 
 for the red-based discs and the  $\bS_\varphi- \D_{p, p'}\cup \D_{p, p''}$. 
 
 So we can assume without loss of generality that the red cusp $p$ is adjacent to red cusps. The recursion formula for $\bS_\varphi$ is written as a sum over the elementary surfaces as in Figure \ref{figure:elesur}. We could write the recursion for the more general integral $\int_{\cal M_{\bS_\varphi}(\T_\bullet, \K_{\bf r}, \L_{\bf r})} f \cdot W_p {\Omega}_{\bS_\varphi}(\T_\bullet, \K_{\bf r}, \L_{\bf r})$ as \cite[Formula (145)]{GS24} following the similar strategy. We write the recursion in the following way for simplicity. 

\begin{figure}[ht]
	\centering
	\includegraphics[scale=0.85]{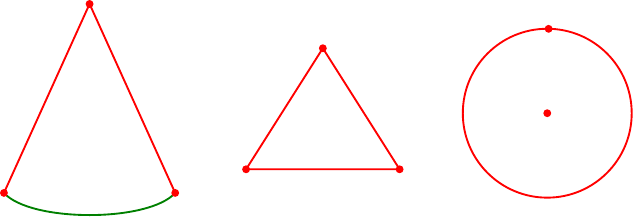}
	\small
	\caption{Cutting off elementary surfaces.}
	\label{figure:elesur}
\end{figure}

Recall from (\ref{FF}) the shorthand 
\be \la{FF1}
f(\K_1,\K_2,\K_3):= \sqrt{\frac{\K_2 \K_3}{\K_1}} \exp\left(-\sqrt{\frac{\K_1 \K_2}{\K_3}}- \sqrt{\frac{\K_2 \K_3}{\K_1}}- \sqrt{\frac{\K_3 \K_1}{\K_2}}\right).
\ee
\begin{theorem} \la{Th7.8}
Given a ${\varphi}-$colored surface $\bS_\varphi$ and a red cusp $p$ with  adjacent red edges,  there is the  recursion formula  for the weighted exponential volume:
\begin{equation}
\begin{aligned}
&\int_{\cal M_{\bS_\varphi}(\T_\bullet, \K_{\bf r}, \L_{\bf r})}  W_p {\Omega}_{\bS_\varphi}(\T_\bullet, \K_{\bf r}, \L_{\bf r})\\
=&\sum_{[\T]\in  [{\cal H}_p^\circ]}   
 \int_0^{\infty}  \int_0^{\infty}    2^{1-\mu_\T} \sqrt{\K_\T}e^{-\frac{l_\T}{2}} exp\left(-\sqrt{\K_\T}\left(e^{\frac{l_\T}{2}}+e^{-\frac{l_\T}{2}} \right)\right) \cdot\\
 &\mathcal{E}_{\bS_\varphi-\T}(\T_\bullet,\K_{\bullet},\L_{\bullet}, \K_\T,l_\T) d\log \K_\T \cdot l_\T d l_\T
\\+& \sum_{[\T]\in  [{\cal H}_p^\partial]}   
 \int_0^{\infty}     \sqrt{\K_\T}\left(e^{\frac{l_\T}{2}}+e^{-\frac{l_\T}{2}} \right) exp\left(-\sqrt{\K_\T}\left(e^{\frac{l_\T}{2}}+e^{-\frac{l_\T}{2}} \right)\right) \cdot\\
 &\mathcal{E}_{\bS_\varphi-\T}(\T_\bullet,\K_{\bullet},\L_{\bullet}, \K_\T) d\log \K_\T
\\+&\sum_{[\tau]\in  [\Delta_p^\circ], {\rm qr \;red} }   
 \int_0^{\infty}  \int_0^{\infty}     f(\K_1, \K_2, \K_3) \mathcal{E}_{\bS_\varphi-\tau}(\T_\bullet,\K_{\bf r}\backslash\{\K_1\},\L_{\bullet},\K_2,\K_3) d\log \K_2 \cdot d\log \K_3
\\+&\sum_{[\tau] \in  [\Delta_p^\partial], {\rm qr \;red} }   
 \int_0^{\infty}      f(\K_1, \K_2, \K_3)  \mathcal{E}_{\bS_\varphi-\tau}(\T_\bullet,\K_{\bullet}\backslash\{\K_1\},\L_{\bullet},\K_2) d\log \K_2
\\
+&\frac{1}{4}\sum_{[\tau]\in  [\Delta_p^\circ], {\rm qr \;green} }  \int_0^{\infty} \int_0^{\infty}  \tilde{\mathcal{E}}_{\D_{\varphi_{q a_t}}}(t,\K_{pq},\K_{p a_t})
\mathcal{E}_{{\bS-p q a_t }_{\varphi'}}(\T_\bullet\backslash \{\T_{qr}\}\cup\{\T_{qr}-t\},\L_{\bullet},\K_{\bullet}\cup\{\K_{p q},\K_{p a_t}\})  d\log \K_{p q}  d\log \K_{p a_t}
\\
+&\frac{1}{4}\sum_{[\tau]\in  [\Delta_p^\circ], {\rm qr \;green} }  \int_0^{\infty} \int_0^{\infty}  \tilde{\mathcal{E}}_{\D_{\varphi_{a_t r}}}(\T_{qr}-t,\K_{p a_t},\K_{pr})
\mathcal{E}_{{\bS-p r a_t }_{\varphi''}}(\T_\bullet\backslash \{\T_{qr}\}\cup\{t\},\L_{\bullet},\K_{\bullet}\cup\{\K_{p r},\K_{p a_t}\})  d\log \K_{p r}  d\log \K_{p a_t}
\\
+&\frac{1}{4}\sum_{[\tau]\in [\Delta_p^\partial], {\rm qr \;green} }  \int_0^{\infty} \int_0^{\infty}  \tilde{\mathcal{E}}_{\D_{\varphi_{q a_t}}}(t,\K_{pq},\K_{p a_t})
\mathcal{E}_{{\bS-p q a_t }_{\varphi'}}(\T_\bullet\backslash \{\T_{qr}\}\cup\{\T_{qr}-t\},\L_{\bullet},\K_{\bullet}\cup\{\K_{p q},\K_{p a_t}\})   d\log \K_{p q} d\log \K_{p a_t}
\\
+&\frac{1}{4}\sum_{[\tau]\in [\Delta_p^\partial], {\rm qr \;green} }   \int_0^{\infty}  \tilde{\mathcal{E}}_{\D_{\varphi_{a_t r}}}(\T_{qr}-t,\K_{p a_t},\K_{pr})
\mathcal{E}_{{\bS-p r a_t }_{\varphi''}}(\T_\bullet\backslash \{\T_{qr}\}\cup\{t\},\L_{\bullet},\K_{\bullet}\backslash\{\K_{p r}\}\cup\{\K_{p a_t}\})   d\log \K_{p a_t}
.
  \end{aligned}
\end{equation}
\end{theorem}

\bpr 
Firstly, let us cut off the red-based half-discs and red-based cylinders to obtain a decorated surface $\bS_\bullet:= \bS_\varphi- \bigcup_{e \in {\cal I}_{\varphi}}\D_{e} \ \cup \ \bigcup_{\gamma \in {\cal J}_{\varphi}} \D^*_{\gamma}$ as in Figure \ref{figure:cutlimitX}. The part $\bigcup_{e \in {\cal I}_{\varphi}}\D_{e} \ \cup \ \bigcup_{\gamma \in {\cal J}_{\varphi}} \D^*_{\gamma}$ is ${\rm Mod}(\bS_\varphi)$ invariant, and ${\rm Mod}(\bS_\varphi)={\rm Mod}(\bS_\bullet)$.
By \cite[Theorem 5.3]{GS24}, restated in Theorem \ref{thm:re}, we have the recursion formula for $\int_{\cal M_{\bS_\bullet}(\T_\bullet,\K_\bullet,\L_\bullet)} W_p \mathbb{E}^{\beta}_{\bS_\bullet}(\T_\bullet,\K_\bullet,\L_\bullet)$.
Then we glue $\bigcup_{e \in {\cal I}_{\varphi}}\D_{e} \ \cup \ \bigcup_{\gamma \in {\cal J}_{\varphi}} \D^*_{\gamma}$ back to the two side of \eqref{eq:re} along the red boundary edges and circles, we obtain a similar recursion formula for $\int_{\cal M_{\bS_\varphi}(\T_\bullet, \K_{\bf r}, \L_{\bf r})} W_p {\Omega}_{\bS_\varphi}(\T_\bullet, \K_{\bf r}, \L_{\bf r})$.
 Furthermore, we plug \eqref{eq:dectr} into the ideal triangles.
Let us elaborate this recursion more clearly.

1. The  two sums provided  by cutting off trouser legs work the same way as in the first and second sums of right hand side of Equation \eqref{eq:re}. When the boundary circle $\gamma$ of $\tau$ is the boundary of the surface and $\varphi$-colored, $\bS_\varphi-\tau$ is a disjoint union of $\bS-\tau$ and a disc $\D^*_\gamma$ with a red-based hole both inheriting the coloring from $\varphi$.

2. The sum provided by cutting off red triangles  works the same way as in the third and fourth sums of right hand side of Equation \eqref{eq:re}.

3. When cutting off triangle $\Delta=pqr \in \Delta_p^\circ$ where $qr$ is blue, we firstly cut off this $\varphi$-colored triangle $\Delta_{\varphi_{qr}}=\Delta \cup \D^{-}_{\varphi_{qr}}$ 
\[\int_0^{\infty}  \int_0^{\infty}     \tilde{\mathcal{E}}_{\Delta_{\varphi_{qr}}}(\T_{qr},\K_{qp},\K_{pr}) \mathcal{E}_{\bS_\varphi-\Delta_{\varphi_{qr}}}(\T_\bullet\backslash \{\T_{qr}\},\L_{\bullet},\K_{\bullet}\cup\{\K_{qp},\K_{pr}\}) d\log \K_{qp} \cdot d\log \K_{pr}.\]
Then by plugging Equation \eqref{eq:dectr} into the above formula, we obtain
\be
\begin{aligned}
\label{eq:revai3}
&\frac{1}{4}\int_0^{\infty} \int_0^{\infty} \tilde{\mathcal{E}}_{\Delta_{\varphi_{q a_t}}}(t,\K_{pq},\K_{p a_t})
\mathcal{E}_{{\bS-p q a_t }_{\varphi'}}(\T_\bullet\backslash \{\T_{qr}\}\cup\{\T_{qr}-t\},\L_{\bullet},\K_{\bullet}\cup\{\K_{p q},\K_{p a_t}\})  d\log \K_{p q}  d\log \K_{p a_t}  +
\\&\frac{1}{4}\int_0^{\infty} \int_0^{\infty} \tilde{\mathcal{E}}_{\Delta_{\varphi_{a_t r}}}(\T_{qr}-t,\K_{p a_t},\K_{pr})
\mathcal{E}_{{\bS-p r a_t }_{\varphi''}}(\T_\bullet\backslash \{\T_{qr}\}\cup\{t\},\L_{\bullet},\K_{\bullet}\cup\{\K_{p r},\K_{p a_t}\})  d\log \K_{p r}  d\log \K_{p a_t}
\end{aligned}
\ee
\begin{figure}[ht]
	\centering
	\includegraphics[scale=0.85]{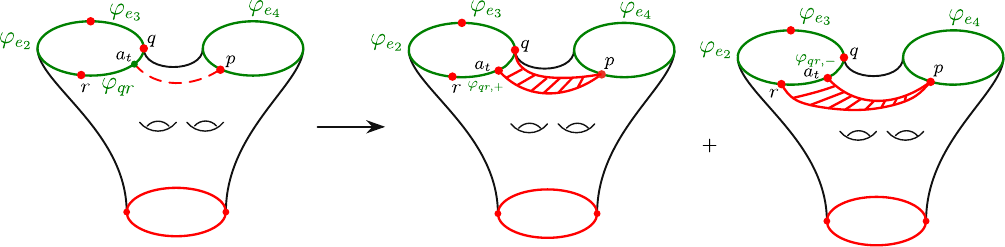}
	\small
	\caption{Cutting off an edge $e_t$.}
	\label{figure:recursionlimit}
\end{figure}
where $\varphi'=\varphi\backslash \{\varphi_{qr}\}\cup \{\varphi_{a_t r}\}$ and $\varphi''=\varphi\backslash \{\varphi_{qr}\}\cup \{\varphi_{q a_t}\}$ as in Figure \ref{figure:recursionlimit}. Similarly, for $\Delta \in \Delta_p^\partial$, $pr\in \partial \bS$, we have
\be
\begin{aligned}
\label{eq:revai4}
&\frac{1}{4} \int_0^{\infty} \int_0^{\infty} \tilde{\mathcal{E}}_{\Delta_{\varphi_{q a_t}}}(t,\K_{pq},\K_{p a_t})
\mathcal{E}_{{\bS-p q a_t }_{\varphi'}}(\T_\bullet\backslash \{\T_{qr}\}\cup\{\T_{qr}-t\},\L_{\bullet},\K_{\bullet}\cup\{\K_{p q},\K_{p a_t}\})   d\log \K_{p q} d\log \K_{p a_t}  +
\\&\frac{1}{4}\int_0^{\infty} \tilde{\mathcal{E}}_{\Delta_{\varphi_{a_t r}}}(\T_{qr}-t,\K_{p a_t},\K_{pr})
\mathcal{E}_{{\bS-p r a_t }_{\varphi''}}(\T_\bullet\backslash \{\T_{qr}\}\cup\{t\},\L_{\bullet},\K_{\bullet}\backslash\{\K_{p r}\}\cup\{\K_{p a_t}\})   d\log \K_{p a_t}.
\end{aligned}
\ee
\epr
\paragraph{Example.} 
\begin{figure}[ht]
	\centering
	\includegraphics[scale=0.65]{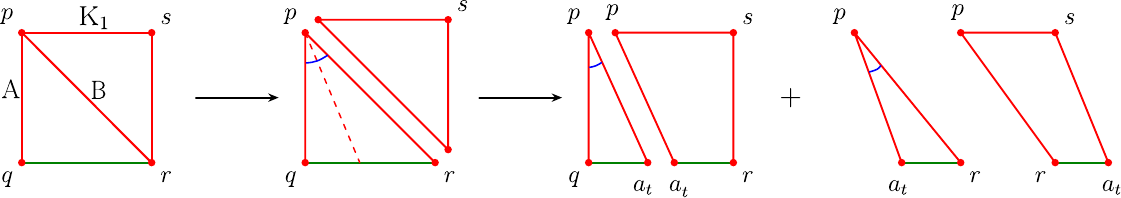}
	\small
	\caption{Recursion for $\A_{1,1,\varphi}$.}
	\label{figure:A11rec}
\end{figure}
Recall the recursion for $\A_{1,1}$ in \cite[Section 7.1]{GS24}. Let us consider the case $\A_{1,1,\varphi_{qr}}$ with parameters for an ideal triangulation as in Figure \ref{figure:A11rec}. The mapping class group $\rm{Mod}(\A_{1,1,\varphi_{qr}})=\mathbb{Z}$. By applying Theorem \ref{Th7.8}, we obtain
\be
\begin{aligned}
&\int_{{\cal M}_{\A_{1,1,\varphi_{qr}}}(\T,\K_1)}  \K_1^{1/2}(e^{l/2}-e^{-l/2}) {\Omega}_{\A_{1,1,\varphi_{qr}}}(\T,\K_1)=\int_0^{\infty}\int_0^{\infty} \tilde{\mathcal{E}}_{\Delta_{\varphi_{q r}}}(\T,\A,\B) \mathcal{E}_{\Delta}(\K_1,\A,\B) d\log \A d\log \B
\\&=\frac{1}{4} \int_0^{\infty} \int_0^{\infty} \tilde{\mathcal{E}}_{\Delta_{\varphi_{q a_t}}}(t,\A,\K_{p a_t})
\mathcal{E}_{{\A_{1,1,\varphi_{qr}}-p q a_t }_{\varphi_{a_t r}}}(\T-t,\K_1,\K_{p a_t},\A)   d\log \A d\log \K_{p a_t}  +
\\&\frac{1}{4}\int_0^{\infty} \int_0^{\infty} \tilde{\mathcal{E}}_{\Delta_{\varphi_{a_t r}}}(\T-t,\K_{p a_t},\B)
\mathcal{E}_{{\A_{1,1,\varphi_{qr}}-p r a_t }_{\varphi_{q a_t}}}(t,\K_1,\B,\K_{p a_t})  d\log \B  d\log \K_{p a_t}.
\end{aligned}
\ee

\section{Discretization of  measures on probabilistic moduli spaces} \la{Sect10}

The weighted exponential volume forms assigned to a weighted colored surface in Section \ref{SEC4} can be viewed as a  discretization of the  measure ${\bf E}$ for decorated surfaces.  In Section \ref{Sect10} we  state this precisely.

\subsubsection{Continuous limit for a half-disc: from  ${\rm Conf}(\A_p, \B_1,...,\B_\N, \A_q)$ to   $\rm{C}_0[0,\T]$.}
\label{section:discont}

Recall Definition \ref{defn:weightDh}. From now on, we assume that  all the weights are equal to $\frac{1}{2}$. Then, for the red-based colored half-disc, the weight factor is 
\be
\label{eq:qredb}
{\cal Q}_{\D^-_{\N+2,\rm{c}}}(\omega)=\frac{\K_{p1}^{\frac{1}{2}} \K_{\N q}^{\frac{1}{2}} \prod_{i=1}^{\N-1} \K_{i,i+1}^{\frac{1}{2}}}{\prod_{i=1}^\N W_i}, 
\ee
and the weighted exponential volume form  is 
\be \la{223}
{\Bbb E}_{\D^-_{\N+2,\rm{c}}}^{(\kappa,\beta)}(\omega):=
 \left(\frac{\K_{p1}^{\frac{1}{2}} \K_{\N q}^{\frac{1}{2}} \prod_{i=1}^{N-1} \K_{i,i+1}^{\frac{1}{2}}}{\prod_{i=1}^\N W_i}\right)^\kappa e^{-\beta(W_p+ W_q)}.
\ee
\vskip 1mm

 We start with a configuration  $(\A_p, \B_1,...,\B_\N, \A_q)$, and run $\N \to \infty$, as illustrated on Figure \ref{figure:interval2}. 
 
\begin{figure}[ht]
	\centering
	\includegraphics[scale=0.65]{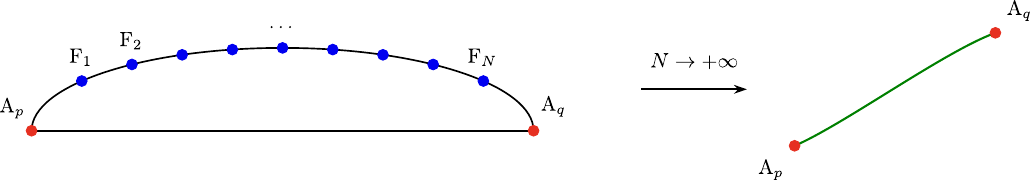}
	\small
	\caption{The weighted colored red-based half-disc $\D^-_{\N+2,c,\omega}$ converges to an interval.}
	\label{figure:interval2}
\end{figure}

Recall the continuous family of decorated flags $[(1-\varphi(t))\A_p+\varphi(t) \A_q]$. Recall $v_i\in \F_i$, $i=1, ..., \N$ and $\lambda=\omega_2(\A_p,\A_q)$ from Section \ref{section:cpAF}. Let $A_p:=v_0$. Set
\be
\label{eq:xs}
\begin{split}x_i:=&\frac{\omega_2(A_p,v_i)}{\omega_2(A_p,v_i)+\omega_2(v_i, A_q)} = (1+ \lambda W(\A_p, v_i, \A_q))^{-1}.\\
\end{split}
\ee
Then   $
x_0=0<x_1<...<x_\N<x_{\N+1}=1. 
$ 
Set $h_\N:=\frac{\T}{\N}$. Let $\varphi$ be a ${\rm C}^1$ monotone function such that 
$ \varphi(0)=0,\varphi(\T)=1$ and $\varphi(i h_\N)=x_i$. Let 
$
d_i:=x_{i+1}-x_i>0$, where  $i=0,...,\N.
$ 
 Then $\frac{d_i}{h_\N}$ is the discrete analogue of $\varphi'$ and 
$
\sum_{i=0}^\N d_i=1. 
$

Suppose that $\xi\in {\rm C}(\lambda, \T)$, see (\ref{124}), that is $\xi$ lies  in the support of the pre-Schwarzian measure $\Theta^{\lambda, \T}_{\sigma^2}$. 
The Malliavin--Shavgulidze map (\ref{MSM}) provides isomorphism (\ref{128}). 

Consider the discrete analogue of $\xi$, that is set 
$
\xi_i:= \log \frac{d_i}{d_0}, \  \forall i=0,...,\N.
$ 
The  condition $\xi(0)=\xi(\T)=0$  in  \eqref{124} implies $d_0=d_\N$. Since $\sum_{i=0}^\N d_i=1$, the discrete analogue of $\int_0^\T e^{\xi(t)} dt$ is
\be
\label{eq:intxi}
h_\N \sum_{i=0}^\N e^{\xi_i}=h_\N \frac{d_0+...+d_\N}{d_0}=\frac{h_\N}{d_0}.
\ee
Thus we get the  boundary conditions $d_0=d_\N=\frac{\pi h_\N}{\lambda}$, matching $\varphi'(\T)=\varphi'(0)=\frac{\pi}{\lambda}$ from (\ref{eq:intdiff}).

\bl
The cluster volume form for $\D_{\N+2,\rm{c}}$ is 
\be
\Omega_{\D_{\N+2,\rm{c}}}=2 d \xi_1 \wedge ...\wedge d\xi_\N \wedge d \log \K_{pq}.
\ee
\el
\bpr

Let $y_i:=\frac{1}{x_i}-\frac{1}{x_{i+1}}$. Note that $x_{\N+1}=1$. Then by  \eqref{eq:xs}, we get 
$
\B_i=\lambda y_i=\K_{pq}^{-\frac{1}{2}}y_i.
$ 
So 
\be
\begin{aligned}
&\bigwedge_{i=1}^\N d\log y_i = \bigwedge_{i=1}^\N d \left(\frac{1}{x_i}\right)/(y_1...y_\N)
=(-1)^\N \frac{\bigwedge_{i=1}^\N d x_i}{\prod_{i=1}^\N x_i^2}/\left(\prod_{i=1}^\N  \frac{d_i}{x_i x_{i+1}}\right)
\\&=(-1)^\N \frac{\bigwedge_{i=1}^\N d x_i}{d_0 d_1...d_\N }
= \frac{\bigwedge_{i=1}^\N d d_i}{d_0 d_1...d_\N }
= \bigwedge_{i=1}^\N d \xi_i,
\end{aligned}
\ee
where the last equality follows $\det \left(\frac{\partial d_i}{\partial \xi_j} \right)= d_0 ... d_\N$.
Thus 
\be
\Omega_{\D_{\N+2,\rm{c}}}=2 \bigwedge_{i=1}^\N d\log \B_i \wedge d\log \K_{pq}=2 \bigwedge_{i=1}^\N d\log y_i \wedge d\log \K_{pq}=2 d \xi_1 \wedge ...\wedge d\xi_\N \wedge d \log \K_{pq}.
\ee
\epr

Proposition \ref{lem:IW} below (for the  $\T=1$ case)  is due to by Budd and Chekhov, cf  \cite{BC}, 
with pushward by Malliavin--Shavgulidze map $\P$. It starts with the following definition. 

\bd
\label{defn:rhok}
Let
\be \la{229}
\begin{split}
&B_\kappa(0):=\int \left(2\cosh \frac{t}{2}\right)^{-\kappa} dt=\frac{\Gamma(\kappa/2)^2}{\Gamma(\kappa)}.\\
&\rho_\kappa(u):=\frac{1}{B_\kappa(0)} \left(2\cosh \frac{u}{2}\right)^{-\kappa},\\\end{split}
\ee
Let $\Xi_{\kappa, \N}$ be the law on $\mathbb{R}^\N$ under which the increments $\Y_{\N,i,\kappa}:=\xi_i-\xi_{i-1}$ for $i=1,...,\N$ are independent with density $\rho_\kappa$.
Precisely 
\be
\Xi_{\kappa, \N}:=\prod_{i=1}^\N \rho_\kappa(\xi_i-\xi_{i-1}) d \xi_1 ... d\xi_\N= B_\kappa(0)^{- \N} \prod_{i=1}^\N \left(2\cosh \frac{\xi_i-\xi_{i-1}}{2}\right)^{-\kappa} d \xi_1 ...d \xi_\N.
\ee
\ed

\bp {\rm \cite{BC}}
\label{lem:IW}
Let $\kappa_\N:=\frac{4\N}{\sigma^2 \T}$. Let $\iota_\N(f)$ be the piecewise affine interpolation with $\iota_\N(f)(i h_\N) = f_i$.
Then  in the weak topology:
\be \la{CWT}
(\iota_\N)_* \Xi_{\kappa_\N, \N} \xrightarrow{\N\rightarrow +\infty} \mathbb{W}^\T_{\sigma^2}.
\ee
\ep

Similar to Proposition \ref{lem:IW}, 
and  using the weight factor for red-based half disc in \eqref{eq:qredb}, we have 
\bp
\label{thm:Nlimit}
Let $\kappa_\N:=\frac{4\N}{\sigma^2 \T}$, $\beta:=\frac{2\pi}{\sigma^2}$.
Define
\be
R_{\N,\T}:=\frac{\kappa_\N}{4\pi} \left( \frac{e\pi \T}{\N}\right)^{\kappa_\N} B_{\kappa_\N}(0)^{-\N}.
\ee
Then the following finite Borel measures on $C_{0}[0,\T]$ converge in the weak topology:
\be
\begin{split}
&(\iota_\N)_* \left(R_{\N,\T} \cdot {\Bbb E}_{D_{\N+2,\rm{c}}}^{\kappa_\N,\beta}(\omega;\K_{pq})\right) \xrightarrow{\N\rightarrow +\infty} \Theta^{\lambda, \T}_{\sigma^2}\\
&(\P\iota_\N)_* \left(R_{\N,\T} \cdot {\Bbb E}_{D_{\N+2,\rm{c}}}^{\kappa_\N,\beta}(\omega;\K_{pq})\right) \xrightarrow{\N\rightarrow +\infty} \Omega^{\K, \T}_{\sigma^2}\\
\end{split}
\ee
\ep

\subsubsection{Continuous limit for blue $\D^*_{\N,\rm{c}}$}
Let $\D^*_{\N,\rm{c}}$ be a once punctured disc with $\N$ blue marked points on the boundary.
Given an ideal triangulation of $\D^*_{\N}$ as in Figure \ref{figure:D13}, we use the notations as in Proposition \ref{Meij}. The cluster Poisson coordinates are $\X_1$,...,$\X_\N$, and 
$e^l:=\prod_{i=1}^\N \X_i$, 
where $l\geq 0$ is the length parameter along the puncture. Set $x_i:= \log \X_i$. Then
$x_1+...+x_\N=l.$ 
The cluster volume form on this hyperplan is $2 d x_1...d x_{\N-1}$.

Let us choose a lift of the crown monodromy to the real line so that the attracting fixed point of the neck geodesic goes to infinity and the lifts $y_0,...,y_\N$ of the cusps in a fundamental domain satisfy
\[y_0=1, \;\;\;y_\N=e^l,\;\;\; y_{i+1}>y_i>0\;\;\; \forall i=1,...,\N-1.\]
Moreover $y_{\N+i}=e^ly_{i}$ extends for all $i\in \mathbb{Z}$. Let $z_i:=y_i-y_{i-1}$ for $i=1,...,\N$. Then
$x_i= \log \frac{z_{i+1}}{z_i}.$ 
Let $u_i:=\sum_{i=1}^i x_i$ and $u_0:=0$. Then 
$z_{i+1}=z_1 e^{u_i}.$ 
Since $\sum_{i=1}^\N z_i=e^l-1$, by plugging the above formula into this one, we obtain 
$z_1=\frac{e^l-1}{\sum_{i=0}^{\N-1} e^{u_i}}.$ 
Therefore
\be
\label{eq:Hl}
\log y_i=\log (y_0+\sum_{i=1}^i z_i )=\log (1+z_1\sum_{j=0}^{i-1} e^{u_j})=\log \left( 1+(e^l-1)\frac{\sum_{j=0}^{i-1} e^{u_j}}{\sum_{j=0}^{\N-1} e^{u_j}} \right).
\ee

A function $\xi \in \rm C_{0,l}[0,\T]$ such that $\xi(i \T/\N)=u_i$ is a continuous analogue of $u_i$. 
The continuous analogue of  \eqref{eq:Hl} is the function $\mathcal{R}_l(\xi)(t)$ in  (\ref{161}), which lies in the space $\mathrm{Diff}_{\rm qp}^1(\R; l, \T)_+$, 
 see (\ref{160}).    
\vskip 1mm

For the moduli space ${\cal M}_{\D^*_{\N}}$, the  $\B$- and $\K$-parameters are $\B_1$,...,$\B_\N$, $\K_1$,...,$\K_\N$.
By (\ref{eq:crownkb}), we obtain
\be
{\cal Q}_{\D^*_{\N,\rm{c}}}(\omega)=\frac{\prod_{i=1}^\N \K_i^{\frac{1}{2}}}{\prod_{i=1}^\N W_i}=\frac{\prod_{i=1}^\N(\B_i \X_i \B_{i+1})^{\frac{1}{2}}}{\prod_{i=1}^\N (\B_i+\B_i \X_i)}=\frac{1}{\prod_{i=1}^\N \left(2\cosh\left(x_i/2\right)\right)}.
\ee 
Fixing $l$, the $(\kappa,\beta)$-weighted exponential volume form 
\[{\Bbb E}^{(\kappa,\beta)}_{\D^*_{\N,\rm{c}}}(\omega;l)=\frac{2}{\prod_{i=1}^\N \left(2\cosh\left(x_i/2\right)\right)^{\kappa}} dx_1...dx_{\N-1}.\]
Set 
\be
\Xi_{\kappa, \N}(l):=B_\kappa(0)^{- \N} {\Bbb E}^{(\kappa,\beta)}_{\D^*_{\N,\rm{c}}}(\omega;l)=\left(\prod_{i=1}^\N \rho_\kappa(x_i)\right) d x_1 ... d x_{\N-1}.
\ee
So using  Proposition \ref{lem:IW} and \eqref{eq:WB}, we get
\bc
Let $\kappa_\N:=\frac{4\N}{\sigma^2 \T}$. For any  $l\in \mathbb{R}$, suppose $\sum_{i=1}^N x_i=l$.
Then, in the week topology,
\be
\begin{split}
&(\iota_\N)_* \Xi_{\kappa_\N, \N}(l) \xrightarrow{\N\rightarrow +\infty} \mathcal{B}^{l,\T}_{\sigma^2};\\
&(\mathcal{R}_l\iota_\N)_* \Xi_{\kappa_\N, \N}(l) \xrightarrow{\N\rightarrow +\infty} {\mathcal{R}_l}_*\mathcal{B}^{l,\T}_{\sigma^2}.\\
\end{split}
\ee

\ec

\begin{figure}[ht]
	\centering
	\includegraphics[scale=0.7]{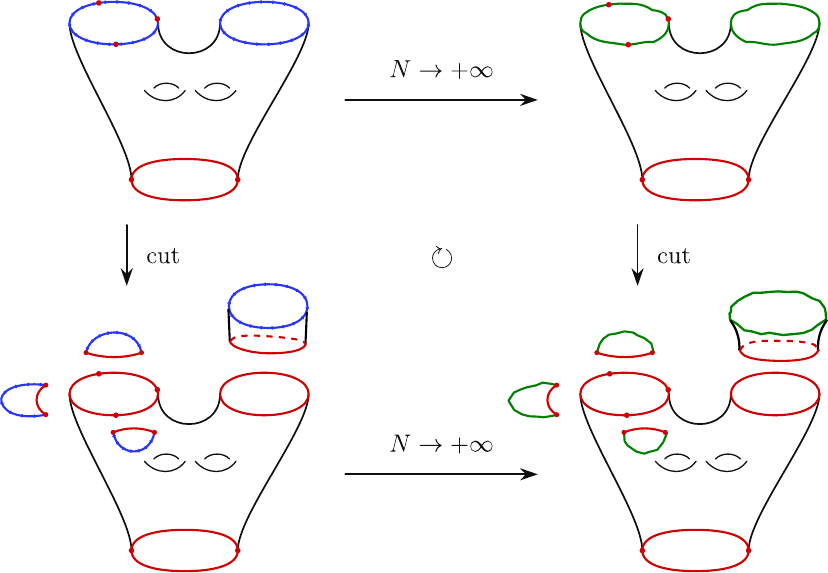}
	\small
	\caption{Cutting and gluing in large $\N$ limit.}
	\label{figure:cutlimit}
\end{figure}

\subsubsection{The general case} 

For each $\varphi-$colored  edge $e$, let $\D_e:=\D_{\N+2,\rm{c}_e}$ be its discrete analogue. Similarly, for each boundary circle $S^1_\gamma$ let $\D^*_\gamma:=\D^*_{\N,\rm{c}_\gamma}$ be its discrete analogue with $\N$ blue points. All weights $\omega_e$, $\omega_\gamma$ are $1/2$.
\vskip 1mm

Recall ${\rm R}_{\N,\T_e}$ and $\mathbb{E}_{\D_e}^{\frac{4\N}{\sigma^2 \T_e},\beta}(\omega_e;\K_e)$  in Theorem \ref{thm:Nlimit}. Consider the weighted colored surface
\be
(\bS_\N, \rm{c}; \omega)=\cup_{e} (\D_e)  \cup_\gamma (\D_\gamma^*) \cup \bS_\bullet
\ee
illustrated in left hand side of Figure \ref{figure:cutlimit}, and define the induced weighted exponential volume form 
\be
\begin{split}
&\mathbb{E}_{\bS,\mathcal{I},\N}(\T_\bullet, \K_{\bf r}, \L_{\bf r}) 
:=   {\Bbb E}^\beta_{\bS}(\K_\bullet, \L_\bullet) \wedge \\
& \bigwedge_{e\in \mathcal{I}_\varphi} (\P\iota_\N)_*\left( R_{\N,\T_e}\cdot \mathbb{E}_{\D_e}^{\frac{4\N}{\sigma^2 \T_e},\beta}(\omega_e;\K_e)   \wedge d \log \K_e \right) \bigwedge_{\gamma\in \mathcal{J}_\varphi} (\mathcal{R}_l\iota_\N)_*\left({\Bbb E}^{\frac{4\N}{\sigma^2 \T_\gamma},\beta}_{\D^*_\gamma}(\omega_\gamma;l_\gamma) \wedge d\theta_\gamma \wedge dl_\gamma\right).\\
\end{split}
\ee
Then by the gluing property and Theorem \ref{thm:Nlimit}, we have
\bt
The diagram in Figure \ref{figure:cutlimit} commutes and we have in the weak topology:
\be
\left(\prod_{e_\in \mathcal{I}_\varphi} R_{\N,\T_e}\right) \left(\prod_{\gamma \in \mathcal{J}_\varphi} \frac{1}{2 B_\kappa(0)^{\N }}\right) \mathbb{E}_{\bS,\mathcal{I},\N}(\T_\bullet, \K_{\bf r}, \L_{\bf r})
 \xrightarrow{\N\rightarrow +\infty} 
{\bf E}_{\bS_\varphi}(\T_\bullet, \K_{\bf r}, \L_{\bf r}).
\ee
\et

\section{Appendix: Multiple $\Gamma-$functions} \la{Section7}

\bd The multiple $\Gamma-$function is given by 
\be
\gamma(s_1, ..., s_m)=2^{s_1+\ldots + s_m-1}\cdot 
\int_{ t_1 \geq \ldots \geq t_m >0} e^{-(t_1+ \ldots + t_m)} t_1^{s_1}  \ldots  t_m^{s_m} d \log t_1 \ldots  d \log t_m.
\ee
\ed
One has the symmetrization formula
\be
2^{s_1+\ldots + s_m-1}\cdot \sum_{\sigma \in S_m}\gamma(s_{\sigma(1)}, ..., s_{\sigma(m)})= \Gamma(s_1) \ldots \Gamma(s_m).
\ee


\bp  The function $\gamma(a,b)$  has the following properties:

\begin{enumerate}

\item 
Symmetrization:
\be
\gamma(a,b) + \gamma(b,a) = 2^{a+b-1}\Gamma(a) \Gamma(b).
\ee

\item  Recursion formulas:
\be  \la{RF1}
\gamma(a,b+1)=2b \gamma(a,b) - \Gamma(a+b).
\ee
\be \la{RF2}
\gamma(a+1,b)=2a \gamma(a,b) + \Gamma(a+b).
\ee

\item  Specialization formulas: 
\be
\begin{aligned}
&\gamma(a,1)=(2^a-1)\Gamma(a).\\
&\gamma(1,b)= \Gamma(b).\\
\end{aligned}
\ee

\item The Barnes integral:
\be
\label{eq:d1compa}
\begin{aligned}
& \gamma(a,b) = \frac{2^{2(a+b-1)}}{\pi i} \int_{-i \infty}^{+i \infty} \frac{\Gamma(x)\Gamma(a+b-x)\Gamma(b-x)}{\Gamma(b+1-x)} d x.
\end{aligned}
\ee

\item The residues. 
\[{\rm Res}_{v=-n} \gamma(u-v,v)= 2^{u-1} \mathrm{Res}_{v=-n} \Gamma(u-v)  \Gamma(v),  \ \ \ \ \ n \in \Z_{\geq 0}.\] 
\end{enumerate}
\ep

\bpr  i) Clear. 
ii) Integrating  by parts for the  $\B$ variable, we get:
\be
\begin{aligned}
&2^{-(a+b-1)}\cdot 
\gamma(a,b)=\int_{ \A \geq \B \geq 0} e^{-(\A+\B)} \A^{a} \B^{b-1} d \log \A d \B
\\&= \int_0^{+\infty} e^{-\A} \A^{a} \left[ e^{-\B} \frac{\B^b}{b} \right]_0^{\A}  d\log \A
+\int_0^{+\infty} e^{-\A} \A^{a} \int_0^{\A} e^{-\B}\frac{\B^b}{b} d \B d\log \A
\\&= \int_0^{+\infty} e^{-2\A} \frac{\A^{a+b}}{b} d\log\A+\frac{\gamma(a,b+1)}{2^{a+b}b}
\\&=  \frac{\Gamma(a+b)}{2^{a+b}b}+\frac{\gamma(a,b+1)}{2^{a+b}b}.
\end{aligned}
\ee
So we get recursion formula (\ref{RF1}).
 
Integrating by parts along the  $\A$ variable, we get:
\be
\begin{aligned}
&2^{-(a+b-1)}\cdot\gamma(a,b)= \int_{\A\geq \B \geq 0} e^{-(\A+\B)} \A^{a-1} \B^{b-1} d \A d \B
\\&= \int_0^{+\infty} e^{-\B} \B^{b-1} \left[ e^{-\A} \frac{\A^a}{a} \right]_B^{+\infty} d\B
+\int_0^{+\infty} e^{-\B} \B^{b-1} \int_B^{+\infty} e^{-\A} \frac{\A^a}{a} d \A d\B
\\&= -\int_0^{+\infty} e^{-2\B} \frac{\B^{a+b-1}}{a} d\B+\frac{\gamma(a+1,b)}{2^{a+b}a}
\\&=  -\frac{\Gamma(a+b)}{2^{a+b}a}+\frac{\gamma(a+1,b)}{2^{a+b}a}.
\end{aligned}
\ee
So we get the recursion formula (\ref{RF2}).
\vskip 1mm

iii) One has 
\be
\begin{aligned}
&\gamma(a,1)=2^{a}\cdot \int_{0}^{+\infty} \int_{0}^A e^{-\B} d \B  e^{-\A} \A^{a} d \log \A 
\\&= 2^{a}\cdot \int_0^{+\infty} e^{-\A} (1-e^{-\A}) A^{a}  d\log \A
=(2^a-1)\Gamma(a).
\end{aligned}
\ee
Similarly, we have 
\be
\begin{aligned}
&\gamma(1,b)=2^b\cdot \int_{0}^{+\infty} \int_{B}^{+\infty} e^{-\A} d \A  e^{-\B} \B^{b} d \log \B 
\\&= 2^b\cdot \int_0^{+\infty} e^{-2\B} B^{b}  d\log \B
=\Gamma(b).
\end{aligned}
\ee

iv) Compare formulas (\ref{eq:d1comp}) and (\ref{gammaba}). 

v)
 Then the residue
\be
\begin{aligned}
&\mathrm{Res}_{v=-n} \gamma(u-v, v)
\\&=2^{u-1} \Gamma(u) \Gamma(u+n)/\Gamma(u) \frac{(-1)^n}{n!}
\\&=2^{u-1} \mathrm{Res}_{v=-n} \Gamma(u) \frac{\Gamma(u-v)}{\Gamma(u)} \Gamma(v)
\\&=2^{u-1} \mathrm{Res}_{v=-n} \Gamma(u-v) \Gamma(v).
\end{aligned}
\ee
\epr

 \end{document}